\documentclass[12pt]{amsart}
\newcommand{\deleted}[1]{}
\newcommand{\delete}[1]{}
\newcommand{\mynotes}[1]{}
\newcommand\notes[1]{}
\usepackage{amsfonts}
\usepackage{stmaryrd}
\usepackage{}
\usepackage{amsfonts}      
\usepackage{txfonts}
\usepackage{amssymb}
\usepackage{eucal}
\usepackage{graphicx}
\usepackage{amsmath}
\usepackage{amscd}
\usepackage[all]{xy}           
\usepackage[active]{srcltx} 
\usepackage{booktabs}
\usepackage{amsfonts,latexsym}
\usepackage{xspace}
\usepackage{epsfig}
\usepackage{float}
\usepackage{color}
\usepackage{fancybox}
\usepackage{colordvi}
\usepackage{multicol}
\usepackage{colordvi}
\usepackage{longtable}
\usepackage{stmaryrd}
\usepackage[colorlinks,final,backref=page,hyperindex,hypertex]{hyperref}

\newcommand\changed[1]{#1}
\changed{}

\newtheorem{theorem}{Theorem}[section]
\theoremstyle{definition}

\deleted{\newtheorem{theorem}{Theorem}[section]
\newtheorem{prop-def}{Proposition-Definition}[section]
\newtheorem{coro-def}{Corollary-Definition}[section]

}

\newcommand{\nc}{\newcommand}
\nc{\tred}[1]{\textcolor{red}{#1}} \nc{\tblue}[1]{\textcolor{blue}{#1}} \nc{\tgreen}[1]{\textcolor{green}{#1}} \nc{\tpurple}[1]{\textcolor{purple}{#1}} \nc{\btred}[1]{\textcolor{red}{\bf #1}} \nc{\btblue}[1]{\textcolor{blue}{\bf #1}} \nc{\btgreen}[1]{\textcolor{green}{\bf #1}} \nc{\btpurple}[1]{\textcolor{purple}{\bf #1}}

\renewcommand{\Bbb}{\mathbb}

\newcommand{\efootnote}[1]{}

\renewcommand{\textbf}[1]{}

\nc{\mlabel}[1]{\label{#1}}  
\nc{\mcite}[1]{\cite{#1}}  
\nc{\mref}[1]{\ref{#1}}  
\nc{\mbibitem}[1]{\bibitem{#1}} 

\delete{
\nc{\mlabel}[1]{\label{#1}  
{\hfill \hspace{1cm}{\bf{{\ }\hfill(#1)}}}}
\nc{\mcite}[1]{\cite{#1}{{\bf{{\ }(#1)}}}}  
\nc{\mref}[1]{\ref{#1}{{\bf{{\ }(#1)}}}}  
\nc{\mbibitem}[1]{\bibitem[\bf #1]{#1}} 
}

\renewcommand\leq{\leqslant}

\renewcommand\bar[1]{\overline{#1}}
\renewcommand\tilde[1]{\widetilde{#1}}

\nc{\rbw}{\mathfrak{R}} \nc{\brp}{\mathrm{brp}} \nc{\lead}{\mathrm{Lead}} \nc{\Id}{\mathrm{Id}} \nc{\Irr}{\mathrm{Irr}} \nc{\vx}{\sigma} \nc{\vy}{\tau} \nc{\dvx}{\sigma^{(1)}} \nc{\dvy}{\tau^{(1)}} \nc{\done}{\vep} \nc{\citep}[1]{\cite{#1}} \nc{\wt}{\mathrm{wt}} \nc{\bre}[1]{|#1|} \nc{\mapmonoid}{\frakM} \nc{\disjoint}{\frakM'}
\nc{\ncpoly}[1]{\langle #1\rangle}  
\nc{\mapm}[1]{\lfloor\!|{#1}|\!\rfloor}
\nc{\diff}[1]{{}^\NC\{ #1 \}} \nc{\disj}[1]{\{{#1}\}'} \nc{\mdisj}[1]{\frakM'(#1)} \nc{\brho}{\bar{\rho}} \nc{\om}{\bar{\frakm}} \nc{\frakn}{\mathfrak n} \nc{\ddeg}[1]{^{(#1)}} \nc{\opset}{X} \nc{\genset}{{Z}} \nc{\NC}{\mathrm{{NC}}} \nc{\leaf}{\mathrm{leaf}} \nc{\twig}{\mathrm{twig}} \nc{\fe}{\mathrm{fl}} \nc{\munderline}[1]{#1} \nc{\bo}{o} \nc{\dep}{\mathrm{depth}} \nc{\ofe}{\mathrm{ofl}} \nc{\dfe}{\mathrm{dfe}} \nc{\fex}{\mathrm{fex}} \nc{\dl}{\mathrm{dlex}} \nc{\db}{\mathrm{db}} \nc{\lex}{\mathrm{lex}} \nc{\clex}{\mathrm{clex}} \nc{\dgp}{\mathrm{dgp}} \nc{\dgx}{\mathrm{dgx}} \nc{\br}{\mathrm{br}} \nc{\obd}{\mathrm{odb}} \nc{\ob}{\mathrm{ob}}

\nc{\bin}[2]{ (_{\stackrel{\scs{#1}}{\scs{#2}}})}  
\nc{\binc}[2]{ \left (\!\! \begin{array}{c} \scs{#1}\\
    \scs{#2} \end{array}\!\! \right )}  
\nc{\bincc}[2]{  \left ( {\scs{#1} \atop
    \vspace{-1cm}\scs{#2}} \right )}  
\nc{\bs}{\bar{S}} \nc{\cosum}{\sqsubset} \nc{\la}{\longrightarrow} \nc{\rar}{\rightarrow} \nc{\dar}{\downarrow} \nc{\dprod}{**} \nc{\dap}[1]{\downarrow \rlap{$\scriptstyle{#1}$}} \nc{\md}{\mathrm{dth}} \nc{\uap}[1]{\uparrow \rlap{$\scriptstyle{#1}$}} \nc{\defeq}{\stackrel{\rm def}{=}} \nc{\disp}[1]{\displaystyle{#1}} \nc{\dotcup}{\ \displaystyle{\bigcup^\bullet}\ } \nc{\gzeta}{\bar{\zeta}} \nc{\hcm}{\ \hat{,}\ } \nc{\hts}{\hat{\otimes}} \nc{\barot}{{\otimes}} \nc{\free}[1]{\bar{#1}} \nc{\uni}[1]{\tilde{#1}} \nc{\hcirc}{\hat{\circ}} \nc{\leng}{\ell} \nc{\lleft}{[} \nc{\lright}{]} \nc{\lc}{\lfloor} \nc{\rc}{\rfloor}
\nc{\lb}{[} 
\nc{\rb}{]} 
\nc{\curlyl}{\left \{ \begin{array}{c} {} \\ {} \end{array}
    \right.  \!\!\!\!\!\!\!}
\nc{\curlyr}{ \!\!\!\!\!\!\!
    \left. \begin{array}{c} {} \\ {} \end{array}
    \right \} }
\nc{\longmid}{\left | \begin{array}{c} {} \\ {} \end{array}
    \right. \!\!\!\!\!\!\!}
\nc{\onetree}{\bullet} \nc{\ora}[1]{\stackrel{#1}{\rar}}
\nc{\ola}[1]{\stackrel{#1}{\la}}
\nc{\ot}{\otimes} \nc{\mot}{{{\boxtimes\,}}} \nc{\otm}{\overline{\boxtimes}} \nc{\sprod}{\bullet} \nc{\scs}[1]{\scriptstyle{#1}} \nc{\mrm}[1]{{\rm #1}} \nc{\msum}{\sum\limits}
\nc{\margin}[1]{\marginpar{\rm #1}}   
\nc{\dirlim}{\displaystyle{\lim_{\longrightarrow}}\,} \nc{\invlim}{\displaystyle{\lim_{\longleftarrow}}\,} \nc{\mvp}{\vspace{0.3cm}} \nc{\tk}{^{(k)}} \nc{\tp}{^\prime} \nc{\ttp}{^{\prime\prime}} \nc{\svp}{\vspace{2cm}} \nc{\vp}{\vspace{8cm}} \nc{\proofbegin}{\noindent{\bf Proof: }}
\nc{\proofend}{$\blacksquare$ \vspace{0.3cm}}
\nc{\modg}[1]{\!<\!\!{#1}\!\!>}
\nc{\intg}[1]{F_C(#1)} \nc{\lmodg}{\!<\!\!} \nc{\rmodg}{\!\!>\!} \nc{\cpi}{\widehat{\Pi}}
\nc{\sha}{{\mbox{\cyr X}}}  
\nc{\shap}{{\mbox{\cyrs X}}} 
\nc{\shpr}{\diamond}    
\nc{\shp}{\ast} \nc{\shplus}{\shpr^+}
\nc{\shprc}{\shpr_c}    
\nc{\msh}{\ast} \nc{\zprod}{m_0} \nc{\oprod}{m_1} \nc{\vep}{\varepsilon} \nc{\labs}{\mid\!} \nc{\rabs}{\!\mid}
\nc{\astarrow}{\overset{\raisebox{-3pt}{$\ast$}}{\rightarrow}}

\nc{\dth}{d} \nc{\mmbox}[1]{\mbox{\ #1\ }} \nc{\fp}{\mrm{FP}} \nc{\rchar}{\mrm{char}} \nc{\Fil}{\mrm{Fil}} \nc{\Mor}{Mor\xspace} \nc{\gmzvs}{gMZV\xspace} \nc{\gmzv}{gMZV\xspace} \nc{\mzv}{MZV\xspace} \nc{\mzvs}{MZVs\xspace} \nc{\Hom}{\mrm{Hom}} \nc{\id}{\mrm{id}} \nc{\im}{\mrm{im}} \nc{\incl}{\mrm{incl}} \nc{\map}{\mrm{Map}} \nc{\mchar}{\rm char} \nc{\nz}{\rm NZ} \nc{\supp}{\mathrm Supp}

\nc{\Alg}{\mathbf{Alg}} \nc{\Bax}{\mathbf{Bax}} \nc{\bff}{\mathbf f} \nc{\bfk}{{\bf k}} \nc{\bfone}{{\bf 1}} \nc{\bfx}{\mathbf x} \nc{\bfy}{\mathbf y}
\nc{\base}[1]{\bfone^{\otimes ({#1}+1)}} 
\nc{\Cat}{\mathbf{Cat}} \delete{}
\nc{\detail}{\marginpar{\bf More detail}
    \noindent{\bf Need more detail!}
    \svp}
\nc{\Int}{\mathbf{Int}} \nc{\Mon}{\mathbf{Mon}}
\nc{\rbtm}{{shuffle }} \nc{\rbto}{{Rota-Baxter }} \nc{\remarks}{\noindent{\bf Remarks: }} \nc{\Rings}{\mathbf{Rings}} \nc{\Sets}{\mathbf{Sets}}

\nc{\BA}{{\Bbb A}} \nc{\CC}{{\Bbb C}} \nc{\DD}{{\Bbb D}} \nc{\EE}{{\Bbb E}} \nc{\FF}{{\Bbb F}} \nc{\GG}{{\Bbb G}} \nc{\HH}{{\Bbb H}} \nc{\LL}{{\Bbb L}} \nc{\NN}{{\Bbb N}} \nc{\KK}{{\Bbb K}} \nc{\QQ}{{\Bbb Q}} \nc{\RR}{{\Bbb R}} \nc{\TT}{{\Bbb T}} \nc{\VV}{{\Bbb V}} \nc{\ZZ}{{\Bbb Z}}

\nc{\cala}{{\mathcal A}} \nc{\calc}{{\mathcal C}} \nc{\cald}{{\mathcal D}} \nc{\cale}{{\mathcal E}} \nc{\calf}{{\mathcal F}} \nc{\calg}{{\mathcal G}} \nc{\calh}{{\mathcal H}} \nc{\cali}{{\mathcal I}} \nc{\call}{{\mathcal L}} \nc{\calm}{{\mathcal M}} \nc{\caln}{{\mathcal N}} \nc{\calo}{{\mathcal O}} \nc{\calp}{{\mathcal P}} \nc{\calr}{{\mathcal R}} \nc{\cals}{{\mathcal S}} \nc{\calt}{{\mathcal T}} \nc{\calw}{{\mathcal W}} \nc{\calk}{{\mathcal K}} \nc{\calx}{{\mathcal X}} \nc{\CA}{\mathcal{A}}

\nc{\fraka}{{\mathfrak a}} \nc{\frakA}{{\mathfrak A}} \nc{\frakb}{{\mathfrak b}} \nc{\frakB}{{\mathfrak B}} \nc{\frakD}{{\mathfrak D}} \nc{\frakH}{{\mathfrak H}} \nc{\frakM}{{\mathfrak M}} \nc{\bfrakM}{\overline{\frakM}} \nc{\frakm}{{\mathfrak m}} \nc{\frakP}{{\mathfrak P}} \nc{\frakN}{{\mathfrak N}} \nc{\frakp}{{\mathfrak p}} \nc{\frakS}{{\mathfrak S}} \nc{\frakx}{{\mathfrak x}} \nc{\ox}{\bar{\frakx}} \nc{\frakX}{{\mathfrak X}} \nc{\fraky}{{\mathfrak y}} \nc\dop{\delta}
\nc{\Reduce}{{\rm Red}}

\font\cyr=wncyr10 \font\cyrs=wncyr7
\nc{\redt}[1]{\textcolor{red}{#1}}

\nc{\li}[1]{\textcolor{red}{\tt Li:#1}} \nc{\sz}[1]{\textcolor{blue}{\tt sz:#1}} \nc{\xg}[1]{\textcolor{purple}{\tt xg:#1}}

\begin{document}
\title[Rota-Baxter operators on semigroup algebras of order two and three]{Classification of Rota-Baxter operators on semigroup algebras of order two and three}

\author{Li Guo}
\address{
Department of Mathematics and Computer Science, Rutgers University, Newark, NJ 07102, USA}
\email{liguo@rutgers.edu}

\author{Markus Rosenkranz}
\address{
 	School of Mathematics, Statistics and Actuarial Science,
	University of Kent,
	Canterbury CT2 7NF, England}
\email{M.Rosenkranz@kent.ac.uk}

\author{Shanghua Zheng}
\address{Department of Mathematics, Lanzhou University, Lanzhou, Gansu 730000, China}
\email{zheng2712801@163.com}

\hyphenpenalty=8000
\date{\today}

\begin{abstract}
  In this paper we determine all the Rota-Baxter operators of weight zero on
  semigroup algebras of order two and three with the help of computer
  algebra. We determine the matrices for these Rota-Baxter operators by directly
  solving the defining equations of the operators. We also produce a Mathematica
  procedure to predict and verify these solutions.
\end{abstract}

\subjclass[2010]{16W99, 16S36, 20M25, 16Z05}


\delete{
\begin{keyword}
Rota-Baxter algebra, semigroup, semigroup algebra, matrix, Mathematica. \end{keyword}
}

\maketitle

\tableofcontents

\hyphenpenalty=8000 \setcounter{section}{0}


\section{Introduction}
\mlabel{sec:int}
Rota-Baxter operators arose from the probability study of G.~Baxter in 1960~\mcite{Ba}, defined by the operator identity
\begin{equation}
P(x)P(y)=P(xP(y))+P(P(x)y)+\lambda P(xy),
\mlabel{eq:rbo}
\end{equation}
where $\lambda$ is a fixed scalar called the weight. When $\lambda=0$, the
operator is a natural algebraic generalization of the integral operator. In the
1960s and 70s, these operators attracted attention from well-known analysts such
as Atkinson~\mcite{At} and combinatorialists such as Cartier and
Rota~\mcite{Ca,Ro1,Ro2}. In the 1980s these operators were studied in integrable
systems as the operator form of the classical Yang-Baxter equations~\mcite{STS},
named after the well-known physicists C.N.~Yang and R.J.~Baxter. Since the late
1990s, the study of Rota-Baxter operators has made great progress both in theory
and in applications in combinatorics, number theory, operads, boundary value
problems and mathematical
physics~\mcite{AGKO,Bai,BGN,BBGN,BCQ,CK,EGK,Gub,GK1,GK3,RR}.

Because of the importance of Rota-Baxter operators, it is useful to
obtain miscellaneous examples. The examples arising naturally from
applications as well as theoretical investigations (e.g.\@ free
Rota-Baxter algebras~\mcite{EG1,GK1}) are infinite dimensional.
In recent years, some progress regarding computations of Rota-Baxter operators on low dimensional algebras has
been achieved, with applications to pre-Lie algebras,
dendriform algebras and the classical Yang-Baxter equation~\mcite{AB,LHB,PBG,TZS}..

This paper studies Rota-Baxter operators on semigroup
algebras. Semigroup algebras, a natural generalization of group
algebras, form an important class of associative algebras arising from
semigroups~\mcite{Ok}. The representation of semigroups leads to a
semigroup algebra satisfying polynomial identities. In this regards,
Rota-Baxter operators on a semigroup algebra can be regarded as an
operated semigroup algebra satisfying an operator identity. It has
been shown that every finite dimensional algebra of finite
representation type over an algebraically closed field is a contracted
semigroup algebra. Recently, semigroup algebras have experienced rapid
development on the theoretical side~\mcite{DSZZ,DE,GX,JO} as well as
in applications to representation theory, cohomology, geometric group
theory, topology, combinatorics, algebraic geometry and number
theory~\mcite{BEN,GW,Ke,LM,Ma,Ma2,May,NP,SS}. Thus studying semigroup
algebras in their canonical basis has a significance of its own.

In this paper, we classify all Rota-Baxter operators on semigroup
algebras of order~$2$ and~$3$. Rota-Baxter operators on associative
algebras of dimension $2$ and $3$ have been determined up to
isomorphism in~\mcite{LHB}. Here we focus on the particular
presentation of such an algebra in terms of the canonical semigroup
algebra basis, because of the aforementioned importance of using such
a basis. Indeed, as one notices by comparing the classifications given
here with the ones in~\mcite{LHB}, the resulting Rota-Baxter operators
take a very different form.

Because of the complex nature of Rota-Baxter operators, determining
their classification by hand is challenging even for low dimensional
algebras, as observed in~\mcite{AB,LHB,TZS}. In such a case, computer algebra provides an indispensable
aid for both predicting and verifying these operators. Nevertheless,
for ensuring theoretical accuracy, it is still necessary to carry out
a rigorous proof of the classification. In
Section~\mref{sec:RBsemitwo} we start by developing the general setup
of the equations that serve as the necessary and sufficient conditions
characterizing a Rota-Baxter operator on a semigroup algebra. We then
classify all Rota-Baxter operators on semigroup algebras of
order~$2$. For Rota-Baxter operators on semigroup algebras of
order~$3$, we carry out the classification in two sections. In
Section~\mref{sec:rbocs}, we determine all Rota-Baxter operators on
commutative semigroup algebras of order~$3$, in
Section~\mref{sec:rboncs} on noncommutative semigroup algebras of
order~$3$.
In Section~\mref{sec:comp}, we provide the Mathematica procedure that has helped us in solving the classification problem.

\allowdisplaybreaks

\section{The general setup and Rota-Baxter operators on semigroup algebras of order two}
\mlabel{sec:RBsemitwo}

In this section, we first formulate the general setup for determining
Rota-Baxter operators of weight zero on a semigroup algebra. We then
apply this setup to classify all Rota-Baxter operators of weight zero
on semigroup algebras of order two.

\subsection{The general setup}

In this subsection, we give the general setup of Rota-Baxter operators
on a semigroup algebra in matrix form. Let $S$ be a finite semigroup
with multiplication $\cdot$ that we often suppress. Thus
$S=\{e_1,\cdots,e_n\}$. Let $\bfk$ be a commutative unitary ring and let
\begin{equation}
\bfk[S]:=\sum_{m=1}^n \bfk e_m=\left\{\sum_{m=1}^n a_m e_m\,\big|\, a_m\in \bfk, 1\leq m\leq n\right \}
\mlabel{eq:sga}
\end{equation}
denote the semigroup algebra of $S$. The order $n$ of the semigroup
$S$ is also said to be the {\bf order} of the semigroup algebra
$\bfk[S]$.

Let $P:\bfk[S]\to \bfk[S]$ be a Rota-Baxter operator of weight zero. Since $P$ is $\bfk$-linear, we have
\begin{equation}
\left(\begin{matrix}
P(e_1)\\
P(e_2)\\
\cdots\\
P(e_n)
\end{matrix}\right)
=
\left(
\begin{matrix}
c_{11}&\cdots & c_{1n}\\
c_{21}&\cdots & c_{2n}\\
\cdots & \cdots & \cdots\\
c_{n1}&\cdots & c_{nn}
\end{matrix}
\right)
\left(\begin{matrix}
e_1\\
e_2\\
\cdots\\
e_n
\end{matrix}\right)\quad (c_{ij}\in\bfk, 1\leq i,j\leq n).
\mlabel{eq:rbma2}
\end{equation}
The matrix $C:=C_P:=\big( c_{ij}\big)_{1\leq i, j\leq n}$ is called the {\bf matrix of $P$}.
Further, $P$ is a Rota-Baxter operator if and only if
\begin{equation}
P(e_i)P(e_j)=P(P(e_i)e_j+e_iP(e_j))\quad ( 1\leq i,j\leq n).
\mlabel{eq:rbobn}
\end{equation}
Let the Cayley (multiplication) table of the semigroup $S$ be given by
\begin{equation}
e_k\cdot e_\ell=\sum_{m=1}^nr_{k\ell}^me_m\quad  (1\leq k,\ell\leq n),
\mlabel{eq:ctablen}
\end{equation}
where~$r^m_{k\ell} \in \{0,1\}$. Then we have
$$
P(e_i)P(e_j)=\sum_{k=1}^n\sum_{\ell=1}^n
c_{ik}c_{j\ell}
e_ke_\ell=\sum_{m=1}^n\sum_{k=1}^n
\sum_{\ell=1}^nr_{k\ell}^m c_{ik}c_{j\ell}e_m
$$
and
{\allowdisplaybreaks
\begin{eqnarray*}
P(P(e_i)e_j+e_iP(e_j))
&=&\sum_{k=1}^nc_{ik}P(e_ke_j)+\sum_{\ell=1}
^nc_{j\ell}P(e_ie_\ell)\\
&=&\sum_{k=1}^n
\sum_{m=1}^nr_{kj}^mc_{ik}P(e_m)+\sum_{\ell=1}
^n\sum_{m=1}
^nr_{i\ell}^mc_{j\ell}P(e_m)\\
&=&\sum_{k=1}^n
\sum_{m=1}^nr_{k j}^mc_{ik}\bigg(\sum_{\ell=1}^n
c_{m\ell}e_\ell\bigg)+\sum_{k=1}
^n\sum_{m=1}
^nr_{ik}^mc_{jk}\bigg(\sum_{\ell=1}^n
c_{m\ell}e_\ell\bigg)\\
&=&\sum_{m=1}^n
\sum_{\ell=1}^n\sum_{k=1}^n (r_{k j}^\ell c_{ik}
+r_{ik}^\ell c_{jk})
c_{\ell m}e_m.
\end{eqnarray*}
}
Thus we obtain
\begin{theorem}
Let $S=\{e_1,\cdots,e_n\}$ be a semigroup with its Cayley table given by Eq.~(\mref{eq:ctablen}). Let $\bfk$ be a commutative unitary ring and let $P:\bfk[S]\to \bfk[S]$ be a linear operator with matrix $C:=C_P=(c_{ij})_{1\leq i, j\leq n}$. Then $P$ is a Rota-Baxter operator of weight zero on $\bfk[S]$ if and only if the following equations hold.
\begin{equation}
\sum_{\ell =1}^n\sum_{k=1}^n
r_{k\ell}^m c_{ik}c_{j\ell}=\sum_{\ell=1}^n
\sum_{ k=1}^n( r_{k j}^\ell c_{ik}+r_{ik}^\ell c_{jk})
c_{\ell m} \quad (1\leq i,j,m\leq n). \mlabel{eq:neijm}\\
\end{equation}
\mlabel{thm:rbon}
\end{theorem}

We will determine the matrices $C_P$ for all Rota-Baxter operators $P$ on $\bfk[S]$ of order two or three.

\subsection{Rota-Baxter operators on semigroup algebras of order 2}
\mlabel{subsec:semitwo}
In this section, we determine all Rota-Baxter operators on semigroup
algebras $\bfk[S]$ of order~$2$.

As is well known~\mcite{Pe}, there are exactly five distinct
nonisomorphic semigroups of order $2$. We use $N_2, L_2, R_2, Y_2$ and
$Z_2$ respectively to denote the null semigroup of order $2$, the left
zero semigroup, right zero semigroup, the semilattice of order $2$ and
the cyclic group of order $2$. Since $L_2$ and $R_2$ are
anti-isomorphic, there are exactly four distinct semigroups of order
$2$, up to isomorphism and anti-isomorphism, namely $N_2$, $Y_2$,
$Z_2$ and $L_2$.  Let $\{e_1,e_2\}$ denote the underlying set of each
semigroup. Then the Cayley tables for these semigroups are as follows:
\begin{table}[!tbhp]
\caption{The Cayley tables of semigroups of order 2 \mlabel{ta:s2}}
\begin{minipage}[t]{0.85\textwidth}
\begin{tabular}{|c|c|c|c|}
\hline
$N_2:=$
\begin{tabular}{c|ccc}
$\cdot$ &$e_1$ & $e_2$\\
\hline
$e_1$ &$e_1$ & $e_1$ \\
$e_2$ &$e_1$ & $e_1$ \\
\end{tabular}
&$Y_2:=$\begin{tabular}{c|ccc}
$\cdot$ &$e_1$ & $e_2$\\
\hline
$e_1$ &$e_1$ & $e_1$ \\
$e_2$ &$e_1$ & $e_2$ \\
\end{tabular}
&$Z_2:=$
\begin{tabular}{c|ccc}
$\cdot$ &$e_1$ & $e_2$\\
\hline
$e_1$ &$e_1$ & $e_2$ \\
$e_2$ &$e_2$ & $e_1$ \\
\end{tabular}
&$L_2:=$
\begin{tabular}{c|ccc}
$\cdot$ &$e_1$ & $e_2$\\
\hline
$e_1$ &$e_1$ & $e_1$ \\
$e_2$ &$e_2$ & $e_2$ \\
\end{tabular}\\
\hline
\end{tabular}
\end{minipage}
\end{table}

\begin{theorem}
Let $\bfk$ be a field of characteristics zero. All Rota-Baxter operators on a semigroup algebra $\bfk[S]$ of order $2$ have their matrices $C_P$ given in  Table~\mref{ta:rbs2},
where all the parameters are in $\bfk$ and RBO $($resp. SA$)$ is the abbreviation of Rota-Baxter operator $($resp. semigroup algebra$)$.
\mlabel{thm:rbs2}
\end{theorem}

{\allowdisplaybreaks
\begin{table}[!tbhp]
\caption{Table of RBOs on $2$-dimensional semigroup algebras \mlabel{ta:rbs2}}
\tiny
\begin{tabular}{|c| c|c|c|p{13cm}|}
\hline
\text{Semigroup of order 2}
&\text{ RBOs on  SA}
&\text{Semigroup of order 2}
&\text{ RBOs on  SA}  \\
\hline
$N_2$
&$\left(\begin{matrix}
a & -a \\
b & -b \\
\end{matrix}\right)$
&$Z_2$
&$\left(\begin{matrix}
0 & 0 \\
0 & 0 \\
\end{matrix}\right)
$\\
\hline
$Y_2$
&$\left(\begin{matrix}
0 & 0 \\
0 & 0 \\
\end{matrix}\right)$
&$L_2$
&$\left(\begin{matrix}
a & -\frac{a^2}{b} \\
b & -a \\
\end{matrix}\right)(b\neq0),
\left(\begin{matrix}
0 & a \\
0 & 0 \\
\end{matrix}\right)
$\\
\hline
\end{tabular}
\end{table}
}
\begin{proof}

  We divide the proof of the theorem into four cases, one for each of the four
  semigroups $S$ in Table~\mref{ta:s2}.  For each case, by
  Theorem~\mref{thm:rbon}, $P$ is a RBO on $\bfk[S]$ if and only if the eight
  equations~(\mref{eq:neijm}) hold (with~$1 \leq i, j, m \leq 2$).  So we just
  need to solve these equations. It is straightforward to verify that what we
  obtain does indeed satisfy all equations. Let $\mathbf{0}_{2\times 2}$ denote
  the $2\times 2$ zero matrix.  \smallskip

\noindent
{\bf Case 1. Let  $S=N_2$:}  In Eq.~(\mref{eq:neijm}), taking $i=j=1$ with $1\leq m\leq 2$ and $i=j=2$ with $1\leq m\leq 2$, we get
\begin{eqnarray}
&(c_{11}+c_{12})^2=2c_{11}(c_{11}+c_{12}),
\mlabel{eq:cs211}\\
&c_{12}(c_{11}+c_{12})=0,\mlabel{eq:cs212}\\
&(c_{21}+c_{22})^2=2c_{11}(c_{21}+c_{22}),\mlabel{eq:cs213}\\
&c_{12}(c_{21}+c_{22})=0.\mlabel{eq:cs214}
\end{eqnarray}
Assume $c_{11}+c_{12}\neq 0$. Then by Eqs.~(\mref{eq:cs211}) and ~(\mref{eq:cs212}), we have $c_{11}=c_{12}=0$, a contradiction to $c_{11}+c_{12}\neq 0$.  Thus $c_{11}+c_{12}=0$. Assume $c_{21}+c_{22}\neq 0$. Then by Eq.~(\mref{eq:cs214}), we get $c_{12}=0$, and so $c_{11}=0$. Then by Eq.~(\mref{eq:cs213}), we have $c_{21}+c_{22}=0$, a contradiction. Thus $c_{21}+c_{22}=0$.  Therefore Eqs.~(\mref{eq:cs211})-(\mref{eq:cs214}) are equivalent to
$$\left\{\begin{array}{cc}
c_{11}+c_{12}=0,\\
c_{21}+c_{22}=0.\\
\end{array}\right.$$
Denoting $a=c_{11}$ and $b=c_{21}$, we see that solutions $(c_{ij})_{1\leq i, j\leq 2}$ of Eqs.~(\mref{eq:cs211})-(\mref{eq:cs214}) are of the form
$$\left(
\begin{matrix}
a&-a\\
b&-b\\
\end{matrix}\right)\quad (a, b\in \bfk).$$
It is straightforward to check that they also satisfy the other equations in Eq.~(\mref{eq:neijm}). Hence these are all the matrices $C_P$ for Rota-Baxter operators on $\bfk[S]$.

\smallskip
\noindent
{\bf Case 2. Let $S=Y_2$:} In Eq.~(\mref{eq:neijm}), taking $i=j=1$ with $1\leq m\leq 2$ and $i=j=2$ with $1\leq m\leq 2$, we obtain
\begin{eqnarray}
&c_{11}^2+2c_{11}c_{12}=2c_{11}(c_{11}+c_{12}),
\mlabel{eq:cs221}\\
&c_{12}^2=2c_{12}(c_{11}+c_{12}),\mlabel{eq:cs222}\\
&c_{21}^2+2c_{21}c_{22}=2c_{21}(c_{11}+c_{22}),\mlabel{eq:cs223}\\
&c_{22}^2=2(c_{12}c_{21}+c_{22}^2).\mlabel{eq:cs224}
\end{eqnarray}
From Eq.~(\mref{eq:cs221}) we have $c_{11}=0$.  Then from Eq.~(\mref{eq:cs222}) we have $c_{12}=0$.  Thus Eq.~(\mref{eq:cs224}) gives $c_{22}=0$. Further by Eq.~(\mref{eq:cs223}), we have $c_{21}=0$. Thus the only solution is the zero solution $\mathbf{0}_{2\times 2}$.

\smallskip
\noindent
{\bf Case 3. $S=Z_2$:} In Eq.~(\mref{eq:neijm}), taking $i=j=1$ with $1\leq m\leq 2$; $i=1,j=2$ with $m=2$ and $i=j=2$ with $m=1$, we obtain
\begin{eqnarray}
&c_{11}^2+c_{12}^2=2c_{11}^2+
2c_{12}c_{21},
\mlabel{eq:cs231}\\
&c_{11}c_{12}=c_{12}( c_{11}+c_{22}),\mlabel{eq:cs232}\\
&c_{11}c_{22}+c_{12}c_{21}
=c_{22}(c_{11}+c_{22})+c_{12}(c_{12}+c_{21})
,\mlabel{eq:cs233}\\
&c_{21}^2+c_{22}^2=2(c_{21}^2+c_{11}c_{22}).
\mlabel{eq:cs234}
\end{eqnarray}
By Eq.~(\mref{eq:cs232}) we have $c_{12}c_{22}=0$. From Eq.~(\mref{eq:cs233}) we get $c_{12}^2+c_{22}^2=0$. Thus $c_{12}=c_{22}=0$. Then  Eqs.~(\mref{eq:cs231}) and ~(\mref{eq:cs234}) give $c_{11}=c_{21}=0$.  Thus the only solution is the zero solution $\mathbf{0}_{2\times 2}.$
\smallskip

\noindent
{\bf Case 4. $S=L_2$:}  In Eq.~(\mref{eq:neijm}), taking $i=j=1$ with $1\leq m\leq 2$ and $i=2,j=1$ with $1\leq m\leq 2$,
we get
\begin{eqnarray}
&c_{11}(c_{11}+c_{12})=c_{11}(2c_{11}+c_{12})
+c_{12}c_{21},
\mlabel{eq:cs241}\\
&c_{12}(c_{11}+c_{12})=c_{12}( 2c_{11}+c_{12}+c_{22}),\mlabel{eq:cs242}\\
&c_{21}(c_{11}+c_{12})
=c_{21}(2c_{11}+c_{12}+c_{22})
,\mlabel{eq:cs243}\\
&c_{22}(c_{11}+c_{12})
=c_{21}c_{12}+c_{22}(c_{11}+c_{12}+c_{22}).
\mlabel{eq:cs244}
\end{eqnarray}
By Eqs.~(\mref{eq:cs241}) and ~(\mref{eq:cs244}), we have $c_{11}^2=c_{22}^2$. By Eqs.~(\mref{eq:cs242}) and ~(\mref{eq:cs243}), we get $c_{12}(c_{11}+c_{22})=0$ and $c_{21}(c_{11}+c_{22})=0$. Assume $c_{11}+c_{22}\neq0$. Then $c_{12}=c_{21}=0$. So by Eq.~(\mref{eq:cs241}), we have $c_{11}=0$ and then $c_{22}=0$, a contradiction. Thus $c_{11}+c_{22}=0$. Then Eqs.~(\mref{eq:cs241})-(\mref{eq:cs244}) are equivalent to the system of equations
$$\left\{
\begin{array}{cc}
c_{11}+c_{22}=0,\\
c_{11}^2+c_{12}c_{21}=0.\\
\end{array}\right.$$
Denoting $a=c_{11}$ and $b=c_{21}$, then $c_{22}=-a$. When $b\neq 0$, then we also have $c_{12}=-\frac{a^2}{b}$.
This gives the solutions
$$\left(
\begin{matrix}
a&-\frac{a^2}{b}\\
b&-a\\
\end{matrix}\right)\quad (b\neq 0, a\in\bfk).$$
On the other hand, when $b=0$, then $c_{11}=c_{22}=0$. Denoting $a=c_{12}$, we get the solutions
$$\left(
\begin{matrix}
0&a\\
0&0\\
\end{matrix}\right)\quad ( a\in\bfk).$$
These solutions to Eqs.~(\mref{eq:cs241})-(\mref{eq:cs244})  also satisfy the other equations in Eq.~(\mref{eq:neijm}). Thus they give all the Rota-Baxter operators on $\bfk[S]$.
\smallskip

This completes the proof of Theorem~\mref{thm:rbs2}.
\end{proof}

\section{Rota-Baxter operators on commutative semigroup algebras of order $3$}
\mlabel{sec:rbocs}
Up to isomorphism and anti-isomorphism, there are 18 semigroups of order $3$~\mcite{CS,Edm,Fo}. The Cayley tables of the 18 semigroups of order $3$ can be found in~\mcite{Fo}. See also~\mcite{CS, LZ,Ple}.
We denote by $CS$ and $NCS$ the class of 12 commutative semigroups and the class of 6 noncommutative semigroups, respectively.

When a semigroup $S$ has order $3$, the equations in Eq.~(\mref{eq:neijm}) for the matrix $C_P$ of a Rota-Baxter operator $P$ on $\bfk[S]$ are given by the following 27 equations.
\begin{equation}
\sum_{k=1}^3\sum_{\ell=1}^3
r_{k\ell}^m c_{ik}c_{j\ell}=\sum_{k=1}^3
\sum_{\ell=1}^3( r_{k j}^\ell c_{ik}+r_{ik}^\ell c_{jk})
c_{\ell m}\quad (1\leq i,j,m\leq 3).\mlabel{eq:3eijm}
\end{equation}

In this section, we determine the Rota-Baxter operators on the semigroup algebras for the 12 commutative semigroups of order 3 in Table~\mref{ta:s3}. Rota-Baxter operators on semigroup algebras for the 6 noncommutative semigroups of order 3 will be determined in Section~\mref{sec:rboncs}.

\subsection{Statement of the classification theorem in the commutative case}
\mlabel{ss:sord3}

A classification of the 12 commutative semigroups of order $3$
is given in Table~\mref{ta:s3}.

\begin{tiny}
\begin{table}[!thbp]
\caption{Table of commutative semigroups of order 3 \mlabel{ta:s3}}
\begin{minipage}[tbph]{0.8\textwidth}
\begin{tabular}{|c|c|c|c|}
\hline
$CS(1):=$
\begin{tabular}{c|cccc}
$\cdot$ &$e_1$ & $e_2$& $e_3$\\
\hline
$e_1$ &$e_1$ & $e_1$ &$e_1$\\
$e_2$ &$e_1$ & $e_1$ &$e_1$\\
$e_3$ &$e_1$ & $e_1$ &$e_1$\\
\end{tabular}
&$CS(2):=$
\begin{tabular}{c|cccc}
$\cdot$ &$e_1$ & $e_2$& $e_3$\\
\hline
$e_1$ &$e_1$& $e_1$ &$e_1$\\
$e_2$ &$e_1$ & $e_1$ &$e_1$\\
$e_3$ &$e_1$ & $e_1$ &$e_2$\\
\end{tabular}
&$CS(3):=$
\begin{tabular}{c|ccccc}
$\cdot$ &$e_1$ & $e_2$& $e_3$\\
\hline
$e_1$ &$e_1$ & $e_1$ &$e_1$\\
$e_2$ &$e_1$ & $e_2$ &$e_1$\\
$e_3$ &$e_1$ & $e_1$ &$e_1$\\
\end{tabular}\\
\hline
$CS(4):=$
\begin{tabular}{c|ccccc}
$\cdot$ &$e_1$ & $e_2$& $e_3$\\
\hline
$e_1$ &$e_1$ & $e_1$ &$e_1$\\
$e_2$ &$e_1$ & $e_2$ &$e_1$\\
$e_3$ &$e_1$ & $e_1$ &$e_3$\\
\end{tabular}
&$CS(5):=$
\begin{tabular}{c|ccccc}
$\cdot$ &$e_1$ & $e_2$& $e_3$\\
\hline
$e_1$ &$e_1$& $e_1$ &$e_1$\\
$e_2$ &$e_1$ & $e_2$ &$e_2$\\
$e_3$ &$e_1$ & $e_2$ &$e_2$\\
\end{tabular}
&$CS(6):=$
\begin{tabular}{c|ccccc}
$\cdot$ &$e_1$ & $e_2$& $e_3$\\
\hline
$e_1$ &$e_1$ & $e_1$ &$e_1$\\
$e_2$ &$e_1$ & $e_2$ &$e_2$\\
$e_3$ &$e_1$ & $e_2$ &$e_3$\\
\end{tabular}\\
\hline
$CS(7):=$
\begin{tabular}{c|ccccc}
$\cdot$ &$e_1$ & $e_2$& $e_3$\\
\hline
$e_1$ &$e_1$ & $e_1$ &$e_1$\\
$e_2$&$e_1$ & $e_2$ &$e_3$\\
$e_3$ &$e_1$ & $e_3$ &$e_1$\\
\end{tabular}
&$CS(8):=$
\begin{tabular}{c|ccccc}
$\cdot$ &$e_1$ & $e_2$& $e_3$\\
\hline
$e_1$ &$e_1$ & $e_1$ &$e_1$\\
$e_2$ &$e_1$ & $e_2$ &$e_3$\\
$e_3$ &$e_1$ & $e_3$ &$e_2$\\
\end{tabular}
&$CS(9):=$
\begin{tabular}{c|ccccc}
$\cdot$ &$e_1$ & $e_2$& $e_3$\\
\hline
$e_1$ &$e_1$ & $e_1$ &$e_3$\\
$e_2$ &$e_1$ & $e_1$ &$e_3$\\
$e_3$ &$e_3$& $e_3$ &$e_1$\\
\end{tabular}\\
\hline
$CS(10):=$
\begin{tabular}{c|ccc}
$\cdot$ &$e_1$ & $e_2$& $e_3$\\
\hline
$e_1$ &$e_1$ & $e_1$ &$e_3$\\
$e_2$ &$e_1$ & $e_2$ &$e_3$\\
$e_3$&$e_3$ & $e_3$ &$e_1$\\
\end{tabular}
&$CS(11):=$
\begin{tabular}{c|ccccc}
$\cdot$ &$e_1$ & $e_2$& $e_3$\\
\hline
$e_1$ &$e_1$ & $e_2$ &$e_2$\\
$e_2$ &$e_2$ & $e_1$ &$e_1$\\
$e_3$&$e_2$ & $e_1$ &$e_1$\\
\end{tabular}
&$CS(12):=$
\begin{tabular}{c|ccccc}
$\cdot$ &$e_1$ & $e_2$& $e_3$\\
\hline
$e_1$ &$e_1$ & $e_2$ &$e_3$\\
$e_2$ &$e_2$ & $e_3$ &$e_1$\\
$e_3$ &$e_3$ & $e_1$ &$e_2$\\
\end{tabular}\\
\hline
\end{tabular}
\end{minipage}
\end{table}
\end{tiny}

We have the following classification of Rota-Baxter operators on order 3
commutative semigroup algebras. The proof will be given in
Section~\mref{ss:proofcs}.

\begin{theorem}
Let $\bfk$ be a field of characteristic zero. The matrices of Rota-Baxter operators on $3$-dimensional commutative semigroup algebras are given in Table~\mref{ta:rbcs}, where all the parameters take values in $\bfk$
and RBO $($resp. CS$)$ is the abbreviation for Rota-Baxter operator $($resp. commutative semigroup$)$.
\mlabel{thm:rboc3}
\end{theorem}

\begin{tiny}
\begin{longtable}
{|p{1cm}|p{3cm}| p{1cm}| p{5.5cm}|}
\caption[RBOs on 3-dimensional noncomumutative semigroup algebras ]{{RBOs on commutative semigroup algebras of order 3} \mlabel{ta:rbcs}}\\
\hline \multicolumn{1}{|c|}{{\bf CS of order 3}} & \multicolumn{1}{c|}{{\bf Matrices of RBOs on semigroup algebras}}
&\multicolumn{1}{|c|}{{\bf CS of order 3}} & \multicolumn{1}{c|}{{\bf Matrices of RBOs on semigroup algebras}}\\
\hline
\endfirsthead

\multicolumn{4}{c}
{{ \tablename\  \thetable{:\ RBOs on commutative semigroup algebras of order 3}}} \\
\hline \multicolumn{1}{|c|}{{\bf CS of order 3 }} &
\multicolumn{1}{c|}{{\bf Matrices of RBOs on semigroup algebras}}
&\multicolumn{1}{|c|}{{\bf CS of order 3 }} &
\multicolumn{1}{c|}{{\bf Matrices of RBOs on semigroup algebras}}\\
\hline
\endhead

\hline \multicolumn{4}{|r|}{{Continued on next page}} \\ \hline\endfoot

\hline
\endlastfoot

\delete{
\begin{tiny}
\begin{table}[!thbp]
\caption{RBOs on 3-dimensional commutative semigroup algebras\mlabel{ta:rbcs}}
\begin{tabular}{|c|c|c|c|}
\hline
\text{CS of order 3}
&\text{RBOs on  SA}
&\text{CS of order 3}
&\text{ RBOs on  SA}  \\
\hline
}
$CS(1)$
&
$C_{1,1}:=\left(
\begin{matrix}
a & b & -a-b\\
c & d & -c-d\\
e & f & -e-f\\
\end{matrix}
\right)
$
&$CS(2)$
&\begin{tabular}{ll}
&$C_{2,1}:=\left(\begin{matrix}
a & -a & 0\\
b & -b & 0\\
c & -c & 0\\
\end{matrix}
\right),$\\
&$C_{2,2}:=\left(\begin{matrix}
a & -a & 0\\
b & -b & 0\\
c & 2(b-a)-c & 2(a-b)\\
\end{matrix}\right)(a\neq b)$\end{tabular}\\
\hline
$CS(3)$
&
$C_{3,1}:=\left(\begin{matrix}
a & 0 & -a\\
b & 0 & -b\\
c & 0 & -c\\
\end{matrix}
\right)$
&
$CS(4)$
&
$C_{4,1}:=\left(\begin{matrix}
0 & 0 & 0\\
0 & 0 & 0\\
0 & 0 & 0\\
\end{matrix}\right)$
\\
\hline
$CS(5)$
&
$C_{5,1}:=\left(\begin{matrix}
0 &a & -a\\
0 &b & -b\\
0 & c & -c\\
\end{matrix}\right)$
&
$CS(6)$
&
$C_{6,1}:=\left(\begin{matrix}
0 & 0 & 0\\
0 & 0 & 0\\
0 & 0 & 0\\
\end{matrix}\right)$
\\
\hline
$CS(7)$
&
$C_{7,1}:=\left(\begin{matrix}
a & 0 & -a\\
b & 0 & -b\\
a & 0 & -a\\
\end{matrix}\right)$
&
$CS(8)$
&
$C_{8,1}:=\left(\begin{matrix}
0 & 0 & 0\\
0 & 0 & 0\\
0 & 0 & 0\\
\end{matrix}\right)$
\\
\hline
$CS(9)$
&
$C_{9,1}:=\left(\begin{matrix}
a & -a  & 0\\
b  & -b & 0\\
c & -c  & 0\\
\end{matrix}\right)$
&
$CS(10)$
&
$C_{10,1}:=\left(\begin{matrix}
0 & 0 & 0\\
0 & 0 & 0\\
0 & 0 & 0\\
\end{matrix}\right)$
\\
\hline
$CS(11)$
&
$C_{11,1}:=\left(\begin{matrix}
0  &a  & -a\\
0  &b  & -b\\
0  &c  & -c\\
\end{matrix}\right)$
&
$CS(12)$
&
$C_{12,1}:=\left(\begin{matrix}
0  & 0  & 0\\
0  & 0  & 0\\
0  &0  & 0\\
\end{matrix}\right)$
\\
\hline
\end{longtable}
\end{tiny}

\subsection{Proof of Theorem~\mref{thm:rboc3}}
\mlabel{ss:proofcs} We will prove Theorem~\mref{thm:rboc3} by considering each
of the $12$ commutative semigroups $CS(i), 1\leq i\leq 12,$ of order 3 in
Table~\mref{ta:s3}. For each semigroup, we solve some of the equations in
Eq.~(\mref{eq:3eijm}) for the Cayley table of the corresponding semigroup. It is
straightforward to verify that what we obtain this way indeed satisfies all the
equations in Eq.~(\mref{eq:3eijm}). Let $\mathbf{0}_{3\times 3}$ denote the
$3\times 3$ zero matrix.

\subsubsection{The proof for $\bfk[CS(1)]$}
\mlabel{ss:cs1}
We prove that the matrices $C_P=(c_{ij})_{1\leq i, j\leq 3}$ of all the Rota-Baxter operators $P$ on the semigroup algebra $\bfk[CS(1)]$ are given by $C_{1,1}$ in Table~\mref{ta:rbcs}.

Applying the Cayley table of $CS(1)$ in Eqs.~(\mref{eq:3eijm}) and taking $i=j=1$ with $1\leq m\leq 3$; $i=1$,$j=2,3$ with $m=1$  and $i=2,3$, $j=1$ with $m=1$, respectively, we obtain
\allowdisplaybreaks{
\begin{eqnarray}
&(c_{11}+c_{12}+c_{13})^2=2c_{11}(c_{11}+c_{12}
+c_{13}),\mlabel{eq:cs311}\\
&c_{12}(c_{11}+c_{12}+c_{13})=0,
\mlabel{eq:cs312}\\
&c_{13}(c_{11}+c_{12}+c_{13})=0,
\mlabel{eq:cs313}\\
&(c_{11}+c_{12}+c_{13})
(c_{21}+c_{22}+c_{23})= c_{11}(c_{11}+c_{12}+c_{13}+
c_{21}+c_{22}+c_{23}),\mlabel{eq:cs314}\\
&(c_{11}+c_{12}+c_{13})(c_{31}+c_{32}+c_{33})
=c_{11}(c_{11}+c_{12}+c_{13}+
c_{31}+c_{32}+c_{33}),\mlabel{eq:cs317}\\
&(c_{21}+c_{22}+c_{23})^2=2c_{11}
(c_{21}+c_{22}+c_{23}),\mlabel{eq:cs3110}\\
&(c_{31}+c_{32}+c_{33})^2
=2c_{11}(c_{31}+c_{32}+c_{33}).
\mlabel{eq:cs31161}
\end{eqnarray}

Assume $c_{11}+c_{12}+c_{13}\neq 0.$
Then by Eqs.~(\mref{eq:cs311}),~(\mref{eq:cs312}) and ~(\mref{eq:cs313}), we get $$2c_{11}=c_{11}+c_{12}+c_{13}\quad \text{and}\quad
c_{12}=c_{13}=0.$$
So we have $c_{11}=c_{12}=c_{13}=0$. This contradicts  $c_{11}+c_{12}+c_{13}\neq 0$.
Thus $c_{11}+c_{12}+c_{13}=0$.
By Eqs.~(\mref{eq:cs314}) and ~
(\mref{eq:cs317}), we have $c_{11}(c_{21}+c_{22}+c_{23})=0$ and $c_{11}(c_{31}+c_{32}+c_{33})=0.$   Then by Eqs.~(\mref{eq:cs3110}) and~(\mref{eq:cs31161}), we get
$(c_{21}+c_{22}+c_{23})^2=0$ and  $(c_{31}+c_{32}+c_{33})^2=0$. So we have $c_{21}+c_{22}+c_{23}=0$ and $
c_{31}+c_{32}+c_{33}=0$. From these discussions, we see that Eqs.~(\mref{eq:cs311})-(\mref{eq:cs31161}) are equivalent to
$$\left\{\begin{array}{cccc}
c_{11}+c_{12}+c_{13}=0,\\
c_{21}+c_{22}+c_{23}=0,\\
c_{31}+c_{32}+c_{33}=0.\\
\end{array}\right.$$
}
So we have
$$\left\{\begin{array}{cccc}
c_{13}=-c_{11}-c_{12},\\
c_{23}=-c_{21}-c_{22},\\
c_{33}=-c_{31}-c_{32}.\\
\end{array}\right.$$
Denote $a =c_{11},b =c_{12},c=c_{21},d=c_{22},
e= c_{31},f= c_{32}$. Thus solutions $(c_{ij})_{1\leq i, j\leq 3}$ of Eqs.~(\mref{eq:cs311})-(\mref{eq:cs31161}) are given by
$$C_{11}=\left(
\begin{matrix}
a & b & -a-b\\
c & d & -c-d\\
e & f & -e-f\\
\end{matrix}
\right)\quad (a,b,c,d,e,f\in\bfk).$$
Since they can be checked to satisfy other equations in Eq.~(\mref{eq:3eijm}), they give the matrices of all the Rota-Baxter operators on $\bfk[CS(1)]$.

\subsubsection{The proof for $\bfk[CS(2)]$}
\mlabel{ss:cs2}
Here we prove that the matrices $C_P=(c_{ij})_{1\leq i, j\leq 3}$ of all the Rota-Baxter operators $P$ on the semigroup algebra $\bfk[CS(2)]$ are given by $C_{2,1}$ and $C_{2,2}$ in Table~\mref{ta:rbcs}.

Applying the Cayley table of $CS(2)$ in Eqs.~(\mref{eq:3eijm}) and taking $i=j=1$ with $1\leq m\leq 3$ and $1\leq i\leq 3, j=2$ with $m=1,2$, we obtain
\begin{eqnarray}
&(c_{11}+c_{12})(c_{11}+c_{12}+c_{13})
+c_{13}(c_{11}+c_{12})=2c_{11}
(c_{11}+c_{12}+c_{13}),\mlabel{eq:cs321}\\
&c_{13}^2=2c_{12}(c_{11}+c_{12}+c_{13}),\mlabel{eq:cs322}\\
&c_{13}(c_{11}+c_{12}+c_{13})=0,\mlabel{eq:cs323}\\
&(c_{11}+c_{12})(c_{21}+c_{22}+c_{23})
+c_{13}(c_{21}+c_{22})=
c_{11}(c_{11}+c_{12}+c_{13}+c_{21}+c_{22}
+c_{23}),\mlabel{eq:cs324}\\
&c_{12}(c_{11}+c_{12}+c_{13}+c_{21}+c_{22}
+c_{23})=0,\mlabel{eq:cs325}\\
&(c_{21}+c_{22})(c_{21}+c_{22}+c_{23})+
c_{23}(c_{21}+c_{22})=2c_{11}(c_{21}+c_{22}
+c_{23}),\mlabel{eq:cs3210}\\
&c_{23}^2=2c_{12}(c_{21}+c_{22}+c_{23}),\mlabel{eq:cs3211}\\
&(c_{31}+c_{32})(c_{31}+c_{32}+c_{33})+
c_{33}(c_{31}+c_{32})=2c_{11}(c_{31}+c_{32}
)+2c_{33}c_{21},\mlabel{eq:cs3216}\\
&c_{33}^2=2c_{12}
(c_{31}+c_{32})+2c_{33}c_{22}.\mlabel{eq:cs3217}
\end{eqnarray}
Assume $c_{11}+c_{12}+c_{13}\neq 0$. Then by Eqs.~(\mref{eq:cs322}) and~(\mref{eq:cs323}), we have $c_{13}=c_{12}=0$. By Eq.~(\mref{eq:cs321}), we have $c_{11}=0$. This contradicts $c_{11}+c_{12}+c_{13}\neq 0$. Thus we have $c_{11}+c_{12}+c_{13}=0$. By Eq.~(\mref{eq:cs322}),  we have $c_{13}=0$. So $c_{11}+c_{12}=0$. Then by Eqs.~(\mref{eq:cs324}) and~(\mref{eq:cs325}), we get $c_{11}(c_{21}+c_{22}+c_{23})=0$ and $c_{12}(c_{21}+c_{22}+c_{23})=0$. From Eq.~(\mref{eq:cs3211}), we obtain $c_{23}=0$. Then Eq.~(\mref{eq:cs3210}) gives $c_{21}+c_{22}=0$. Adding Eqs.~(\mref{eq:cs3216}) and ~(\mref{eq:cs3217}), we get $(c_{31}+c_{32}+c_{33})^2=0$. So $c_{31}+c_{32}+c_{33}=0$. By Eq.~(\mref{eq:cs3217}), we have $c_{33}^2=2c_{33}(c_{22}-c_{12})$.
Thus $c_{33}=0$ or $c_{33}=2(c_{22}-c_{12}).$
Denote $a=c_{11}, b=c_{21}, c=c_{31}$. Then $c_{12}=-a$ and $c_{22}=-b$. We consider two cases.

\smallskip
\noindent
{\bf Case 1. Suppose $c_{33}=0$:} Then $c_{31}+c_{32}=0$. Then we have $c_{32}=-c$. Thus we get the solutions
$$C_{21}=\left(
\begin{matrix}
a & -a & 0\\
b & -b & 0\\
c & -c & 0\\
\end{matrix}
\right) \quad (a, b\in \bfk).$$

\smallskip
\noindent
{\bf Case 2. Suppose $c_{33}\neq 0$:}
Then $c_{33}=2(c_{22}-c_{12})$. Thus  $$c_{33}=2(a-b)\quad \text{and} \quad c_{31}+c_{32}=-c_{33}=2(c_{12}-c_{22})
=2(b-a).$$
So $c_{32}=2(b-a)-c$. Thus we get the solutions
$$C_{22}=\left(
\begin{matrix}
a & -a & 0\\
b & -b & 0\\
c & 2(b-a)-c & 2(a-b)\\
\end{matrix}
\right)\quad (a, b \in \bfk, a\neq b).$$
They also satisfy the other equations in Eq.~(\mref{eq:3eijm}) hence give all the Rota-Baxter operators on $\bfk[CS(2)]$.

\subsubsection{The proof for $\bfk[CS(3)]$}
\mlabel{ss:cs3}
Here we prove that the matrices $C_P=(c_{ij})_{1\leq i, j\leq 3}$ of all the Rota-Baxter operators $P$ on the semigroup algebra $\bfk[CS(3)]$ are given by $C_{3,1}$ in Table~\mref{ta:rbcs}.

Applying the Cayley table of $CS(3)$ in Eqs.~(\mref{eq:3eijm}) and taking $i=1,j=1$ with $m=2$; $i=1,j=2$ with $m=1,3$; $i=2,j=2$ with $1\leq m\leq 3$ and $i=3,j=3$ with $1\leq m\leq 3$, we obtain
\begin{eqnarray}
&c_{12}^2=2c_{12}(c_{11}+c_{12}+c_{13}),
\mlabel{eq:cs332}\\
&c_{13}(c_{21}+
c_{22}+c_{23})+c_{12}c_{23}=
c_{11}(c_{11}+c_{13}),\mlabel{eq:cs334}\\
&c_{13}(c_{11}+c_{13}+c_{21}+c_{22}+c_{23}
)+c_{12}c_{23}=0,\mlabel{eq:cs336}\\
&(c_{21}+c_{23})(c_{21}+2c_{22}+c_{23})
=2(c_{11}(c_{21}+
c_{23})+c_{22}c_{21}),\mlabel{eq:cs3310}\\
&c_{22}^2=2(c_{12}(c_{21}+c_{23})+c_{22}^2)
,\mlabel{eq:cs3311}\\
&c_{13}(c_{21}+c_{23})+c_{22}c_{23}=0,
\mlabel{eq:cs3312}\\
&(c_{31}+c_{33})(c_{31}+c_{32}+c_{33})
+c_{32}(c_{31}+c_{33})=2c_{11}(c_{31}+
c_{32}+c_{33}),\mlabel{eq:cs3316}\\
&c_{32}^2=2c_{12}(c_{31}+c_{32}+c_{33}),
\mlabel{eq:cs3317}\\
&c_{13}(c_{31}+c_{32}+c_{33})=0.
\mlabel{eq:cs3318}
\end{eqnarray}
Eq.~(\mref{eq:cs336}) gives
$$c_{13}(c_{21}+c_{22}+c_{23})+c_{12}c_{23}=
-c_{13}(c_{11}+c_{13}).$$
Thus by Eq.~(\mref{eq:cs334}), we have
$(c_{11}+c_{13})^2=0,$ and so $c_{11}+c_{13}=0$. Then by Eq.~(\mref{eq:cs332}), we have $c_{12}=0$. Thus Eq.~(\mref{eq:cs3311}) gives $c_{22}=0$. By Eq.~(\mref{eq:cs3312}), we have $c_{13}(c_{21}+c_{23})=0$, and so $c_{11}(c_{21}+c_{23})=0$. Then Eq.~(\mref{eq:cs3310}) gives $c_{21}+c_{23}=0$. By Eq.~(\mref{eq:cs3317}), we have $c_{32}=0$. Then by Eqs.~(\mref{eq:cs3316}) and ~(\mref{eq:cs3318}), we can obtain $c_{31}+c_{33}=0$.
\delete{
Thus Eqs.~(\mref{eq:cs332})-(\mref{eq:cs3318})
are equivalent to
$$\left\{\begin{array}{cc}
c_{12}=c_{22}=c_{32}=0,\\
c_{11}+c_{13}=0,\\
c_{21}+c_{23}=0,\\
c_{31}+c_{33}=0.\\
\end{array}\right.$$
}
Let $a=c_{11}, b=c_{21}, c=c_{31}$. Thus solutions of Eqs.~(\mref{eq:cs332})-(\mref{eq:cs3318}) are given by
$$C_{3,1}=
\left(\begin{matrix}
a&0&-a\\
b&0&-b\\
c&0&-c\\
\end{matrix}\right)\quad (a,b,c\in \bfk).$$
It can be checked that they also satisfy the other equations in Eq.~(\mref{eq:3eijm}) and hence give all the Rota-Baxter operators on $\bfk[CS(3)]$.

\subsubsection{The proof for $\bfk[CS(4)]$}
\mlabel{ss:cs4}

We prove that the matrices $C_P=(c_{ij})_{1\leq i, j\leq 3}$ of all the Rota-Baxter operators $P$ on the semigroup algebra $\bfk[CS(4)]$ are given by $C_{4,1}$ in Table~\mref{ta:rbcs}.

Applying the Cayley table of $CS(4)$ in Eqs.~(\mref{eq:3eijm}) and
taking $i=j=1$ with $1\leq m\leq 3$; $i=j=2$ with $1\leq m\leq 3$ and $i=j=3$ with $1\leq m\leq 3$, we obtain
\begin{eqnarray}
&c_{11}^2=2c_{12}c_{13},
\mlabel{eq:cs341}\\
&c_{12}^2+2c_{12}(c_{11}+c_{13})=0,
\mlabel{eq:cs342}\\
&c_{13}^2+2c_{13}(c_{11}+c_{12})=0,
\mlabel{eq:cs343}\\
&c_{21}(c_{21}+c_{23})+c_{22}c_{23}+c_{23}
(c_{21}+c_{22})=2c_{11}(c_{21}+c_{23})
,\mlabel{eq:cs3410}\\
&c_{22}^2+2c_{12}(c_{21}+c_{23})=0,
\mlabel{eq:cs3411}\\
&c_{23}^2=2c_{13}(c_{21}+c_{23})+2c_{22}
c_{23},\mlabel{eq:cs3412}\\
&c_{31}(c_{31}+c_{32})+c_{32}
(c_{31}+c_{32})+c_{32}c_{33}
=2c_{11}(c_{31}+c_{32})
,\mlabel{eq:cs3416}\\
&c_{32}^2=2c_{12}(c_{31}+c_{32})+2c_{33}c_{32}
,\mlabel{eq:cs3417}\\
&c_{33}^2+2c_{13}(c_{31}+c_{32})=0.\mlabel{eq:cs3418}
\end{eqnarray}
By Eqs.~(\mref{eq:cs342}) and ~(\mref{eq:cs343}) and using Eq.(\mref{eq:cs341}), we have
$(c_{11}+c_{12})^2=0 $ and $(c_{11}+c_{13})^2=0$. So $c_{11}+c_{12}=0$ and $c_{11}+c_{13}=0$. Then Eq.~(\mref{eq:cs341}) gives $c_{11}^2=2c_{11}^2$. Thus $c_{11}=0$. Then we get $c_{11}=c_{12}=c_{13}=0$.
By Eq.~(\mref{eq:cs3411}), we have $c_{22}=0$. Thus Eqs.~(\mref{eq:cs3410})and ~(\mref{eq:cs3412}) give $c_{23}=0$ and $c_{21}=0$. Eq.~(\mref{eq:cs3418}) gives $c_{33}=0$ and so Eqs.~(\mref{eq:cs3416}) and ~(\mref{eq:cs3417}) give $c_{31}=c_{32}=0$. Thus the system in Eq.~(\mref{eq:3eijm}) only has the zero solution
$C_{4,1}=\mathbf{0}_{3\times 3}.$

\subsubsection{The proof for $\bfk[CS(5)]$}
\mlabel{ss:cs5}
We now prove that the matrices $C_P=(c_{ij})_{1\leq i, j\leq 3}$ of the Rota-Baxter operators $P$ on the semigroup algebra $\bfk[CS(5)]$ are given by $C_{5,1}$ in Table~\mref{ta:rbcs}.

Applying the Cayley table of $CS(5)$ in Eqs.~(\mref{eq:3eijm}) and taking $i=j=1$ with $1\leq m\leq 3$; $i=j=2$ with $1\leq m\leq 3$ and $i=j=3$ with $1\leq m\leq 3$, we get
\begin{eqnarray}
&c_{11}^2=0,\mlabel{eq:cs351}\\
&c_{13}^2=2c_{12}c_{11}+c_{12}^2,
\mlabel{eq:cs352}\\
&c_{13}(c_{11}+c_{12}+c_{13})=0,
\mlabel{eq:cs353}\\
&c_{21}^2=2c_{21}c_{11}
,\mlabel{eq:cs3510}\\
&(c_{22}+c_{23})^2=2c_{21}c_{12}+2c_{22}
(c_{22}+c_{23}),\mlabel{eq:cs3511}\\
&c_{21}c_{13}+c_{23}(c_{22}+c_{23})=0,
\mlabel{eq:cs3512}\\
&c_{31}(c_{31}+c_{32}+c_{33})+c_{31}(
c_{32}+c_{33})=2c_{11}c_{31}+2c_{21}(c_{
32}+c_{33}),\mlabel{eq:cs3516}\\
&(c_{32}+c_{33})^2=2 c_{12}c_{31}+2c_{22}
(c_{32}+c_{33}),\mlabel{eq:cs3517}\\
&c_{13}c_{31}+c_{23}(c_{32}+c_{33})=0.
\mlabel{eq:cs3518}
\end{eqnarray}
Note that Eq.~(\mref{eq:cs351}) gives $c_{11}=0$. Then by Eqs.~(\mref{eq:cs352}) and ~(\mref{eq:cs353}), we have $c_{12}^2=c_{13}^2$ and $c_{13}(c_{12}+c_{13})=0$. If $c_{12}+c_{13}\neq 0$, then $c_{13}=0$, and so $c_{12}=0$. This contradicts $c_{12}+c_{13}\neq 0$. Thus we get $c_{12}+c_{13}=0$. By Eq.~(\mref{eq:cs3510}), we have $c_{21}=0$. So Eqs.~(\mref{eq:cs3511}) and ~(\mref{eq:cs3512}) give $c_{22}^2=c_{23}^2$ and $c_{23}(c_{22}+c_{23})=0$. If $c_{22}+c_{23}\neq 0$, then $c_{23}=0$ and so $c_{22}=0$, contradiction. Thus $c_{22}+c_{23}=0$. By Eq.~(\mref{eq:cs3518}), we have $c_{12}c_{31}+c_{22}(c_{32}+c_{33}
)=0$.  Then by Eq.~(\mref{eq:cs3517}), $(c_{32}+c_{33})^2=0$. So $c_{32}+c_{33}=0$. Thus Eq.~(\mref{eq:cs3516}) gives $c_{31}=0$. Then we know that
Eqs.~(\mref{eq:cs351})-(\mref{eq:cs3518}) are equivalent to
$$
c_{11}=c_{21}=c_{31}=0,\
c_{12}+c_{13}=0,\
c_{22}+c_{23}=0,\
c_{32}+c_{33}=0.
$$
Let $a=c_{12}, b=c_{22}, c=c_{32}$. Then we get solutions of Eqs.~(\mref{eq:cs351})-(\mref{eq:cs3518})
$$C_{5,1}=\left(\begin{matrix}
0&a&-a\\
0&b&-b\\
0&c&-c\\
\end{matrix}\right).$$

\subsubsection{The proof for $\bfk[CS(6)]$}
\mlabel{ss:cs6}
We prove that the matrices $C_P=(c_{ij})_{1\leq i, j\leq 3}$ of the Rota-Baxter operators $P$ on the semigroup algebra $\bfk[CS(6)]$ are given by $C_{6,1}$ in Table~\mref{ta:rbcs}.

Applying the Cayley table of $CS(6)$ in Eqs.~(\mref{eq:3eijm}) and taking $i=j=1$ with $1\leq m\leq 3$; $i=j=2$ with $1\leq m\leq 3$ and $i=j=3$ with $1\leq m\leq 3$, we get
\begin{eqnarray}
&c_{11}^2=0,\mlabel{eq:cs361}\\
&c_{12}^2+2c_{11}c_{12}=0,\mlabel{eq:cs362}\\
&c_{13}^2+2c_{13}(c_{11}+c_{12})=0,
\mlabel{eq:cs363}\\
&c_{21}^2=2c_{11}c_{21},\mlabel{eq:cs3610}\\
&c_{22}^2+2c_{12}c_{21}=0,\mlabel{eq:cs3611}\\
&c_{23}^2+2c_{22}c_{23}+2c_{13}c_{21}=0,
\mlabel{eq:cs3612}\\
&c_{31}c_{32}+c_{31}(c_{31}+c_{32})=2c_{11}
c_{31}+2c_{21}c_{32},\mlabel{eq:cs3616}\\
&c_{32}^2=2c_{31}c_{12}+2c_{32}c_{22},\mlabel{eq:cs3617}\\
&c_{33}^2+2c_{13}c_{31}+2c_{23}c_{32}=0.\mlabel{eq:cs3618}
\end{eqnarray}
We note that Eqs.~(\mref{eq:cs361})-(\mref{eq:cs363})
give $c_{11}=c_{12}=c_{13}=0$. Then by Eqs.~(\mref{eq:cs3610})-(\mref{eq:cs3612}), we have $c_{21}=c_{22}=c_{23}=0$. Thus  Eqs.~(\mref{eq:cs3617}) and ~(\mref{eq:cs3618}) give $c_{32}=c_{33}=0$. Then Eq.~(\mref{eq:cs3616}) gives $c_{31}=0$.  Thus Eqs.~(\mref{eq:cs361})-(\mref{eq:cs3618}) only have the zero solution
$C_{6,1}=\mathbf{0}_{3\times 3}$,
as claimed.

\subsubsection{The proof for $\bfk[CS(7)]$}
\mlabel{ss:cs7}
Here we prove that the matrices $C_P=(c_{ij})_{1\leq i, j\leq 3}$ of the Rota-Baxter operators $P$ on the semigroup algebra $\bfk[CS(7)]$ are given by $C_{7,1}$ in Table~\mref{ta:rbcs}.

Applying the Cayley table of $CS(7)$ in Eqs.~(\mref{eq:3eijm}) and
taking $i=j=1$ with $1\leq m\leq 3$; $i=j=2$ with $1\leq m\leq 3$ and $i=j=3$ with $1\leq m\leq 3$, we obtain
\begin{eqnarray}
&c_{11}^2=c_{13}^2,\mlabel{eq:cs371}\\
&c_{12}^2+2c_{12}(c_{11}+c_{13})=0,
\mlabel{eq:cs372}\\
&c_{13}(c_{11}+c_{13})=0,
\mlabel{eq:cs373}\\
&c_{21}(c_{21}+2c_{23})+c_{23}^2=
2c_{21}c_{11}+2c_{23}c_{31},
\mlabel{eq:cs3710}\\
&c_{22}^2+2c_{21}c_{12}+2c_{23}c_{32}=0,
\mlabel{eq:cs3711}\\
&c_{21}c_{13}+c_{23}c_{33}=0,
\mlabel{eq:cs3712}\\
&c_{31}(c_{31}+2c_{33})+c_{33}^2=
2c_{11}(c_{31}+c_{33}),\mlabel{eq:cs3716}\\
&c_{32}^2+2c_{12}(c_{31}+c_{33})
=0,\mlabel{eq:cs3717}\\
&c_{13}(c_{31}+c_{33})=0.\mlabel{eq:cs3718}
\end{eqnarray}
Assume  $c_{11}+c_{13}\neq 0$. Then by Eq.~(\mref{eq:cs373}), we have $c_{13}=0$. So by Eq.~(\mref{eq:cs371}), we get $c_{11}=0$, a contradiction. This means that $c_{11}+c_{13}=0$. Thus Eq.~(\mref{eq:cs372}) gives $c_{12}=0$. Then by Eq.~(\mref{eq:cs3717}), we have $c_{32}=0$. So Eq.~(\mref{eq:cs3711}) gives $c_{22}=0$. By $c_{11}+c_{13}=0$ and Eq.~(\mref{eq:cs3718}), we have $c_{11}(c_{31}+c_{32})=0$. So by Eq.~(\mref{eq:cs3716}) we obtain $(c_{31}+c_{33})^2=0$. This means that $c_{31}+c_{33}=0$. Thus by Eq.~(\mref{eq:cs3712}), we get $c_{21}c_{11}+c_{23}c_{31}=0$ and then by Eq.~(\mref{eq:cs3710}), we get $c_{21}+c_{23}=0$. Thus Eqs.~(\mref{eq:cs371})-(\mref{eq:cs3718})
are equivalent to
$c_{12}=c_{22}=c_{32}=0$ together with
$$
c_{11}+c_{13}=0,\
c_{21}+c_{23}=0,\
c_{31}+c_{33}=0,\
c_{13}(c_{11}+c_{33})=0,\
c_{21}(c_{13}-c_{33})=0,\
c_{31}(c_{11}+c_{33})=0.
$$
We consider two cases.

\smallskip
\noindent
{\bf Case 1. $c_{11}=0$:} Then we have $c_{13}=0$ and $c_{31}=c_{33}=0$. Denote $b=c_{21}$. Then $c_{23}=-b$. Then we obtain the solutions
$$C_{7,1,1}=\left(\begin{matrix}
0&0&0\\
b&0&-b\\
0&0&0\\
\end{matrix}\right)\quad (b\in\bfk).$$
\smallskip

\noindent
{\bf Case 2. $c_{11}\neq 0$:} Then we have $c_{13}\neq 0$. So we have $c_{11}+c_{33}=0$. Thus $c_{11}=-c_{33}$. This means that $c_{33}\neq 0$. Then Eqs.~(\mref{eq:cs371})-(\mref{eq:cs3718})
are equivalent to
$$
c_{12}=c_{22}=c_{32}=0,\
c_{11}+c_{13}=0,\
c_{21}+c_{23}=0,\
c_{11}=c_{31},\
c_{13}=c_{33}.
$$
Denote $a=c_{11}, b=c_{21}$. Then we obtain the solutions of Eqs.~(\mref{eq:cs371})-(\mref{eq:cs3718}):
$$C_{7,1,2}=\left(\begin{matrix}
a&0&-a\\
b&0&-b\\
a&0&-a\\
\end{matrix}\right)\quad (a,b\in\bfk, a\neq 0).$$

In summary, the solutions of Eqs.~(\mref{eq:cs371})-(\mref{eq:cs3718}) are
$$C_{7,1}=\left(\begin{matrix}
a&0&-a\\
b&0&-b\\
a&0&-a\\
\end{matrix}\right)\quad (a ,b\in\bfk).$$

\subsubsection{The proof for $\bfk[CS(8)]$}
\mlabel{ss:cs8}
We prove that the matrices $C_P=(c_{ij})_{1\leq i, j\leq 3}$ of the Rota-Baxter operators $P$ on the semigroup algebra $\bfk[CS(8)]$ are given by $C_{8,1}$ in Table~\mref{ta:rbcs}.

Applying the Cayley table of $CS(8)$ in Eqs.~(\mref{eq:3eijm}) and
taking $i=j=1$ with $1\leq m\leq 3$, $i=j=2$ with $1\leq m\leq 3$; $i=2,j=3$ with $m=3$ and $i=j=3$ with $1\leq m\leq 2$, we obtain
\begin{eqnarray}
&c_{11}^2=0,\mlabel{eq:cs381}\\
&c_{13}^2=c_{12}^2+2c_{12}(c_{11}+c_{13}),\mlabel{eq:cs382}\\
&c_{13}(c_{11}+c_{13})=0,\mlabel{eq:cs383}\\
&c_{21}(c_{21}+c_{23})+c_{21}c_{23}=
2(c_{21}c_{11}+c_{23}c_{31}),\mlabel{eq:cs3810}\\
&c_{23}^2=c_{22}^2+2(c_{21}c_{12}+c_{23}
c_{32}),\mlabel{eq:cs3811}\\
&c_{21}c_{13}+c_{23}c_{33}=0,\mlabel{eq:cs3812}\\
&c_{13}(c_{21}+c_{31})+c_{23}^2+c_{33}^2=0,
\mlabel{eq:cs3815}\\
&c_{31}(c_{31}+2c_{33})=2c_{31}c_{11}+2
c_{33}c_{21},\mlabel{eq:cs3816}\\
&c_{33}^2=2c_{31}c_{12}+2c_{33}c_{22}
+c_{32}^2.\mlabel{eq:cs3817}
\end{eqnarray}
Note that Eq.~(\mref{eq:cs381}) gives $c_{11}=0$.  By Eqs.~(\mref{eq:cs382}) and ~(\mref{eq:cs383}), we have $c_{12}=c_{13}=0$. Thus by Eqs.~(\mref{eq:cs3812}) and ~(\mref{eq:cs3815}), we have $c_{23}c_{33}=0$ and $c_{23}^2+c_{33}^2=0$. Then we get $c_{23}=c_{33}=0$. So by Eqs.~(\mref{eq:cs3810}) and ~(\mref{eq:cs3811}), we have $c_{21}=c_{22}=0$. Thus Eqs.~(\mref{eq:cs3816}) and~(\mref{eq:cs3817}) give $c_{31}=c_{32}=0$. So the only solution is the zero solution giving by the zero matrix
$C_{8,1}=\mathbf{0}_{3\times 3},$
proving the claim.

\subsubsection{The proof for $\bfk[CS(9)]$}
\mlabel{ss:cs9}

We prove that the matrices $C_P=(c_{ij})_{1\leq i, j\leq 3}$ of the Rota-Baxter operators $P$ on the semigroup algebra $\bfk[CS(9)]$ are given by $C_{9,1}$ in Table~\mref{ta:rbcs}.

Applying the Cayley table of $CS(9)$ in Eqs.~(\mref{eq:3eijm}) and
taking $i=j=1$ with $1\leq m\leq 3$; $i=1,j=3$ with $m=3$; $i=j=2$ with $1\leq m\leq 3$ and $i=j=3$ with $1\leq m\leq 2$, we obtain
\begin{eqnarray}
&c_{12}^2+c_{13}^2=c_{11}^2+2c_{13}c_{31},\mlabel{eq:cs391}\\
&c_{13}c_{33}=0,\mlabel{eq:cs392}\\
&c_{12}(c_{11}+c_{12})+ c_{13}c_{32}=0,\mlabel{eq:cs393}\\
&c_{13}^2+c_{33}^2=0,\mlabel{eq:cs399}\\
&(c_{21}+c_{22})^2+c_{23}^2=2c_{11}
(c_{21}+
c_{22})+2c_{23}c_{31},\mlabel{eq:cs3910}\\
&c_{12}(c_{21}+c_{22})+c_{23}c_{32}=0,
\mlabel{eq:cs3911}\\
&c_{23}(c_{21}+c_{22})=c_{13}(c_{21}+c_{
22})+c_{23}c_{33},\mlabel{eq:cs3912}\\
&c_{32}^2+c_{33}^2
=2c_{33}c_{11}+c_{31}^2,\mlabel{eq:cs3916}\\
&c_{33}c_{12}+c_{32}(c_{31}+c_{32})=0.
\mlabel{eq:cs3917}
\end{eqnarray}
Note that Eqs.~(\mref{eq:cs392}) and ~(\mref{eq:cs399}) give $c_{13}=c_{33}=0$. Then by Eqs.~(\mref{eq:cs391}) and ~(\mref{eq:cs393}), we obtain $c_{11}+c_{12}=0$. Similarly, Eqs.~(\mref{eq:cs3916}) and ~(\mref{eq:cs3917}) give $c_{31}+c_{32}=0$. By Eq.~(\mref{eq:cs3912}), we get $c_{23}(c_{21}+c_{22})=0$. Assume $c_{21}+c_{22}\neq 0$. Then we have $c_{23}=0$. So by Eq.~(\mref{eq:cs3911}), $c_{11}=c_{12}=0$. Thus Eq.~(\mref{eq:cs3910}) gives $c_{21}+c_{22}=0$, a contradiction. So we must have $c_{21}+c_{22}=0$. Then by Eq.~(\mref{eq:cs3911}), we have $c_{23}c_{32}=0$. So $c_{23}c_{31}=0$. Thus Eq.~(\mref{eq:cs3910}) gives $c_{23}=0$. Then Eqs.~(\mref{eq:cs391})-(\mref{eq:cs3917}) are equivalent to
$$
c_{13}=c_{23}=c_{33}=0,\
c_{11}+c_{12}=0,\
c_{21}+c_{22}=0,\
c_{31}+c_{32}=0.
$$
Denote $a=c_{11}, b=c_{21}, c=c_{31}$. Then we get the solutions of Eqs.~(\mref{eq:cs391})-(\mref{eq:cs3917}) are
$$C_{9,1}=\left(\begin{matrix}
a&-a&0\\
b&-b&0\\
c&-c&0\\
\end{matrix}\right)\quad (a,b,c\in\bfk).
$$

\subsubsection{The proof for $\bfk[CS(10)]$}
\mlabel{ss:cs10}
We prove that the matrices $C_P=(c_{ij})_{1\leq i, j\leq 3}$ of the Rota-Baxter operators $P$ on the semigroup algebra $\bfk[CS(10)]$ are given by $C_{10,1}$ in Table~\mref{ta:rbcs}.

Applying the Cayley table of $CS(10)$ in Eqs.~(\mref{eq:3eijm}) and
taking $i=j=1$ with $1\leq m\leq 3$; $i=1,j=3$ with $m=3$; $i=j=2$ with $1\leq m\leq 3$ and $i=j=3$ with $1\leq m\leq 2$, we obtain
\begin{eqnarray}
&c_{13}^2=c_{11}^2+2c_{13}c_{31},
\mlabel{eq:cs3101}\\
&c_{12}^2+2c_{11}c_{12}+2c_{13}c_{32}
=0,\mlabel{eq:cs3102}\\
&c_{13}c_{33}=0,\mlabel{eq:cs3103}\\
&c_{13}^2+c_{33}^2=0,\mlabel{eq:cs3109}\\
&c_{21}^2+c_{23}^2=2c_{21}c_{11}+2c_{23}
c_{31},\mlabel{eq:cs31010}\\
&c_{22}^2+2c_{21}c_{12}+2c_{23}c_{32}
=0,
\mlabel{eq:cs31011}\\
&c_{21}c_{23}=c_{21}c_{13}+c_{23}c_{33},\mlabel{eq:cs31012}\\
&c_{33}^2=2c_{33}c_{11}+c_{31}^2,
\mlabel{eq:cs31016}\\
&c_{32}^2+2c_{33}c_{12}+2c_{32}c_{31}=0.
\mlabel{eq:cs31017}
\end{eqnarray}
First Eqs.~(\mref{eq:cs3103}) and ~(\mref{eq:cs3109}) give $c_{13}=c_{33}=0$. Thus by Eqs.~(\mref{eq:cs3101}) and ~(\mref{eq:cs3102}), we have $c_{11}=c_{12}=0.$  By Eqs.~(\mref{eq:cs31016}) and ~(\mref{eq:cs31017}), we get $c_{31}=c_{32}=0$. Then by Eqs.~(\mref{eq:cs31010}) and~(\mref{eq:cs31012}), we have $c_{21}^2
+c_{23}^2=0$ and $c_{21}c_{23}=0$. Thus we have $c_{21}=c_{23}=0$. Further, by Eq.~(\mref{eq:cs31011}), we have $c_{22}=0$. We see that the system only has the zero solution
given by the zero matrix $C_{10,1}=\mathbf{0}_{3\times 3}.$

\subsubsection{The proof for $\bfk[CS(11)]$}
\mlabel{ss:cs11}
We prove that the matrices $C_P=(c_{ij})_{1\leq i, j\leq 3}$ of the Rota-Baxter operators $P$ on the semigroup algebra $\bfk[CS(11)]$ are given by $C_{11,1}$ in Table~\mref{ta:rbcs}.

Applying the Cayley table of $CS(11)$ in Eqs.~(\mref{eq:3eijm}) and
taking $i=j=1$ with $1\leq m\leq 3$; $i=1,j=2$ with $2\leq m\leq 3$; $i=j=2$ with $1\leq m\leq 2$ and $i=j=3$ with $1\leq m\leq 3$, we obtain
\begin{eqnarray}
&(c_{12}+c_{13})^2=c_{11}^2+2c_{21}
(c_{12}+c_{13}),\mlabel{eq:cs3111}\\
&c_{11}c_{13}=c_{22}(c_{12}+c_{13}),
\mlabel{eq:cs3112}\\
&c_{11}c_{13}+c_{23}(c_{12}+c_{13})=0,
\mlabel{eq:cs3113}\\
&c_{11}c_{23}+c_{21}c_{13}=c_{12}(c_{12}
+c_{13})+c_{22}(c_{22}+c_{23}),
\mlabel{eq:cs3115}\\
&c_{13}(c_{12}+c_{13}+c_{21})+
c_{23}(c_{11}+c_{22}+c_{23})=0,
\mlabel{eq:cs3116}\\
&(c_{22}+c_{23})^2=2c_{11}(c_{22}+c_{23})
+c_{21}^2,\mlabel{eq:cs31110}\\
&c_{21}c_{23}=c_{12}(c_{22}+c_{23}),
\mlabel{eq:cs31111}\\
&c_{31}^2+(c_{32}+c_{33})^2=2c_{11}(c_{32}
+c_{33})+2c_{31}c_{21},\mlabel{eq:cs31116}\\
&c_{31}(c_{32}+c_{33})=c_{12}(c_{32}+c_{33}
)+c_{31}c_{22},\mlabel{eq:cs31117}\\
&c_{13}(c_{32}+c_{33})+c_{31}c_{23}=0.
\mlabel{eq:cs31118}
\end{eqnarray}
By Eq.~(\mref{eq:cs3116}), we get $$c_{13}c_{21}+c_{23}c_{11}=-c_{13}(c_{12}
+c_{13})-c_{23}(c_{22}+c_{23}).$$
Then by Eq.~(\mref{eq:cs3115}), we get $(c_{12}+c_{13})^2+(c_{22}+c_{23})^2=0$.
Eqs.~(\mref{eq:cs3112}) and ~(\mref{eq:cs3113}) give $(c_{22}+c_{23})(c_{12}+c_{13})=0$. Thus we get $c_{22}+c_{23}=0$ and $c_{12}+c_{13}=0$. Then by Eqs.~(\mref{eq:cs3111}) and~(\mref{eq:cs31110}), we get $c_{11}=c_{21}=0$. By Eq.~(\mref{eq:cs31118}), we have $c_{12}(c_{32}+c_{33})+c_{31}c_{22}=0$.
So by Eq.~(\mref{eq:cs31117}), we have $c_{31}(c_{32}+c_{33})=0$. With Eq.~(\mref{eq:cs31116}) and $c_{11}=c_{21}=0$, we get $c_{31}=0$ and $c_{32}+c_{33}=0$.
Thus Eqs.~(\mref{eq:cs3111})-(\mref{eq:cs31118})
are equivalent to
$$
c_{11}=c_{21}=c_{31}=0,\
c_{12}+c_{13}=0,\
c_{22}+c_{23}=0,\
c_{32}+c_{33}=0.
$$
Let $a=c_{12}, b=c_{22}$ and $c= c_{32}$. Then we obtain the solutions of Eqs.~(\mref{eq:cs3111})-(\mref{eq:cs31118})
$$C_{11,1}=\left(\begin{matrix}
0&a&-a\\
0&b&-b\\
0&c&-c\\
\end{matrix}\right)\quad (a,b,c\in\bfk).$$

\subsubsection{The proof for $\bfk[CS(12)]$}
\mlabel{ss:cs12}
We finally prove that the matrices $C_P=(c_{ij})_{1\leq i, j\leq 3}$ of the Rota-Baxter operators $P$ on the semigroup algebra $\bfk[CS(12)]$ are given by $C_{12,1}$ in Table~\mref{ta:rbcs}.

Applying the Cayley table of $CS(12)$ in Eqs.~(\mref{eq:3eijm}) and
taking $i=j=1$ with $1\leq m\leq 3$; $i=1,j=2$ with $m=2,3$; $i=1,j=3$ with $m=2,3$; $i=j=2$ with $1\leq m\leq 3$; $i=2,j=3$ with $m=2,3$ and $i=j=3$ with $1\leq m\leq 3$, we obtain
\begin{eqnarray}
&2c_{12}c_{13}=c_{11}^2+2c_{12}
c_{21}+2c_{13}c_{31},\mlabel{eq:cs3121}\\
&c_{13}^2=2
c_{12}c_{22}+2c_{13}c_{32},\mlabel{eq:cs3122}\\
&c_{12}^2=2
c_{12}c_{23}+2c_{13}c_{33},\mlabel{eq:cs3123}\\
&c_{13}c_{23}=c_{12}c_{13}+c_{22}^2+c_{32}
(c_{12}+c_{23}),\mlabel{eq:cs3125}\\
&c_{12}c_{22}=c_{13}^2+c_{23}c_{22}+c_{33}
(c_{12}+c_{23}),\mlabel{eq:cs3126}\\
&c_{13}c_{33}=c_{12}^2+c_{22}(c_{13}+c_{32}
)+c_{32}c_{33},\mlabel{eq:cs3128}\\
&c_{12}c_{32}=c_{12}c_{13}+c_{23}(c_{13}+
c_{32})+c_{33}^2,\mlabel{eq:cs3129}\\
&2c_{22}c_{23}=2c_{23}c_{11}+c_{21}^2+2c_{22}
c_{31},\mlabel{eq:cs31210}\\
&c_{23}^2=2c_{23}c_{12}+2c_{22}c_{32},\mlabel{eq:cs31211}\\
&c_{22}^2=2c_{23}c_{13}+2c_{22}c_{33},\mlabel{eq:cs31212}\\
&c_{23}c_{33}=c_{12}(c_{22}+c_{33})+
c_{22}c_{23}+c_{32}^2,\mlabel{eq:cs31214}\\
&c_{22}c_{32}=c_{13}(c_{22}+c_{33})
+c_{23}^2+c_{33}c_{32},\mlabel{eq:cs31215}\\
&2c_{32}c_{33}=2c_{11}c_{32}+2c_{33}c_{21}
+c_{31}^2,\mlabel{eq:cs31216}\\
&c_{33}^2=2c_{12}c_{32}+2c_{33}c_{22},\mlabel{eq:cs31217}\\
&c_{32}^2=2c_{13}c_{32}+2c_{33}c_{23}.\mlabel{eq:cs31218}
\end{eqnarray}
Denote $a=c_{13}, b=c_{23}$ and $c=c_{33}$.
By Eqs.~(\mref{eq:cs31211}) and ~(\mref{eq:cs31217}), we get
$$\left\{\begin{array}{cc}
2c_{12}c_{22}c_{32}=b^2c_{12}-2bc_{12}^2,\\
2c_{12}c_{22}c_{32}=c^2c_{22}-2c c_{22}^2.\\
\end{array}\right.$$
So
\begin{equation}
2bc_{12}^2 - b^2c_{12}=2cc_{22}^2- c^2c_{22}.
\mlabel{eq:cs312im}
\end{equation}
By Eqs.~(\mref{eq:cs3123}) and ~(\mref{eq:cs31212}), we have
$$\left\{
\begin{array}{cc}
4abc=2bc_{12}^2-4b^2c_{12}^2,\\
4abc=2cc_{22}^2-4c^2c_{22}.\\
\end{array}\right.$$
Thus
$$\left\{
\begin{array}{cc}
2bc_{12}^2-b^2c_{12}=4abc+3b^2c_{12},\\
2cc_{22}^2-c^2c_{22}=4abc+3c^2c_{22}.\\
\end{array}\right.$$
Using Eq.~(\mref{eq:cs312im}), we obtain  $b^2c_{12}=c^2c_{22}$, and so $ bc_{12}^2=cc_{22}^2$. Similarly, by Eqs.~(\mref{eq:cs3122}) and ~(\mref{eq:cs31211}), and Eqs.~(\mref{eq:cs3123}) and ~(\mref{eq:cs31218}), we get $b^2c_{12}=a^2c_{32}$ and $bc_{12}^2=ac_{32}^2$, respectively. This means that
\begin{equation}
b^2c_{12}=c^2c_{22}=a^2c_{32}
\mlabel{eq:cs312im2}
\end{equation}
and
\begin{equation} bc_{12}^2=cc_{22}^2=ac_{32}^2.
\mlabel{eq:cs312im3}
\end{equation}
We consider two cases.

\smallskip
\noindent
{\bf Case 1. abc=0:} There are three subcases to consider.

\smallskip
\noindent
{\bf Subcase 1. a=b=c=0:} Then by  Eqs.~(\mref{eq:cs3123}),~
(\mref{eq:cs31214}) and ~(\mref{eq:cs3125}), we have
$c_{12}=0,c_{32}=0$ and $c_{22}=0$. By Eqs.~(\mref{eq:cs3121}), ~(\mref{eq:cs31210}) and~(\mref{eq:cs31216}), we obtain $c_{11}=c_{21}=c_{31}=0$. Thus the only solution is
given by the zero matrix $C_{12,1}:=\mathbf{0}_{3\times 3}$.

\smallskip
\noindent
{\bf Subcase 2. Two of $a,b,c$ are $0$:} Without loss of generality, we may assume  $a=b=0$ and $c\neq 0$. From Eq.~(\mref{eq:cs31218}), we obtain $c_{32}=0$. Then by Eq.~(\mref{eq:cs3129}), we have $c=0$, a contradiction. Thus this case can not occur.

\smallskip
\noindent
{\bf Subcase 3. One of $a,b, c$ is $0$:} We may assume without loss of generality that $a=0$, $b\neq0$ and $c\neq 0$. By Eq.~(\mref{eq:cs3122}), we have $c_{12}c_{22}=0$. Then $c_{12}=0$ or $c_{22}=0$. If $c_{12}=0$, then by Eq.~(\mref{eq:cs3126}), we have $c_{33}+c_{22}=0$, and so $c_{22}=-c\neq 0$. By Eq.~(\mref{eq:cs31212}), $3c^2=0$. Thus $c=0$. This contradicts $c\neq 0$. If $c_{22}=0$, then by Eq.~(\mref{eq:cs3126}), $b=-c_{12}$. By Eq.~(\mref{eq:cs3123}), we have $3b^2=0$. Thus $b=0$, a contradiction.
We see that this case also can not occur.

\smallskip
\noindent
{\bf Case 2. abc $\neq 0$:} Then $a\neq 0,b\neq 0$ and $c\neq 0.$ If one of $c_{12}, c_{22}, c_{32}$ is $0$, we may assume $c_{12}= 0$, and  then
by Eqs.~(\mref{eq:cs312im2} and~(\mref{eq:cs312im3}), we obtain $c_{22}=c_{32}=0$. So Eq.~(\mref{eq:cs3122}) gives $a^2=0$, a contradiction. Thus $c_{12}, c_{22}$ and $c_{32}$ are nonzero. Since $b^2c_{12}=c^2c_{22}$ and $bc_{12}^2=cc_{22}^2$, we have
$ b^2c_{12}^2=bcc_{22}^2.$ So $c^2c_{22}c_{12}=bcc_{22}^2$. Thus $cc_{12}=bc_{22}$. Similarly, we have $ac_{12}=bc_{32}$ and $ac_{22}=cc_{32}.$ So $b^3c_{12}=c^2bc_{22}=c^3c_{12}$. Thus $a^3=c^3$. Similarly, we have $a^3=b^3$. By Eqs.~(\mref{eq:cs31211}) and ~(\mref{eq:cs31215}),  we get
$$b^2+2bc_{12}+2ac_{22}+2cc_{32}+2ac=0.$$
Since $ac_{22}=cc_{32}$, we have $b^2+2bc_{12}+4ac_{22}+2ac=0$. So $b^3+2b^2c_{12}+4abc_{22}+2abc=0$. Since $b^2c_{12}=c^2c_{22}$, we have
\begin{equation}
2(c^2+2ab)c_{22}=-b^3-2abc=-c^3-2abc.
\mlabel{eq:cs312im4}
\end{equation}
We further divide into two subcases.

\smallskip
\noindent
{\bf Subcase 1. $c^2+2ab\neq 0$:} By Eq.~(\mref{eq:cs312im4}), we have $c_{22}=-\frac{c}{2}$. From $cc_{12}=
bc_{22}$ and $ac_{22}=cr_{32}$, we obtain $c_{12}=-\frac{b}{2}$ and $c_{32}=-\frac{a}{2}$. By Eq.~(\mref{eq:cs3122}), we have $a^2=\frac{bc}{2}-a^2$. So $a^2=\frac{bc}{4}$. By Eq.~(\mref{eq:cs31218}), we get $\frac{5}{4}a^2=2bc$. Thus we get $a^2=\frac{8}{5}bc$. So by $a^2=\frac{bc}{4}$, we get $bc=0$, a contradiction.

\smallskip
\noindent
{\bf Subcase 2. $c^2+2ab=0$:} Then $c^4=4a^2b^2$. Adding Eqs.~(\mref{eq:cs31212}) and ~(\mref{eq:cs31217}), we obtain $c_{22}^2+2c_{12}c_{32}=0$. So we have $a^2b^2c_{22}^2+2a^2b^2c_{12}c_{32}=0$.
By Eq.~(\mref{eq:cs312im2}), we have $(a^2b^2+2c^4)c_{22}^2=0$. Since $c_{22}\neq 0$, we have $a^2b^2+2c^4=0$. Then by $c^4=4a^2b^2$, we have $9(ab)^2=0$, and so $ab=0$, again a contradiction.
\smallskip

In summary, the only solution for $\bfk[CS(12)]$ is the zero solution
$\mathbf{0}_{3\times 3}$, as claimed.
\smallskip

Now the proof of Theorem~\mref{thm:rboc3} is completed.

\section{Rota-Baxter operators on noncommutative semigroup algebras of order 3}
\mlabel{sec:rboncs}

In this Section, we classify all Rota-Baxter operators on noncommutative
semigroup algebras of order $3$.

\subsection{Statement of the classification theorem in the noncommutative case}
\mlabel{ss:stnc}

A classification of the six noncommutative semigroups of order $3$, up to
isomorphism and anti-isomorphism, is given in Table~\mref{ta:ncs3}.

\begin{table}[!tpbh]
\caption{Cayley table of noncommutative semigroups of order 3 \mlabel{ta:ncs3}}
\begin{minipage}[tbph]{0.75\textwidth}
\tiny
\begin{tabular}{|c|c|c|}
\hline
$NCS(1):=$
\begin{tabular}{c|ccccc}
$\cdot$ &$e_1$ & $e_2$& $e_3$\\
\hline
$e_1$ &$e_1$ & $e_1$ &$e_1$\\
$e_2$ &$e_1$ & $e_2$ &$e_1$\\
$e_3$ &$e_1$ & $e_3$ &$e_1$\\
\end{tabular}
&$NCS(2):=$
\begin{tabular}{c|ccccc}
$\cdot$ &$e_1$ & $e_2$& $e_3$\\
\hline
$e_1$ &$e_1$ & $e_1$ &$e_1$\\
$e_2$ &$e_1$ & $e_2$ &$e_1$\\
$e_3$ &$e_3$ & $e_3$ &$e_3$\\
\end{tabular}
&$NCS(3):=$
\begin{tabular}{c|ccccc}
$\cdot$ &$e_1$ & $e_2$& $e_3$\\
\hline
$e_1$ &$e_1$ & $e_1$&$e_1$\\
$e_2$ &$e_1$ & $e_2$ &$e_2$\\
$e_3$ &$e_1$ & $e_3$ &$e_3$\\
\end{tabular}
\\
\hline
$NCS(4):=$
\begin{tabular}{c|ccccc}
$\cdot$ &$e_1$ & $e_2$& $e_3$\\
\hline
$e_1$ &$e_1$ & $e_1$ &$e_1$\\
$e_2$ &$e_2$ & $e_2$ &$e_2$\\
$e_3$&$e_1$ & $e_1$ &$e_1$\\
\end{tabular}
&
$NCS(5):=$
\begin{tabular}{c|ccccc}
$\cdot$ &$e_1$ & $e_2$& $e_3$\\
\hline
$e_1$ &$e_1$ & $e_1$ &$e_1$\\
$e_2$ &$e_2$ & $e_2$ &$e_2$\\
$e_3$ &$e_3$ & $e_3$ &$e_3$\\
\end{tabular}
&$NCS(6)$:=\begin{tabular}{c|cccc}
$\cdot$ &$e_1$ & $e_2$& $e_3$\\
\hline
$e_1$ &$e_1$ & $e_1$ &$e_1$\\
$e_2$ &$e_1$ & $e_2$ &$e_3$\\
$e_3$ &$e_3$ & $e_3$ &$e_3$\\
\end{tabular}
\\
\hline
\end{tabular}
\end{minipage}
\end{table}
\vspace{.5cm}

For Rota-Baxter operators on the corresponding semigroup algebras, we have the following classification theorem whose proof will be given in Section~\mref{ss:proofnc}
\begin{theorem}
Let $\bfk$ be a field of characteristic zero that is closed under taking square root. The matrices of the Rota-Baxter operators
on noncommutative semigroup algebras of order three are given in
Table~\mref{ta:rbncs1}, where all the parameters take values in~$\bfk$ and~$i$
denotes~$\sqrt{-1}$ as usual.
\mlabel{thm:rbonc3}
\end{theorem}

\begin{tiny}
\begin{longtable}{|p{1.5cm}|p{11.3cm}|}
\caption[RBOs on 3-dimensional noncomumutative semigroup algebras]{{RBOs on noncomumutative semigroup algebras of order 3} \mlabel{ta:rbncs1}}\\
\hline \multicolumn{1}{|c|}{{\bf Semigroups}} & \multicolumn{1}{c|}{{\bf Matrices of Rota-Baxter operators on semigroup algebras}}\\
\hline
\endfirsthead

\multicolumn{2}{c}%
{{ \tablename\  \thetable{:\ Rota-Baxter operators on noncomumutative semigroup algebras of order 3}}} \\
\hline \multicolumn{1}{|c|}{{\bf Semigroups }} &
\multicolumn{1}{c|}{{\bf Matrices of Rota-Baxter operators on semigroup algebras}}\\
\hline
\endhead

\hline \multicolumn{2}{|r|}{{Continued on next page}} \\ \hline\endfoot

\hline
\endlastfoot

$NCS(1)$
& \begin{minipage}[t]{0.85\textwidth}
\noindent
$
N_{1,1}=\left(\begin{matrix}
0 & 0 & 0\\
0 &0  & 0\\
0 & a & 0\\
\end{matrix}
\right),
N_{1,2}=\left(\begin{matrix}
0 & 0 & 0\\
0 &0  & 0\\
-a & a & 0\\
\end{matrix}\right)(a\neq0),
N_{1,3}=\left(\begin{matrix}
0 & 0 & 0\\
-b+c &-c  & b\\
\frac{c^2-bc}{b} & -\frac{c^2}{b} & c\\
\end{matrix}\right)(a\neq0, b\neq0),\\
N_{1,4}=\left(\begin{matrix}
0 & 0 & 0\\
-a &-b  & a\\
-b & -\frac{b^2}{a} & b\\
\end{matrix}
\right)\,(a\neq0),
N_{1,5}=\left(\begin{matrix}
a & 0 & -a\\
-b &0  & b\\
a & 0 & -a\\
\end{matrix}\right)\,(a\neq 0)
$
\end{minipage}\\
\hline
$NCS(2)$
& \begin{minipage}[t]{0.85\textwidth}
\noindent
$N_{2,1}=\left(
\begin{matrix}
0 & 0 & 0\\
c & 0 & -c\\
0 & 0 & 0\\
\end{matrix}\right)(a\neq 0),
N_{2,2}=\left(\begin{matrix}
0 & 0 & 0\\
0 & 0 & 0\\
a & 0 & 0\\
\end{matrix}\right),
N_{2,3}=\left(\begin{matrix}
0 & 0 & 0\\
0 & 0 & 0\\
0 & b & 0\\
\end{matrix}\right),\\
N_{2,4}=\left(\begin{matrix}
b & 0 & 0\\
b & 0 & 0\\
-\frac{b^2}{a} & 0 & -b\\
\end{matrix}\right)(a\neq 0, b\neq 0, a+b\neq0),
N_{2,5}=\left(
\begin{matrix}
0 & 0 & a\\
0& 0 & 0\\
0 & 0 & 0\\
\end{matrix}\right)(a\neq 0),\\
N_{2,6}=\left(
\begin{matrix}
0 & 0 & a\\
0 & 0 & a\\
0 & 0 & 0\\
\end{matrix}\right)(a\neq 0),
N_{2,7}=\left(\begin{matrix}
-a & 0 & a\\
c &0  & -c\\
-a & 0 & a\\
\end{matrix}\right)(a\neq0),\\
N_{2,8}=\left(\begin{matrix}
-a & b & a\\
-a & b & a\\
b-a & \frac{b(a-b)}{a} & a-b\\
\end{matrix}\right)(a\neq 0,b\neq 0),
N_{2,9}=\left(\begin{matrix}
0 & b & 0\\
 0 & 0 & 0\\
0 & 0 & 0\\
\end{matrix}\right)
(b\neq0)\\
$
\end{minipage}\\
\hline
$NCS(3)$
&\begin{minipage}[t]{0.85\textwidth}
\noindent
$N_{3,1}=\left(\begin{matrix}
0 & a & -a\\
0 & b & -b\\
0 & b & -b\\
\end{matrix}\right)(a\neq 0),
N_{3,2}=\left(\begin{matrix}
0 & 0 & 0\\
0 & -\sqrt{ab}i & b\\
0 & a & \sqrt{ab}i\\
\end{matrix}\right)(a\neq 0,b\neq 0),\\
N_{3,3}=\left(\begin{matrix}
0 & 0 & 0\\
0 & \sqrt{ab}i & b\\
0 & a & -\sqrt{ab}i\\
\end{matrix}\right)(a\neq 0,b\neq 0),
N_{3,4}=\left(\begin{matrix}
0&0&0\\
0&0&0\\
0&a&0\\
\end{matrix}\right)(a\neq 0),\\
N_{3,5}=\left(\begin{matrix}
0&0&0\\
0&0&0\\
-a&a&0\\
\end{matrix}\right)(a\neq 0),
N_{3,6}=\left(\begin{matrix}
0&0&0\\
c&-(b+c)&b\\
\frac{c(b+c)}{b}&-\frac{(b+c)^2}{b}&b+c\\
\end{matrix}\right)(a\neq 0, b\neq 0, b+c\neq 0),\\
N_{3,7}=\left(\begin{matrix}
0&0&0\\
0&0&b\\
0&0&0\\
\end{matrix}\right),
N_{3,8}=\left(\begin{matrix}
0&0&0\\
-b&0&b\\
0&0&0\\
\end{matrix}\right)(b\neq 0)\\
$
\end{minipage}\\
\hline
$NCS(4)$
& \begin{minipage}[t]{0.85\textwidth}
\noindent
$N_{4,1}=\left(\begin{matrix}
-a & 0 & a\\
b & 0 & -b\\
c  & 0 & -c\\
\end{matrix}\right)(a\neq 0, a+b\neq 0),
N_{4,2}=\left(\begin{matrix}
-a & 0 & a\\
-a & 0 & a\\
c  & d & -c-d\\
\end{matrix}\right)(a\neq 0),\\
N_{4,3}=\left(\begin{matrix}
0 & 0 & 0\\
-a & 0 & a\\
b  & 0 & -b\\
\end{matrix}\right)(a\neq 0),
N_{4,4}=\left(\begin{matrix}
0 & 0 & 0\\
0 & 0 & 0\\
b  & c & -(b+c)\\
\end{matrix}\right)(c\neq 0),\\
N_{4,5}=\left(\begin{matrix}
0&0&0\\
a&0&0\\
b&0&-b\\
\end{matrix}\right),
N_{4,6}=\left(\begin{matrix}
 a&-(a+b)& b\\
 a&-(a+b)&b\\
 c&d&-(c+d)\\
 \end{matrix}\right)(b\neq 0, a+b\neq 0),\\
N_{4,7}=\left(
 \begin{matrix}
 a&b&0\\
 -\frac{a^2}{b}&-a&0\\
 c&b&a-c\\
 \end{matrix}\right)(a+b\neq 0, b\neq 0),
N_{4,8}=\left(\begin{matrix}
 -b&b&0\\
 -b&b&0\\
 c&d&-c-d\\
 \end{matrix}\right)(b\neq 0)\\
 $
\end{minipage}\\
\hline
$NCS(5)$
& \begin{minipage}[t]{0.85\textwidth}
\noindent
$N_{5,1}=\left(
\begin{matrix}
\frac{ac}{F} & \frac{ab}{F} & a\\
\frac{cd}{F} & \frac{bd}{F} & d\\
c & b & F\\
\end{matrix}
\right),
N_{5,2}=\left(
\begin{matrix}
-\frac{ac}{F} & -\frac{ab}{F} & a\\
-\frac{cd}{F} & -\frac{bd}{F} & d\\
c & b & -F\\
\end{matrix}
\right)(a,b,c,d\in\bfk\setminus\{0\}), F:=\sqrt{-(ac+bd)},\\
N_{5,3}=\left(\begin{matrix}
0 &  0 & 0\\
-\frac{c\sqrt{d}}{\sqrt{b}}i & -\sqrt{b d}i & d\\
c &b & \sqrt{bd}i\\
\end{matrix}\right),
N_{5,4}=\left(\begin{matrix}
0 & 0 & 0\\
\frac{c\sqrt{d}}{\sqrt{b}}i & \sqrt{b d}i & d\\
c &b & -\sqrt{b d}i\\
\end{matrix}\right)(b,c,d\in\bfk\setminus\{0\}
),\\
N_{5,5}=\left(\begin{matrix}
-\sqrt{ac}i &  0 & a\\
-\frac{\sqrt{c}d}{\sqrt{a}}i & 0 & d\\
c &0 & \sqrt{ac}i\\
\end{matrix}\right),
N_{5,6}=\left(\begin{matrix}
\sqrt{ac}i &  0 & a\\
\frac{\sqrt{c}d}{\sqrt{a}}i & 0 & d\\
c &0 & -\sqrt{ac}i\\
\end{matrix}\right)
(a,c,d\in\bfk\setminus\{0\}),\\
N_{5,7}=\left(
\begin{matrix}
0 & -\frac{a\sqrt{b} i}{\sqrt{d}} & a\\
0 & -\sqrt{bd}i & d\\
0 & b & \sqrt{bd}i\\
\end{matrix}
\right),
N_{5,8}=\left(
\begin{matrix}
0 & \frac{a\sqrt{b} i}{\sqrt{d}} & a\\
0 & \sqrt{bd}i & d\\
0 & b & -\sqrt{bd}i\\
\end{matrix}
\right)(a,b,c,d\in\bfk\setminus\{0\}),\\
N_{5,9}=\left(\begin{matrix}
-\sqrt{ac}i & -\frac{\sqrt{a}bi}{\sqrt{c}} & a\\
0 & 0 & 0\\
c &b & \sqrt{ac}i\\
\end{matrix}\right),
N_{5,10}=\left(\begin{matrix}
\sqrt{ac}i & \frac{\sqrt{a}bi}{\sqrt{c}} & a\\
0 & 0 & 0\\
c &b & -\sqrt{ac}i\\
\end{matrix}\right)(a,b,c\in\bfk\setminus\{0\}),\\
N_{5,11}=\left(\begin{matrix}
0 & 0 & 0\\
0 & -\sqrt{db}i & d\\
0 &b & \sqrt{db}i\\
\end{matrix}\right),
N_{5,12}=\left(\begin{matrix}
0 & 0 & 0\\
0 &\sqrt{d b} i & d\\
0 &b & -\sqrt{db}i\\
\end{matrix}\right)(b\neq0,d\neq0),\\
N_{5,13}=\left(\begin{matrix}
-\frac{e b}{c} & -\frac{e b^2}{c^2} & 0\\
e & \frac{e b}{c} & 0\\
c & b & 0\\
\end{matrix}\right)(b\neq 0,c\neq0),
N_{5,14}=\left(\begin{matrix}
\frac{ae}{d} & -\frac{a^2e}{d^2}  & a\\
e & -\frac{ae}{d} & d\\
0 &0 & 0\\
\end{matrix}\right)(a\neq 0,d\neq0),\\
N_{5,15}=\left(\begin{matrix}
-\sqrt{ac}i & 0 & a\\
0 & 0 & 0\\
c &0 & \sqrt{ac}i\\
\end{matrix}\right),
N_{5,16}=\left(\begin{matrix}
\sqrt{ac}i & 0 & a\\
0 & 0 & 0\\
c&0 & -\sqrt{ac}i\\
\end{matrix}\right)(a\neq 0,c\neq0),\\
N_{5,17}=\left(\begin{matrix}
0 & 0 & 0\\
e &0  & d\\
0 &0 & 0\\
\end{matrix}\right)(d\neq0),
N_{5,18}=\left(\begin{matrix}
0 & 0 & 0\\
e &0  & 0\\
c &0 & 0\\
\end{matrix}\right)(c\neq0),
N_{5,19}=\left(\begin{matrix}
0 & e & 0\\
0 &0  & 0\\
0 &b & 0\\
\end{matrix}\right)(b\neq0),\\
N_{5,20}=\left(\begin{matrix}
0 & e & a\\
0 &0  & 0\\
0 &0 & 0\\
\end{matrix}\right)(a\neq0),
N_{5,21}=\left(\begin{matrix}
-\sqrt{ef}i & e & 0\\
f & \sqrt{ef}i & 0\\
0 &0 &0\\
\end{matrix}\right),
N_{5,22}=\left(\begin{matrix}
\sqrt{ef}i & e & 0\\
f & -\sqrt{ef}i & 0\\
0 &0 & 0\\
\end{matrix}\right)
$\end{minipage}\\
\hline
$NCS(6)$& \begin{minipage}[t]{0.85\textwidth}
\noindent
$N_{6,1}=\left(\begin{matrix}
a & 0 & a\\
0 & 0 & 0\\
-a & 0 & -a\\
\end{matrix}\right)(a\neq 0),
N_{6,2}=
\left(\begin{matrix}
a & 0 & a\\
-c & 0 & c\\
-a & 0 & -a\\
\end{matrix}\right)(a\neq 0),\\
N_{6,3}=\left(\begin{matrix}
\sqrt{ab}i & 0 & a\\
0 & 0 & 0\\
b & 0 & -\sqrt{ab}i\\
\end{matrix}\right)(a\neq 0, b\neq 0),
N_{6,4}=\left(\begin{matrix}
-\sqrt{ab}i & 0 & a\\
0 & 0 & 0\\
b & 0 & \sqrt{ab}i\\
\end{matrix}\right)(a\neq 0, b\neq 0),\\
N_{6,5}=
\left(\begin{matrix}
0 & 0 & a\\
0 & 0 & 0\\
0 & 0 & 0\\
\end{matrix}\right)(a\neq 0),
N_{6,6}=
\left(\begin{matrix}
0 & 0 & 0\\
a & 0 & -a\\
0 & 0 & 0\\
\end{matrix}\right),
N_{6,7}=
\left(\begin{matrix}
0 & 0 & 0\\
0 & 0 & 0\\
a & 0 & 0\\
\end{matrix}\right)(a\neq 0),\\
N_{6,8}=
\left(\begin{matrix}
0 & 0 & 0\\
0 & 0 & 0\\
a & -a & 0\\
\end{matrix}\right)(a\neq 0),
N_{6,9}=
\left(\begin{matrix}
\frac{ab}{a-b} & 0 & \frac{a^2}{b-a}\\
0 & 0 & 0\\
\frac{b^2}{a-b} & b & \frac{ab}{b-a}\\
\end{matrix}\right)(a\neq 0)\\
$
\end{minipage}\\
\hline
\end{longtable}
\end{tiny}

\delete{
\begin{longtable}{|l|l|l|}
\caption[Feasible triples for a highly variable Grid]{Feasible triples for
highly variable Grid, MLMMH.} \label{grid_mlmmh} \\
\hline \multicolumn{1}{|c|}{hj\textbf{Time (s)}} & \multicolumn{1}{c|}{hh\textbf{ Triple chosen}}
& \multicolumn{1}{c|}{hh\textbf{Other feasible triples}} \\
\hline
\endfirsthead

\multicolumn{3}{c}%
{{\bfseries \tablename\ \thetable{} -- continued from previous page}} \\
\hline \multicolumn{1}{|c|}{\textbf{Time (s)}} &
\multicolumn{1}{c|}{\textbf{Triple chosen}}
&\multicolumn{1}{c|}{\textbf{Other feasible triples}} \\
\hline
\endhead

\hline \multicolumn{3}{|r|}{{Continued on next page}} \\ \hline
\endfoot

\hline \hline
\endlastfoot

0 & (1, 11, 13725) & (1, 12, 10980), (1, 13, 8235), (2, 2, 0), (3, 1, 0) \\
2745 & (1, 12, 10980) & (1, 13, 8235), (2, 2, 0), (2, 3, 0), (3, 1, 0) \\
5490 & (1, 12, 13725) & (2, 2, 2745), (2, 3, 0), (3, 1, 0) \\
8235 & (1, 12, 16470) & (1, 13, 13725), (2, 2, 2745), (2, 3, 0), (3, 1, 0) \\
10980 & (1, 12, 16470) & (1, 13, 13725), (2, 2, 2745), (2, 3, 0), (3, 1, 0) \\
13725 & (1, 12, 16470) & (1, 13, 13725), (2, 2, 2745), (2, 3, 0), (3, 1, 0) \\
16470 & (1, 13, 16470) & (2, 2, 2745), (2, 3, 0), (3, 1, 0) \\
19215 & (1, 12, 16470) & (1, 13, 13725), (2, 2, 2745), (2, 3, 0), (3, 1, 0) \\
21960 & (1, 12, 16470) & (1, 13, 13725), (2, 2, 2745), (2, 3, 0), (3, 1, 0) \\
24705 & (1, 12, 16470) & (1, 13, 13725), (2, 2, 2745), (2, 3, 0), (3, 1, 0) \\
27450 & (1, 12, 16470) & (1, 13, 13725), (2, 2, 2745), (2, 3, 0), (3, 1, 0) \\
30195 & (2, 2, 2745) & (2, 3, 0), (3, 1, 0) \\
32940 & (1, 13, 16470) & (2, 2, 2745), (2, 3, 0), (3, 1, 0) \\
35685 & (1, 13, 13725) & (2, 2, 2745), (2, 3, 0), (3, 1, 0) \\
38430 & (1, 13, 10980) & (2, 2, 2745), (2, 3, 0), (3, 1, 0) \\
41175 & (1, 12, 13725) & (1, 13, 10980), (2, 2, 2745), (2, 3, 0), (3, 1, 0) \\
43920 & (1, 13, 10980) & (2, 2, 2745), (2, 3, 0), (3, 1, 0) \\
46665 & (2, 2, 2745) & (2, 3, 0), (3, 1, 0) \\
49410 & (2, 2, 2745) & (2, 3, 0), (3, 1, 0) \\
52155 & (1, 12, 16470) & (1, 13, 13725), (2, 2, 2745), (2, 3, 0), (3, 1, 0) \\
54900 & (1, 13, 13725) & (2, 2, 2745), (2, 3, 0), (3, 1, 0) \\
57645 & (1, 13, 13725) & (2, 2, 2745), (2, 3, 0), (3, 1, 0) \\
60390 & (1, 12, 13725) & (2, 2, 2745), (2, 3, 0), (3, 1, 0) \\
63135 & (1, 13, 16470) & (2, 2, 2745), (2, 3, 0), (3, 1, 0) \\
65880 & (1, 13, 16470) & (2, 2, 2745), (2, 3, 0), (3, 1, 0) \\
68625 & (2, 2, 2745) & (2, 3, 0), (3, 1, 0) \\
71370 & (1, 13, 13725) & (2, 2, 2745), (2, 3, 0), (3, 1, 0) \\
74115 & (1, 12, 13725) & (2, 2, 2745), (2, 3, 0), (3, 1, 0) \\
76860 & (1, 13, 13725) & (2, 2, 2745), (2, 3, 0), (3, 1, 0) \\
79605 & (1, 13, 13725) & (2, 2, 2745), (2, 3, 0), (3, 1, 0) \\
82350 & (1, 12, 13725) & (2, 2, 2745), (2, 3, 0), (3, 1, 0) \\
85095 & (1, 12, 13725) & (1, 13, 10980), (2, 2, 2745), (2, 3, 0), (3, 1, 0) \\
87840 & (1, 13, 16470) & (2, 2, 2745), (2, 3, 0), (3, 1, 0) \\
90585 & (1, 13, 16470) & (2, 2, 2745), (2, 3, 0), (3, 1, 0) \\
93330 & (1, 13, 13725) & (2, 2, 2745), (2, 3, 0), (3, 1, 0) \\
96075 & (1, 13, 16470) & (2, 2, 2745), (2, 3, 0), (3, 1, 0) \\
98820 & (1, 13, 16470) & (2, 2, 2745), (2, 3, 0), (3, 1, 0) \\
101565 & (1, 13, 13725) & (2, 2, 2745), (2, 3, 0), (3, 1, 0) \\
104310 & (1, 13, 16470) & (2, 2, 2745), (2, 3, 0), (3, 1, 0) \\
107055 & (1, 13, 13725) & (2, 2, 2745), (2, 3, 0), (3, 1, 0) \\
109800 & (1, 13, 13725) & (2, 2, 2745), (2, 3, 0), (3, 1, 0) \\
112545 & (1, 12, 16470) & (1, 13, 13725), (2, 2, 2745), (2, 3, 0), (3, 1, 0) \\
115290 & (1, 13, 16470) & (2, 2, 2745), (2, 3, 0), (3, 1, 0) \\
118035 & (1, 13, 13725) & (2, 2, 2745), (2, 3, 0), (3, 1, 0) \\
120780 & (1, 13, 16470) & (2, 2, 2745), (2, 3, 0), (3, 1, 0) \\
123525 & (1, 13, 13725) & (2, 2, 2745), (2, 3, 0), (3, 1, 0) \\
126270 & (1, 12, 16470) & (1, 13, 13725), (2, 2, 2745), (2, 3, 0), (3, 1, 0) \\
129015 & (2, 2, 2745) & (2, 3, 0), (3, 1, 0) \\
131760 & (2, 2, 2745) & (2, 3, 0), (3, 1, 0) \\
134505 & (1, 13, 16470) & (2, 2, 2745), (2, 3, 0), (3, 1, 0) \\
137250 & (1, 13, 13725) & (2, 2, 2745), (2, 3, 0), (3, 1, 0) \\
139995 & (2, 2, 2745) & (2, 3, 0), (3, 1, 0) \\
142740 & (2, 2, 2745) & (2, 3, 0), (3, 1, 0) \\
145485 & (1, 12, 16470) & (1, 13, 13725), (2, 2, 2745), (2, 3, 0), (3, 1, 0) \\
148230 & (2, 2, 2745) & (2, 3, 0), (3, 1, 0) \\
150975 & (1, 13, 16470) & (2, 2, 2745), (2, 3, 0), (3, 1, 0) \\
153720 & (1, 12, 13725) & (2, 2, 2745), (2, 3, 0), (3, 1, 0) \\
156465 & (1, 13, 13725) & (2, 2, 2745), (2, 3, 0), (3, 1, 0) \\
159210 & (1, 13, 13725) & (2, 2, 2745), (2, 3, 0), (3, 1, 0) \\
161955 & (1, 13, 16470) & (2, 2, 2745), (2, 3, 0), (3, 1, 0) \\
164700 & (1, 13, 13725) & (2, 2, 2745), (2, 3, 0), (3, 1, 0) \\
115290 & (1, 13, 16470) & (2, 2, 2745), (2, 3, 0), (3, 1, 0) \\
118035 & (1, 13, 13725) & (2, 2, 2745), (2, 3, 0), (3, 1, 0) \\
120780 & (1, 13, 16470) & (2, 2, 2745), (2, 3, 0), (3, 1, 0) \\
123525 & (1, 13, 13725) & (2, 2, 2745), (2, 3, 0), (3, 1, 0) \\
126270 & (1, 12, 16470) & (1, 13, 13725), (2, 2, 2745), (2, 3, 0), (3, 1, 0) \\
129015 & (2, 2, 2745) & (2, 3, 0), (3, 1, 0) \\
131760 & (2, 2, 2745) & (2, 3, 0), (3, 1, 0) \\
\end{longtable}

\begin{center}
\begin{longtable}{cccccc}
 \caption{The longtable package defines an environment that has most of the features of the tabular environment.}\\
 \hline
 \textbf{First entry} & \textbf{Second entry} & \textbf{Third entry} & \textbf{Fourth entry} & \textbf{Fifth entry} & \textbf{Sixth entry}\\
 \hline
 \endfirsthead
 \multicolumn{6}{c}%
 {\tablename\ \thetable\ -- \textit{Continued from previous page}} \\
 \hline
 \textbf{First entry} & \textbf{Second entry} & \textbf{Third entry} & \textbf{Fourth entry} & \textbf{Fifth entry} & \textbf{Sixth entry}\\
 \hline
 \endhead
 \hline \multicolumn{6}{r}{\textit{Continued on next page}} \\
 \endfoot
 \hline
 \endlastfoot
 1 & 2 & 3 & 4 & 5 & 6
 \end{longtable}
 \end{center}
}

\delete{
{\allowdisplaybreaks
\begin{table}[!thbp]
\caption{RBOs on 3-dimensional noncommutative semigroup algebras \mlabel{ta:rbncs2}}
\centering
\tiny
\begin{tabular}{|c|p{11.5cm}|c|}
\hline
\text{noncommutative semigroup of order 3} &\text{\quad\quad
Rota-Baxter operators on 3-dimensional noncommutative semigroup algebras}  \\
\hline
$NCS(3,4):=$
\begin{tabular}{c|ccccc}
$\cdot$ &$e_1$ & $e_2$& $e_3$\\
\hline
$e_1$ &$e_1$ & $e_1$ &$e_1$\\
$e_2$ &$e_2$ & $e_2$ &$e_2$\\
$e_3$&$e_1$ & $e_1$ &$e_1$\\
\end{tabular}
& \begin{minipage}[t]{0.85\textwidth}
\noindent
$N_{41}=\left(\begin{matrix}
-a & 0 & a\\
b & 0 & -b\\
c  & 0 & -c\\
\end{matrix}\right)(a\neq 0, a+b\neq 0)
\bigcup N_{42}=\left(\begin{matrix}
-a & 0 & a\\
-a & 0 & a\\
c  & d & -c-d\\
\end{matrix}\right)(a\neq 0)\\
\bigcup
N_{43}=\left(\begin{matrix}
0 & 0 & 0\\
-a & 0 & a\\
b  & 0 & -b\\
\end{matrix}\right)(a\neq 0)
\bigcup
N_{44}=\left(\begin{matrix}
0 & 0 & 0\\
0 & 0 & 0\\
b  & c & -(b+c)\\
\end{matrix}\right)(c\neq 0)
\\
\bigcup
N_{45}=\left(\begin{matrix}
0&0&0\\
a&0&0\\
b&0&-b\\
\end{matrix}\right)
\bigcup
N_{46}=\left(\begin{matrix}
 a&-(a+b)& b\\
 a&-(a+b)&b\\
 c&d&-(c+d)\\
 \end{matrix}\right)(b\neq 0, a+b\neq 0)\\
\bigcup N_{47}=\left(
 \begin{matrix}
 a&b&0\\
 -\frac{a^2}{b}&-a&0\\
 c&b&a-c\\
 \end{matrix}\right)(a+b\neq 0, b\neq 0)
\bigcup
N_{48}=\left(\begin{matrix}
 -b&b&0\\
 -b&b&0\\
 c&d&-c-d\\
 \end{matrix}\right)(b\neq 0,c,d\in\bfk)\\
 $
\end{minipage}\\
\hline
NCS(3,5):=
\begin{tabular}{c|ccccc}
$\cdot$ &$e_1$ & $e_2$& $e_3$\\
\hline
$e_1$ &$e_1$ & $e_1$ &$e_1$\\
$e_2$ &$e_2$ & $e_2$ &$e_2$\\
$e_3$ &$e_3$ & $e_3$ &$e_3$\\
\end{tabular}
& \begin{minipage}[t]{0.85\textwidth}
\noindent
$\Bigg\{N_{51}=\left(
\begin{matrix}
\frac{ac}{F} & \frac{ab}{F} & a\\
\frac{cd}{F} & \frac{bd}{F} & d\\
c & b & F\\
\end{matrix}
\right)
\cup N_{52}=\left(
\begin{matrix}
-\frac{ac}{F} & -\frac{ab}{F} & a\\
-\frac{cd}{F} & -\frac{bd}{F} & d\\
c & b & -F\\
\end{matrix}
\right)(a\neq0,b\neq0,c\neq0,d\neq0, F:=\sqrt{-(ac+bd)}\Bigg\}\\
\,
\bigcup \Bigg\{N_{53}=\left(\begin{matrix}
0 &  0 & 0\\
-\frac{c\sqrt{d}}{\sqrt{b}}i & -\sqrt{b d}i & d\\
c &b & \sqrt{bd}i\\
\end{matrix}\right)
\,\cup\,
N_{54}=\left(\begin{matrix}
0 & 0 & 0\\
\frac{c\sqrt{d}}{\sqrt{b}}i & \sqrt{b d}i & d\\
c &b & -\sqrt{b d}i\\
\end{matrix}\right)(b\neq0,c\neq 0,d\neq0)\Bigg\}\\
\,\bigcup\,
\Bigg\{N_{55}=\left(\begin{matrix}
-\sqrt{ac}i &  0 & a\\
-\frac{\sqrt{c}d}{\sqrt{a}}i & 0 & d\\
c &0 & \sqrt{ac}i\\
\end{matrix}\right)
\,\cup\,
N_{56}=\left(\begin{matrix}
\sqrt{ac}i &  0 & a\\
\frac{\sqrt{c}d}{\sqrt{a}}i & 0 & d\\
c &0 & -\sqrt{ac}i\\
\end{matrix}\right)
(a\neq 0,c\neq 0,d\neq0)\Bigg\}\\
\,\bigcup\,
\Bigg\{N_{57}=\left(
\begin{matrix}
0 & -\frac{a\sqrt{b} i}{\sqrt{d}} & a\\
0 & -\sqrt{bd}i & d\\
0 & b & \sqrt{bd}i\\
\end{matrix}
\right)
\,\cup\,
N_{58}=\left(
\begin{matrix}
0 & \frac{a\sqrt{b} i}{\sqrt{d}} & a\\
0 & \sqrt{bd}i & d\\
0 & b & -\sqrt{bd}i\\
\end{matrix}
\right)(a\neq0,b\neq0,d\neq0)\Bigg\}\\
\bigcup \Bigg\{N_{59}=\left(\begin{matrix}
-\sqrt{ac}i & -\frac{\sqrt{a}bi}{\sqrt{c}} & a\\
0 & 0 & 0\\
c &b & \sqrt{ac}i\\
\end{matrix}\right)
\,\cup\,
N_{510}=\left(\begin{matrix}
\sqrt{ac}i & \frac{\sqrt{a}bi}{\sqrt{c}} & a\\
0 & 0 & 0\\
c &b & -\sqrt{ac}i\\
\end{matrix}\right)(a\neq0,b\neq0,c\neq0)
\Bigg\}\\
\,\bigcup \Bigg\{
N_{511}=\left(\begin{matrix}
0 & 0 & 0\\
0 & -\sqrt{db}i & d\\
0 &b & \sqrt{db}i\\
\end{matrix}\right)\,\cup\,
N_{512}=\left(\begin{matrix}
0 & 0 & 0\\
0 &\sqrt{d b} i & d\\
0 &b & -\sqrt{db}i\\
\end{matrix}\right)(b\neq0,d\neq0)\Bigg\}\\
\,\bigcup \Bigg\{
N_{513}=\left(\begin{matrix}
-\frac{e b}{c} & -\frac{e b^2}{c^2} & 0\\
e & \frac{e b}{c} & 0\\
c & b & 0\\
\end{matrix}\right)(b\neq 0,c\neq0)
\,\cup\,
N_{514}=\left(\begin{matrix}
\frac{ae}{d} & -\frac{a^2e}{d^2}  & a\\
e & -\frac{ae}{d} & d\\
0 &0 & 0\\
\end{matrix}\right)(a\neq 0,d\neq0)\Bigg\}\\
\,\bigcup \Bigg\{
N_{515}=\left(\begin{matrix}
-\sqrt{ac}i & 0 & a\\
0 & 0 & 0\\
c &0 & \sqrt{ac}i\\
\end{matrix}\right)
\,\cup\,
N_{516}=\left(\begin{matrix}
\sqrt{ac}i & 0 & a\\
0 & 0 & 0\\
c&0 & -\sqrt{ac}i\\
\end{matrix}\right)(a\neq 0,c\neq0)\Bigg\}\\
\,\bigcup \Bigg\{
N_{517}=\left(\begin{matrix}
0 & 0 & 0\\
e &0  & d\\
0 &0 & 0\\
\end{matrix}\right)(d\neq0)
\,\cup\,
N_{518}=\left(\begin{matrix}
0 & 0 & 0\\
e &0  & 0\\
c &0 & 0\\
\end{matrix}\right)(c\neq0)\Bigg\}\\
\,\bigcup \Bigg\{
N_{519}=\left(\begin{matrix}
0 & e & 0\\
0 &0  & 0\\
0 &b & 0\\
\end{matrix}\right)(b\neq0)
\,\cup\,
N_{520}=\left(\begin{matrix}
0 & e & a\\
0 &0  & 0\\
0 &0 & 0\\
\end{matrix}\right)(a\neq0)\Bigg\}\\
\,\bigcup \Bigg\{
N_{521}=\left(\begin{matrix}
-\sqrt{ef}i & e & 0\\
f & \sqrt{ef}i & 0\\
0 &0 &0\\
\end{matrix}\right)
\,\cup\,
N_{522}=\left(\begin{matrix}
\sqrt{ef}i & e & 0\\
f & -\sqrt{ef}i & 0\\
0 &0 & 0\\
\end{matrix}\right)\Bigg\}
$\end{minipage}\\
\hline
$NCS(3,6)$:=\begin{tabular}{c|cccc}
$\cdot$ &$e_1$ & $e_2$& $e_3$\\
\hline
$e_1$ &$e_1$ & $e_1$ &$e_1$\\
$e_2$ &$e_1$ & $e_2$ &$e_3$\\
$e_3$ &$e_3$ & $e_3$ &$e_3$\\
\end{tabular}& \begin{minipage}[t]{0.85\textwidth}
\noindent
$N_{61}=\left(\begin{matrix}
a & 0 & a\\
0 & 0 & 0\\
-a & 0 & -a\\
\end{matrix}\right)(a\neq 0)
\bigcup N_{62}=
\left(\begin{matrix}
a & 0 & a\\
-c & 0 & c\\
-a & 0 & -a\\
\end{matrix}\right)(a\neq 0)\\
N_{63}=\left(\begin{matrix}
\sqrt{ab}i & 0 & a\\
0 & 0 & 0\\
b & 0 & -\sqrt{ab}i\\
\end{matrix}\right)(a\neq 0, b\neq 0)
\bigcup
N_{64}=\left(\begin{matrix}
-\sqrt{ab}i & 0 & a\\
0 & 0 & 0\\
b & 0 & \sqrt{ab}i\\
\end{matrix}\right)(a\neq 0, b\neq 0)\\
\bigcup N_{65}=
\left(\begin{matrix}
0 & 0 & a\\
0 & 0 & 0\\
0 & 0 & 0\\
\end{matrix}\right)(a\neq 0)
\bigcup N_{66}=
\left(\begin{matrix}
0 & 0 & 0\\
a & 0 & -a\\
0 & 0 & 0\\
\end{matrix}\right)
\bigcup N_{67}=
\left(\begin{matrix}
0 & 0 & 0\\
0 & 0 & 0\\
a & 0 & 0\\
\end{matrix}\right)(a\neq 0)\\
\bigcup N_{68}=
\left(\begin{matrix}
0 & 0 & 0\\
0 & 0 & 0\\
a & -a & 0\\
\end{matrix}\right)(a\neq 0)
\bigcup N_{69}=
\left(\begin{matrix}
\frac{ab}{a-b} & 0 & \frac{a^2}{b-a}\\
0 & 0 & 0\\
\frac{b^2}{a-b} & b & \frac{ab}{b-a}\\
\end{matrix}\right)(a\neq 0)\\
$
\end{minipage}\\
\hline
\end{tabular}
\end{table}
}
}

\subsection{Proof of Theorem~\mref{thm:rbonc3}}
\mlabel{ss:proofnc}

We will prove Theorem~\mref{thm:rbonc3} in Section~\mref{ss:ncs1}-\mref{ss:ncs6}, one for each of the six semigroups in Table~\mref{ta:ncs3}.

\subsubsection{The proof for $\bfk[NCS(1)]$}
\mlabel{ss:ncs1}
We first prove that the matrices $(c_{ij})_{1\leq i, j\leq 3}$ of the Rota-Baxter operators $P$ on the semigroup algebra $\bfk[NCS(1)]$ are given by $N_{1,i},1\leq i\leq 5$, in Table~\mref{ta:rbncs1}.

Applying the Cayley table of $NCS(1)$
\delete{
\begin{equation*}\tiny
NCS(1):=
\begin{array}{c|ccc}
\cdot &e_1 & e_2& e_3\\
\hline
e_1 &e_1 & e_1 &e_1\\
e_2 &e_1 & e_2 &e_1\\
e_3 &e_1 & e_3 &e_1\\
\end{array},
\end{equation*}}in Eqs.~(\mref{eq:3eijm}) and then taking $i=1,j=1,2$ with $1\leq m\leq 3$; $i=j=2$ with $1\leq m\leq 3$; $i=2,j=3$ with $m=1$; $i=3,j=2$ with $1\leq m\leq 3$ and $i=j=3$ with $m=1,3$, we obtain
\begin{eqnarray}
&c_{11}^2=c_{13}^2+c_{12}c_{13},
\mlabel{eq:ncs311}\\
&c_{12}^2=2c_{12}(c_{11}+c_{12}+c_{13})
,\mlabel{eq:ncs312}\\
&2c_{13}^2+c_{13}(2c_{11}+c_{12})=0,
\mlabel{eq:ncs313}\\
&c_{12}c_{23}+c_{13}(c_{21}+c_{23})
=c_{11}^2+c_{13}c_{31},\mlabel{eq:ncs314}\\
&c_{12}(c_{11}+c_{21}+c_{22}+c_{23})
+c_{13}c_{32}=0,\mlabel{eq:ncs315}\\
&c_{13}(c_{11}+c_{21}+c_{23})+c_{12}c_{23}
+c_{13}c_{33}=0,\mlabel{eq:ncs316}\\
&c_{21}(c_{21}+c_{23})+c_{23}(c_{21}
+c_{22}+c_{23})=c_{11}(2c_{21}+c_{23})
+c_{23}c_{31},\mlabel{eq:ncs3113}\\
&c_{22}^2+c_{23}c_{32}+c_{12}(2c_{21}
+c_{23})=0,\mlabel{eq:ncs3114}\\
&c_{23}(c_{22}+c_{33})+c_{13}(2c_{21}
+c_{23})=0,\mlabel{eq:ncs3115}\\
&(c_{21}+c_{22}+c_{23})
(c_{31}+c_{33})=c_{11}(c_{21}+c_{22}+
c_{23}+c_{31}+c_{33}),\mlabel{eq:ncs3116}\\
&(c_{31}+c_{32}+c_{33})
(c_{21}+c_{23})=c_{11}(c_{31}+c_{21}+
c_{23})+c_{32}c_{21}+c_{31}c_{33},
\mlabel{eq:ncs3122}\\
&c_{12}(c_{31}+c_{21}+c_{23})
+c_{32}(c_{33}+c_{22})=0,\mlabel{eq:ncs3123}\\
&c_{33}^2+c_{23}c_{32}+c_{13}(c_{31}
+c_{21}+c_{23})=0,\mlabel{eq:ncs3124}\\
&(c_{31}+c_{33})(c_{31}+c_{32}+
c_{33})=c_{11}(2c_{31}+c_{32}+2c_{33})
,\mlabel{eq:ncs3125}\\
&c_{13}(2c_{31}+c_{32}+
2c_{33})=0
.\mlabel{eq:ncs3127}
\end{eqnarray}
By Eq.~(\mref{eq:ncs313}) we have
$c_{12}c_{13}=-2c_{13}^2-2c_{11}c_{13}.$
Then by Eq.~(\mref{eq:ncs311}) we get
$(c_{11}+c_{13})^2=0.$ So we get $c_{11}+c_{13}=0$. Then
Eq.~(\mref{eq:ncs312}) gives $c_{12}=0$. We divide the rest of the proof into two cases depending on whether or not $c_{11}=0$.

\smallskip
\noindent
{\bf Case 1. $c_{11}=0$:} Then $c_{13}=0$.
\delete{Thus Eqs.~(\mref{eq:ncs311})-(\mref{eq:ncs3125})
are equivalent to the following equations.
\begin{eqnarray}
&c_{11}=c_{12}=c_{13}=0,
\mlabel{eq:ncs31e1}\\
&(c_{21}+c_{23})^2=c_{23}(c_{31}-c_{22}),
\mlabel{eq:ncs31e2}\\
&c_{22}^2+c_{23}c_{32}=0,
\mlabel{eq:ncs31e3}\\
&c_{23}(c_{22}+c_{33})=0,
\mlabel{eq:ncs31e4}\\
&(c_{21}+c_{22}+c_{23})(c_{31}+c_{33})=0,
\mlabel{eq:ncs31e}\\
&(c_{21}+c_{23})(c_{31}+c_{32}+c_{33})=
c_{32}c_{21}+c_{31}c_{33}
,\mlabel{eq:ncs31e5}\\
&c_{32}(c_{22}+c_{33})=0,
\mlabel{eq:ncs31e6}\\
&c_{32}c_{23}+c_{33}^2=0,
\mlabel{eq:ncs31e7}\\
&(c_{31}+c_{33})(c_{31}+c_{32}+c_{33})=0.
\mlabel{eq:ncs31e8}
\end{eqnarray}
}
There are two subcases to consider.

\smallskip
\noindent
{\bf Subcases 1. $c_{23}=0$:} Then by Eqs.~(\mref{eq:ncs3114}) and~(\mref{eq:ncs3124}), we have $c_{22}=c_{33}=0$. So by Eq.~(\mref{eq:ncs3113}), we get $c_{21}=0$.
\delete{
Then Eqs.~(\mref{eq:ncs31e1})-(\mref{eq:ncs31e8})
are equivalent to
$$\left\{\begin{array}{cc}
c_{11}=c_{12}=c_{13}=0,\\
c_{21}=c_{22}=c_{23}=0,\\
c_{33}=0,\\
c_{31}(c_{31}+c_{32})=0.\\
\end{array}\right.$$
}
Denote $a=c_{32}$, where $a\in\bfk$.
If $c_{31}=0$, then we get the solution
$$N_{1,1}=\left(
\begin{matrix}
0&0&0\\
0&0&0\\
0&a&0\\
\end{matrix}\right)\quad (a\in\bfk).$$
If $c_{31}\neq 0$, then by Eq.~(\mref{eq:ncs3125}), we have $c_{31}+c_{32}=0$. Thus we get $c_{31}=-c_{32}=-a$. So $a\neq0$. Thus we obtain the solution
$$N_{1,2}=\left(\begin{matrix}
0&0&0\\
0&0&0\\
-a&a&0\\
\end{matrix}\right)\quad (a\in\bfk, a\neq 0).$$

\smallskip
\noindent
{\bf Subcase 2. $c_{23}\neq 0$:} Denote $a=c_{23}$. Then $a\neq 0$. By Eq.~(\mref{eq:ncs3115}), we have $c_{22}+c_{33}=0$.
Denote $b=c_{33}$.  Then $c_{22}=-b$. Thus Eq.~(\mref{eq:ncs3114}) gives $c_{32}=-\frac{b^2}{a}$.
We subdivide further into two cases.
\begin{enumerate}
\item
If $c_{31}+c_{33}\neq 0$, then by Eq.~(\mref{eq:ncs3125}), we have $c_{31}+c_{32}+c_{33}=0$. So we get
$$c_{31}=-c_{32}-c_{33}=\frac{b^2-ba}
{a}\,.$$
Note that if $c_{33}=b=0$, then $c_{31}=0$, a contradiction. Thus $b\neq 0$, and so $c_{32}\neq 0$.
By Eq.~(\mref{eq:ncs3116}), we have $c_{21}+c_{22}+c_{23}=0$. Thus  $c_{21}=b-a$. Thus we get the solutions
$$N_{1,3}=\left(
\begin{matrix}
0&0&0\\
b-a&-b&a\\
\frac{b(b-a)}{a}&-\frac{b^2}{a}&b\\
\end{matrix}
\right)\quad (a, b\in\bfk\setminus \{0\}).$$
\item
If $c_{31}+c_{33}=0,$ then $c_{31}=-b$. Since $c_{22}=-b$, $c_{31}-c_{22}=0$.
Then by Eq.~(\mref{eq:ncs3113}), we have $(c_{21}+c_{23})^2=c_{23}(c_{31}-c_{22})=0$.
Thus $c_{21}=-a$.
 Thus we obtain the solutions
$$N_{1,4}=\left(
\begin{matrix}
0&0&0\\
-a&-b&a\\
-b&-\frac{b^2}{a}&b\\
\end{matrix}\right)\quad (a, b\in\bfk, a\neq 0).$$
\end{enumerate}

\smallskip
\noindent
{\bf Case 2. $c_{11}\neq 0$:} Denote $a=c_{11}$. Then $c_{13}=-a\neq 0$. Since $c_{12}=0$,  Eq.~(\mref{eq:ncs315}) becomes $c_{13}c_{32}=0$. Then $c_{32}=0$. So by Eq.~(\mref{eq:ncs3114}), we obtain $c_{22}=0$. Further, by Eq.~(\mref{eq:ncs3127}), we have $c_{31}+c_{33}=0$. Thus by Eq.~(\mref{eq:ncs3116}), $c_{21}+c_{23}=0$.

We divide into two subcases to consider.

\smallskip
\noindent
{\bf Subcase 1. $c_{23}=0$:}   Then   $c_{21}=0$.
By Eq.~(\mref{eq:ncs314}), $c_{11}^2+c_{13}c_{31}=0$. Then we obtain $c_{31}=a$. So we have $c_{33}=-a$.
Thus we obtain the solutions
$$N_{1,5_1}:=\left(
\begin{matrix}
a&0&-a\\
0&0&0\\
a&0&-a\\
\end{matrix}\right)\quad (a\in \bfk, a\neq 0).$$

\smallskip
\noindent
{\bf Subcase 2. $c_{23}\neq 0$:} Denote $b=c_{23}$. Then $b\neq 0$.  Thus by $c_{21}+c_{23}=0$, we have $c_{21}=-b$ and  by Eq.~(\mref{eq:ncs3115}), $c_{13}c_{21}+c_{23}c_{33}=0$, we get
$c_{33}=-a$. Then we solutions
$$N_{1,5_2}:=
\left(\begin{matrix}
a&0&-a\\
-b&0&b\\
a&0&-a\\
\end{matrix}\right)\quad (a, b\in\bfk\setminus \{0\}).$$
In summary, we get
get solutions of Eqs.~(\mref{eq:ncs311})-(\mref{eq:ncs3127})
$$N_{1,5}=
\left(\begin{matrix}
a&0&-a\\
-b&0&b\\
a&0&-a\\
\end{matrix}\right)\quad (a, b\in \bfk, a\neq 0).$$
It can be checked that they also satisfy the other equations in Eq.~(\mref{eq:3eijm}) and hence give matrices of Rota-Baxter operators on $\bfk[NCS(1)]$.

\subsubsection{The proof for $\bfk[NCS(2)]$}
\mlabel{ss:ncs2}
We next prove that the matrices of Rota-Baxter operators $P$ on the semigroup algebra $\bfk[NCS(2)]$ are given by $N_{2,i}, 1\leq i\leq 9$, in Table~\mref{ta:rbncs1}.

Applying the Cayley table of $NCS(2)$
\delete{
\begin{equation*}
\begin{array}{c|ccc}
\cdot &e_1 & e_2& e_3\\
\hline
e_1 &e_1 & e_1 &e_1\\
e_2 &e_1 & e_2 &e_1\\
e_3 &e_3 & e_3 &e_3\\
\end{array},
\end{equation*}}in Eqs.~(\mref{eq:3eijm}) and then
taking $i=j=1$ with $1\leq m\leq 3$; $i=1,j=2$ with $m=3$; $i=2,1\leq j\leq 3$ with $1\leq m\leq 3$;  $i=3,j=2$ with $m=1$ and $i=j=3$ with $1\le m\leq 3$, we obtain
\begin{eqnarray}
&c_{12}c_{13}=c_{11}^2+c_{13}c_{31},
\mlabel{eq:ncs321}\\
&c_{12}^2+2c_{11}c_{12}+c_{13}(c_{12}+
c_{32})=0,\mlabel{eq:ncs322}\\
&c_{13}(c_{11}+c_{12}+c_{33})=0,
\mlabel{eq:ncs323}\\
&c_{13}(c_{11}+c_{33})+c_{12}c_{23}=0,
\mlabel{eq:ncs326}\\
&c_{13}(c_{21}+c_{22})=c_{11}(
c_{11}+c_{13})+c_{23}c_{31},
\mlabel{eq:ncs3210}\\
&c_{12}(c_{21}+c_{22}+c_{11}+c_{13}
)+c_{23}c_{32}=0,\mlabel{eq:ncs3211}\\
&c_{23}(c_{11}+c_{13})=
c_{13}(c_{21}+c_{22}+c_{11}+c_{13})+c_{23}
c_{33},\mlabel{eq:ncs3212}\\
&c_{21}(c_{21}+c_{23})+c_{22}c_{23}=
c_{11}(2c_{21}+c_{23})+c_{23}c_{31},
\mlabel{eq:ncs3213}\\
&c_{22}^2+c_{12}(2c_{21}+c_{23})+c_{23}
c_{32}=0,\mlabel{eq:ncs3214}\\
&c_{23}(c_{21}+c_{23})=
c_{13}(2c_{21}+c_{23})+c_{23}(c_{22}+
c_{33}),\mlabel{eq:ncs3215}\\
&(c_{21}+c_{22})(c_{31}+c_{33})
=c_{11}(c_{21}+c_{22}+c_{31}+c_{33})+
c_{23}c_{31},\mlabel{eq:ncs3216}\\
&c_{12}(c_{21}+c_{22}+c_{31}+c_{33})+
c_{23}c_{32}=0,\mlabel{eq:ncs3217}\\
&c_{23}c_{31}=c_{13}(c_{21}+c_{22}+c_{31}
+c_{33}),\mlabel{eq:ncs3218}\\
&c_{32}c_{23}=c_{31}(c_{11}+c_{33}),
\mlabel{eq:ncs3222}\\
&c_{32}(c_{31}+c_{33})=c_{11}(
c_{31}+c_{32})+c_{31}c_{33},
\mlabel{eq:ncs3225}\\
&c_{12}(c_{31}+c_{32})+c_{32}(
c_{31}+2c_{33})=0,\mlabel{eq:ncs3226}\\
&c_{33}^2+c_{13}(c_{31}+c_{32})=0.\mlabel{eq:ncs3227}
\end{eqnarray}
We divide the proof into two cases.

\smallskip
\noindent
{\bf Case 1. $c_{12}=0$:} Then Eqs.~(\mref{eq:ncs3211}) and ~(\mref{eq:ncs3214}) give $c_{22}=0$.
There are two subcases to consider.

\smallskip
\noindent
{\bf Subcase 1. $c_{13}=0$:} Then
by Eqs.~(\mref{eq:ncs321}) and~(\mref{eq:ncs3227}), we have $c_{11}=0$ and $c_{33}=0$. Thus Eqs.~(\mref{eq:ncs3210}) and ~(\mref{eq:ncs3212}) give $c_{23}c_{31}=0$ and $c_{23}c_{33}=0$.
Thus by Eqs.~(\mref{eq:ncs3213})and ~(\mref{eq:ncs3215}), we have $c_{21}(c_{21}+c_{23})=0$ and $c_{23}(c_{21}+c_{23})=0.$ So $c_{21}+c_{23}=0$.
\begin{enumerate}
\item
Let $c_{21}\neq0$. Denote $c=c_{21}$. Then by $c_{21}+c_{23}=0$, we have $c_{23}=-c$. So $c_{23}\neq 0$. Thus $c_{31}=0$. By Eq.~(\mref{eq:ncs3211}), we have $c_{32}=0$.
Thus we get the solutions
$$N_{2,1}=\left(
\begin{matrix}
0&0&0\\
c&0&-c\\
0&0&0\\
\end{matrix}\right)\quad(c\in \bfk, c\neq 0).$$
\item
Let $c_{21}=0$. Then $c_{23}=0$. By Eq.~(\mref{eq:ncs3225}), we have $c_{32}c_{31}=0$. Thus $c_{32}=0$ or $c_{31}=0$. Let $c_{31}=a$, and let $c_{32}=b$, where $a,b\in\bfk$. Thus $ab=0$.  Then we get the solutions
$$
N_{2,2}=
\left(\begin{matrix}
0&0&0\\
0&0&0\\
a&0&0\\
\end{matrix}\right)\quad (a\in\bfk)\quad\text{and}
\quad
N_{2,3}=
\left(\begin{matrix}
0&0&0\\
0&0&0\\
0&b&0\\
\end{matrix}\right)\quad (b\in\bfk).$$
\end{enumerate}

\smallskip
\noindent
{\bf Subcase 2. $c_{13}\neq 0$:} Denote $a=c_{13}$.
By Eq.~(\mref{eq:ncs322}), we have $c_{13}c_{32}=0$ and hence $c_{32}=0$.
Eq.~(\mref{eq:ncs323}) gives $c_{13}(c_{11}+c_{33})=0$ and hence $c_{11}+
c_{33}=0$. Denote $b=c_{11}$. Then $c_{33}=-b$. So Eq.~(\mref{eq:ncs321}) gives $c_{31}=-\frac{b^2}{a}$.
By Eq.~(\mref{eq:ncs3218}), we have
$c_{21}=-\frac{b^2}{a^2}c_{23}
+\frac{b^2}{a}+b$. Using Eq.~(\mref{eq:ncs3216}), we obtain
\begin{equation}
b^2(a+b)^2c_{23}=ab^2(a+b)^2.
\mlabel{eq:ncs32r1}
\end{equation}
We continue by subdividing further.
\begin{enumerate}
\item
Let $c_{11}+c_{13}\neq 0$. Then $a+b\neq 0$.
\begin{enumerate}
\item
If $b\neq 0$, then by Eq.~(\mref{eq:ncs32r1}), we have $c_{23}=a$.  Thus $c_{21}=b$. So we obtain the solutions
$$N_{2,4}=\left(
\begin{matrix}
b&0&a\\
b&0&a\\
-\frac{b^2}{a}& 0&-b\\
\end{matrix}\right)\quad (a\neq 0, b\neq 0, a+b\neq 0).$$
\item
If $b=0$, then $c_{11}=c_{33}=0$, $c_{21}=0$ and $c_{31}=0$. By Eq.~(\mref{eq:ncs3215}), we have $c_{23}^2=c_{13}c_{23}$. Thus we get $c_{23}=0$ or $c_{23}=a$.
Then  we obtain the solutions
$$N_{2,5}=\left(\begin{matrix}
0&0&a\\
0&0&0\\
0&0&0\\
\end{matrix}\right)\quad\text{and}
\quad
N_{2,6}=\left(\begin{matrix}
0&0&a\\
0&0&a\\
0&0&0\\
\end{matrix}\right)
\quad(a\neq 0).$$
\end{enumerate}
\item
Let $c_{11}+c_{13}=0$. Then $a+b=0$. So $b=-a\neq0$. Thus $c_{11}=-a$ and $c_{33}=-c_{11}=a$. By Eq.~(\mref{eq:ncs321}), we have $c_{31}=-\frac{c_{11}^2}{c_{13}}=-a.$
Since $c_{21}=-\frac{b^2}{a^2}c_{23}
+\frac{b^2}{a}+b$, we have $c_{21}=-c_{23}$. Let $c_{21}=c$, where $c\in\bfk$. Then we obtain the Rota-Baxter operators
$$N_{2,7}=\left(
\begin{matrix}
-a&0&a\\
c&0&-c\\
-a&0&a\\
\end{matrix}\right)\quad(a, c\in\bfk, a\neq 0).$$
\end{enumerate}

\smallskip
\noindent
{\bf Case 2. $c_{12}\neq 0$:} Denote $b=c_{12}$. Then $b\neq 0$. Also denote $a=c_{13}$.
By Eqs.~(\mref{eq:ncs3211}) and ~(\mref{eq:ncs3217}), we get $c_{12}(c_{11}+c_{13}-c_{31}-c_{33})=0$. Thus
\begin{equation}
c_{11}+c_{13}=c_{31}+c_{33}.
\mlabel{eq:ncs32im1}
\end{equation}
There are two subcases to consider.

\smallskip
\noindent
{\bf Subcase 1. $a\neq 0$:} Then Eqs.~(\mref{eq:ncs323}) and ~(\mref{eq:ncs326}) give $c_{12}(c_{23}-c_{13})=0$. So $c_{23}=c_{13}=a$. By Eq.~(\mref{eq:ncs323}), we have $c_{11}+c_{12}+c_{33}=0.$  Using Eqs.~(\mref{eq:ncs321}) and~(\mref{eq:ncs323}), we obtain $c_{11}^2+c_{11}c_{13}+c_{13}(c_{31}+
c_{33})=0$. So by Eq.~(\mref{eq:ncs32im1}), we get $(c_{11}+c_{13})^2=0.$  Thus $c_{11}+c_{13}=0$. Then we get $c_{11}=-c_{13}=-a$, and so $c_{33}=-c_{11}-c_{12}=a-b$.
Since $c_{31}+c_{33}=c_{11}+c_{13}=0$, we have $c_{31}=b-a$. By Eq.~(\mref{eq:ncs3210}), we get $c_{13}(c_{21}+c_{22})=c_{23}c_{31}$. Since $c_{13}=c_{23}$, we have $c_{21}+c_{22}=c_{31}=b-a$. By Eq.~(\mref{eq:ncs3222}), we have
$$c_{32}=\frac{c_{31}(c_{11}+c_{33})}{c_{23}}
=\frac{b(a-b)}{a}.$$
Thus Eq.~(\mref{eq:ncs3214}) gives $c_{22}=b$.  So $c_{21}=-a$. Then we obtain the solutions
$$N_{2,8}=\left(
\begin{matrix}
-a&b&a\\
-a&b&a\\
b-a&\frac{b(a-b)}{a}&a-b\\
\end{matrix}\right)\quad(a\neq 0,b\neq0).$$

\smallskip
\noindent
{\bf Subcase 2 $a= 0$.} Then $c_{13}=0$. By Eqs.~(\mref{eq:ncs321}) and~(\mref{eq:ncs3227}), we have $c_{11}=c_{33}=0$. By Eq.~(\mref{eq:ncs326}), we have $c_{12}c_{23}=0$, and so $c_{23}=0$.
By Eq.~(\mref{eq:ncs32im1}),
we have $c_{31}=0$. Thus Eq.~(\mref{eq:ncs3213}) gives $c_{21}=0$. So Eq.~(\mref{eq:ncs3214}) gives $c_{22}=0$. By Eq.~(\mref{eq:ncs3226}), we have $c_{32}=0$.
Thus we obtain the solutions
$$N_{2,9}=\left(
\begin{matrix}
0&b&0\\
0&0&0\\
0&0&0\\
\end{matrix}\right)\quad(b\in\bfk\setminus\{0\}).$$

\subsubsection{The proof for $\bfk[NCS(3)]$}
\mlabel{ss:ncs3}
We now show that the matrices of the Rota-Baxter operators $P$ on the semigroup algebra $\bfk[NCS(3)]$ are given by $N_{3,i}, 1\leq i\leq 8,$ in Table~\mref{ta:rbncs1}.

Applying the Cayley table of $NCS(3)$ in Eqs.~(\mref{eq:3eijm}) and then
taking $i=j=1$ with $1\leq m\leq 3$; $i=1,j=2$ with $m=1$; $i=1,j=3$ with $m=2$;
$i=2,j=1$ with $m=2$; $i=2, j=2$ with $m=1, 2$; $i=2, j=3$ with $m=1, 2$; $i=3,j=1$ with $m=2$; $i=3,j=2$ with $2\leq m\leq 3$ and $i=3,j=3$ with $m=1$, we obtain
\begin{eqnarray}
&c_{11}=0,\mlabel{eq:ncs331}\\
&c_{12}^2+c_{12}c_{13}=0,
\mlabel{eq:ncs332}\\
&c_{13}(c_{12}+c_{13})=0,
\mlabel{eq:ncs333}\\
&c_{12}c_{21}+c_{13}c_{31}=0,
\mlabel{eq:ncs334}\\
&c_{12}(c_{31}+c_{22})+c_{13}c_{32}
=0,\mlabel{eq:ncs338}\\
&c_{12}(c_{21}+c_{22}+c_{23})=0,
\mlabel{eq:ncs3310}\\
&c_{21}(c_{21}+c_{23})=c_{23}c_{31},
\mlabel{eq:ncs3312}\\
&c_{22}^2 +2c_{12}c_{21}+c_{23}c_{32}=0,
\mlabel{eq:ncs3313}\\
&c_{21}c_{22}=c_{31}(c_{21}+c_{22}),
\mlabel{eq:ncs3315}\\
&c_{22}^2+c_{12}(c_{21}+c_{31})+c_{23}
c_{32}=0,\mlabel{eq:ncs3316}\\
&c_{12}(c_{31}+c_{32}+c_{33})=0,
\mlabel{eq:ncs3318}\\
&c_{32}(c_{22}+c_{33})+c_{12}(c_{31}+
c_{21})=0,\mlabel{eq:ncs3321}\\
&c_{33}^2+c_{13}(c_{31}+c_{21})+c_{32}
c_{23}=0,\mlabel{eq:ncs3322}\\
&c_{31}^2+c_{31}c_{32}=c_{32}c_{21}.
\mlabel{eq:ncs3323}
\end{eqnarray}
First note that Eqs.~(\mref{eq:ncs331})-(\mref{eq:ncs333})
give $c_{11}=0$ and $c_{12}+c_{13}=0$.
We consider the two cases when $c_{12}=0$ and $c_{12}\neq 0$.
\smallskip

\noindent
{\bf Case 1. $c_{12}\neq 0$:} Denote $a=c_{12}$. Then we get $c_{13}=-c_{12}=-a$. By Eq.~(\mref{eq:ncs334}), we get $c_{21}=c_{31}$. Thus Eq.~(\mref{eq:ncs3312}) gives $c_{21}^2+c_{21}c_{23}=c_{23}c_{21}.$
So $c_{21}=c_{31}=0$. Denote $b=c_{22}$. Then by Eq.~(\mref{eq:ncs3310}), we get $c_{22}+c_{23}=0$. Thus $c_{23}=-b$.
By Eq.~(\mref{eq:ncs3318}), we have $c_{32}+c_{33}=0$. From Eq.~(\mref{eq:ncs338}), we obtain $c_{12}c_{22}+c_{13}c_{32}=0$. Thus we get $c_{22}=c_{32}$. So $c_{32}=b$ and  then $c_{33}=-b$. Then we obtain the solutions
$$N_{3,1}=\left(\begin{matrix}
0&a&-a\\
0&b&-b\\
0&b&-b\\
\end{matrix}\right)\quad(a, b\in\bfk, a\neq 0).$$

\smallskip
\noindent
{\bf Case 2 $c_{12}=0$.} Then $c_{13}=0$. Let $b=c_{23}$. There are two subcases to consider.

\smallskip
\noindent
{\bf Subcase 1. $c_{32}\neq 0$:} Denote $a=c_{32}$. Then $a\neq 0$.  By Eq.~(\mref{eq:ncs3321}), we have $c_{32}(c_{22}+c_{33})=0$.
Since $c_{32}\neq 0$, we have $c_{22}=-c_{33}$. By Eq.~(\mref{eq:ncs3313}), we have $c_{22}^2=-c_{23}c_{32}=-ab.$
Thus $c_{22}=\pm\sqrt{ab}i$ and $c_{33}=\mp\sqrt{ab}i$.
\begin{enumerate}
\item
Let $c_{21}=0$. Then by Eq.~(\mref{eq:ncs3312}), we have $c_{23}c_{31}=0$.
\begin{enumerate}
\item
Let $c_{23}\neq 0$. Then $b\neq0$. So we have $c_{31}=0$. Thus we obtain the solutions
$$N_{3,2}=\left(
\begin{matrix}
0&0&0\\
0&\sqrt{ab}i&b\\
0&a&-\sqrt{ab}i\\
\end{matrix}\right)\quad \text{and}\quad
N_{3,3}=\left(
\begin{matrix}
0&0&0\\
0&-\sqrt{ab}i&b\\
0&a&\sqrt{ab}i\\
\end{matrix}\right)\quad (a\neq0,b\neq0).$$
\item
Let $c_{23}=0$. Then $b=0$. By Eq.~(\mref{eq:ncs3316}), we have $c_{22}=0$. Since $c_{13}=0$ and $c_{23}=0$, Eqs.~(\mref{eq:ncs3312}) and ~(\mref{eq:ncs3322}) give $c_{21}=0$ and $c_{33}=0$. By Eq.~(\mref{eq:ncs3323}), we have $c_{31}(c_{31}+a)=0$. Thus $c_{31}=0$ or $c_{31}=-a$.  So we get solutions
$$N_{3,4}=\left(\begin{matrix}
0&0&0\\
0&0&0\\
0&a&0\\
\end{matrix}\right)\quad (a\neq 0)\quad \text{and}\quad N_{3,5}=\left(\begin{matrix}
0&0&0\\
0&0&0\\
-a&a&0\\
\end{matrix}\right)\quad (a\neq0).
$$
\end{enumerate}
\item
Let $c_{21}\neq 0$. Denote $c=c_{21}$. If  $c_{23}=0$, then by Eq.~(\mref{eq:ncs3312}), we have $c_{21}=0$, a contradiction. Thus we have $c_{23}=b\neq 0$. Applying Eq.~(\mref{eq:ncs3312}) again, we obtain $c_{31}=\frac{c(c+b)}{b}$.
By Eq.~(\mref{eq:ncs3315}), we have $c_{22}(c_{21}-c_{31})=c_{21}c_{31}$.
Since $c_{21}-c_{31}=-\frac{c^2}{b}\neq 0$, we have  $c_{22}=-(b+c)$. Since $c_{22}=-c_{33}$, we have $c_{33}=b+c$.
Then by $c_{12}=0$,  Eq.~(\mref{eq:ncs3313}) gives $c_{32}=-\frac{c_{22}^2}{c_{23}}
=-\frac{(b+c)^2}{b}$. By our assumption that $c_{32}\neq 0$, we have $b+c\neq 0$.
Thus we get the solutions
$$N_{3,6}=\left(
\begin{matrix}
0&0&0\\
c&-(c+b)&b\\
\frac{c(c+b)}{b}&-\frac{(c+b)^2}{b}&c+b\\
\end{matrix}\right)\quad (a\neq 0, b\neq0, b+c\neq0).$$
\end{enumerate}

\smallskip
\noindent
{\bf Subcase 2. $c_{32}=0$:} Then Eq.~(\mref{eq:ncs3313}) gives $c_{22}=0$. Then by Eq.~(\mref{eq:ncs3322}) and~(\mref{eq:ncs3323}), we have $c_{31}=c_{33}=0$. Thus Eq.~(\mref{eq:ncs3312}) gives $c_{21}(c_{21}+b)=0$. Thus $c_{21}=0$ or $c_{21}=-b$. So we get the solutions
$$N_{3,7}=\left(\begin{matrix}
0&0&0\\
0&0&b\\
0&0&0\\
\end{matrix}\right)\quad (b\in\bfk)
\quad\text{and}\quad N_{3,8}=
\left(\begin{matrix}
0&0&0\\
-b&0&b\\
0&0&0\\
\end{matrix}\right)\quad (b\neq0).$$

\subsubsection{The proof for $\bfk[NCS(4)]$}
\mlabel{ss:ncs4}
We prove that the matrices of the Rota-Baxter operators the semigroup algebra $\bfk[NCS(4)]$ are given by $N_{4,i}, 1\leq i\leq 8, $ in Table~\mref{ta:rbncs1}.

Applying the Cayley table of $NCS(4)$ in Eqs.~(\mref{eq:3eijm}) and then
taking $i=j=1$ with $1\leq m\leq 3$; $i=1,j=2$ with $m=1,3$; $i=1,j=3$ with $m=1$; $i=2, j=1$ with $1\leq m\leq 2$;
$i=2, j=2$ with $m= 3$; $i=2, j=3$ with $m=3$; $i=3,j=1$ with $1\leq m\leq 3$ and $i=3,j=2$ with $1\leq m\leq 2$, we obtain
\begin{eqnarray}
&c_{13}(c_{11}+
c_{12}+c_{13})=c_{11}(c_{11}+c_{13})
+c_{12}c_{21}
,\mlabel{eq:ncs341}\\
&c_{12}(c_{11}+c_{13}+c_{22})=0,
\mlabel{eq:ncs342}\\
&c_{13}(2c_{11}+c_{12}+2c_{13})+c_{12}
c_{23}=0,\mlabel{eq:ncs343}\\
&c_{13}(c_{21}+c_{22}+c_{23})=c_{11}(
c_{11}+c_{13})+c_{12}c_{21},
\mlabel{eq:ncs344}\\
&c_{13}(c_{11}+c_{13}+c_{21}+c_{22}+c_{23}
)+c_{12}c_{23}=0,\mlabel{eq:ncs346}\\
&c_{13}(c_{31}+c_{32}+c_{33})=c_{11}
(c_{11}+c_{13})+c_{12}c_{21},\mlabel{eq:ncs347}\\
&c_{23}(c_{12}+c_{13})=c_{21}(c_{11}+
c_{22}),\mlabel{eq:ncs3410}\\
&c_{22}^2+c_{12}(c_{21}+c_{23})=0,
\mlabel{eq:ncs3411}\\
&c_{13}(c_{21}+c_{23})+
c_{23}(c_{21}+2c_{22}+c_{23})=0,
\mlabel{eq:ncs3415}\\
&c_{13}(c_{21}+c_{23})+c_{23}(c_{22}+
c_{31}+c_{32}+c_{33})=0,
\mlabel{eq:ncs3418}\\
&(c_{31}+c_{33})(c_{12}+c_{13})=c_{11}(
c_{11}+c_{12}+c_{13})+c_{32}c_{21},
\mlabel{eq:ncs3419}\\
&c_{32}(c_{11}+c_{12}+c_{13})=
c_{12}(c_{31}+c_{33}+c_{11}+c_{12}+c_{13}
)+c_{32}c_{22},\mlabel{eq:ncs3420}\\
&c_{13}(c_{31}+c_{33}+c_{11}+c_{12}+c_{13}
)+c_{32}c_{23}=0,\mlabel{eq:ncs3421}\\
&(c_{31}+c_{33})(c_{21}+c_{22}+c_{23})=
c_{11}(c_{31}+c_{33}+c_{21}+c_{22}+c_{23}
)+c_{32}c_{21},\mlabel{eq:ncs3422}\\
&c_{32}(c_{21}+c_{22}+c_{23})=c_{12}(c_{31}
+c_{33}+c_{21}+c_{22}+c_{23})+c_{32}c_{22}
.\mlabel{eq:ncs3423}
\end{eqnarray}
We divide the proof into two cases.

\smallskip
\noindent
{\bf Case 1. $c_{12}=0$:} Then Eq.~(\mref{eq:ncs3411}) gives $c_{22}=0$. By Eq.~(\mref{eq:ncs341}), we have $c_{11}^2=c_{13}^2$. Thus by Eq.~(\mref{eq:ncs343}), we get $(c_{11}+c_{13})^2=0$. So $c_{11}+c_{13}=0$.
There are two subcases to consider.

\smallskip
\noindent
{\bf Subcase 1. $c_{13}\neq 0$:}
Denote $a=c_{13}$. Then $a\neq 0$. So $c_{11}=-c_{13}=-a$.
By Eqs.~(\mref{eq:ncs344}), ~(\mref{eq:ncs346}) and ~(\mref{eq:ncs347}), we get $0=c_{11}+c_{13}=c_{21}+c_{23}=c_{31}+c_{32}
+c_{33}$. Let $c_{21}=b$, where $b\in\bfk$. Then $c_{23}=-b$. So Eq.~(\mref{eq:ncs3422}) gives $c_{11}(c_{31}+c_{33})+c_{32}c_{21}=0$, and then $c_{32}(c_{21}-c_{11})=0$. Thus we have $(a+b)c_{32}=0$.
Subdividing further, we have
\begin{enumerate}
\item
Suppose $a+b\neq 0$.  Then $c_{32}=0$. So $c_{31}+c_{33}=0$. Let $c_{31}=c$, where $c\in\bfk$. Then $c_{33}=-c$. Thus we get the solutions
$$N_{4,1}=\left(
\begin{matrix}
-a&0&a\\
b&0&-b\\
c&0&-c\\
\end{matrix}\right)\quad (a, b\in \bfk, a\neq 0, a+b\neq 0).$$
\item
Suppose $a+b=0$. Let $c_{31}=c$, and let $c_{32}=d$, where $c,d\in\bfk$. Then we get $c_{33}=-c-d$.  Thus we obtain the solutions
$$
N_{4,2}=\left(\begin{matrix}
-a&0&a\\
-a&0&a\\
c&d&-c-d\\\end{matrix}\right)
\quad (a, c, d\in\bfk, a\neq 0).$$
\end{enumerate}
\smallskip

\noindent
{\bf Subcases 2. $c_{13}=0$:} Then we have $c_{11}=0$. Thus Eq.~(\mref{eq:ncs3415}) gives $c_{23}(
c_{21}+c_{23})=0$. Using Eq.~(\mref{eq:ncs3423}), we have $c_{32}(c_{21}+c_{23})=0$. By Eq.~(\mref{eq:ncs3418}), we get $c_{23}(c_{31}+c_{32}+c_{33})=0$. By Eqs.~(\mref{eq:ncs3419}) and ~(\mref{eq:ncs3421}), we have $c_{32}c_{31}=0$ and $c_{32}c_{23}=0$.
\begin{enumerate}
\item
Suppose $c_{23}\neq 0$. Denote $a=c_{23}$. Then $a\neq 0$. Thus $c_{21}+c_{23}=0$. So $c_{21}=-a$ and $c_{32}=0$. Thus $c_{31}+c_{33}=0$. Let $c_{31}=b$, where $ b\in\bfk$. So $c_{33}=-b$. Thus we obtain the solutions
$$N_{4,3}=\left(
\begin{matrix}
0&0&0\\
-a&0&a\\
b&0&-b\\
\end{matrix}\right)\quad (a, b\in\bfk, a\neq 0).$$
\item
Let $c_{23}=0$.
Then Eqs.~(\mref{eq:ncs341})-(\mref{eq:ncs3423}) are equivalent to the following system of equations.
$$\left\{
\begin{array}{cc}
c_{32}c_{21}=0,\\
c_{32}(c_{31}+c_{32}+c_{33})=0,\\
(c_{31}+c_{33})(c_{31}+c_{32}+c_{33})
=0.\\
\end{array}\right.$$
Denote $b=c_{31}$.
\begin{enumerate}
\item
When $c_{32}\neq 0$, then $c_{21}=0$ and $c_{31}+c_{32}+c_{33}=0$.
  Take $c=c_{32}$. Then $c_{33}=-c-b$. Then we get the solutions
$$N_{4,4}=\left(\begin{matrix}
0&0&0\\
0&0&0\\
b&c&-b-c\\
\end{matrix}\right)\quad (b, c\in\bfk, c\neq 0).$$
\item
When $c_{32}=0$, then $c_{31}+c_{33}=0$. Then $c_{33}=-b$. Denote $a=c_{21}$. Then we get the solutions
$$N_{4,5}=\left(\begin{matrix}
0&0&0\\
a&0&0\\
b&0&-b\\
\end{matrix}\right)\quad (a, b\in\bfk).$$
\end{enumerate}
\end{enumerate}
\smallskip

\noindent
{\bf Case 2. $c_{12}\neq 0$:} Then Eq.~(\mref{eq:ncs342}) gives $c_{11}+c_{13}+c_{22}=0$. By Eqs.~(\mref{eq:ncs341}), (\mref{eq:ncs344}) and ~(\mref{eq:ncs347}), we have $$c_{13}(c_{11}+c_{12}+c_{13})=c_{13}(c_{21}
+c_{22}+c_{23})=c_{13}(c_{31}+c_{32}+c_{
33}).$$
We divide this case into two subcases.

\smallskip
\noindent
{\bf Subcase 1. $c_{13}\neq 0$:} Denote $b=c_{13}$ and $a=c_{11}$. Then we get $c_{22}=-a-b$ and $c_{11}+c_{12}+c_{13}=c_{21}+c_{22}+c_{23}
=c_{31}+c_{32}+c_{33}$. By Eq.~(\mref{eq:ncs343}), we get $c_{13}(c_{11}+c_{12}+c_{13})=-(c_{13}(c_{11}
+c_{13})+c_{12}c_{23})$.  In Eq.~(\mref{eq:ncs341}), replacing $c_{13}(c_{11}+c_{12}+c_{13})$ by $-(c_{13}(c_{11}+c_{13})+c_{12}c_{23})$, we obtain $(c_{11}+c_{13})^2+c_{12}
(c_{21}+c_{23})=0$. From $c_{11}+c_{12}+c_{13}=c_{21}+c_{22}
 +c_{23}$ and Eq.~(\mref{eq:ncs342}), we get $(c_{11}+c_{12}+c_{13})^2=0$. So $c_{11}+c_{12}+c_{13}=0$. Thus we have $c_{12}=-(a+b)\neq 0$ and $c_{21}+c_{23}=-c_{22}=a+b$. Thus
 Eqs.~(\mref{eq:ncs3410}) and ~(\mref{eq:ncs3411}) give $c_{21}=a$ and $c_{23}=b$.
 Denote $c=c_{31}$ and $d=c_{32}$. Since $c_{31}+c_{32}+c_{33}=0$,  we have $c_{33}=-c-d$. Then we obtain the following solutions of Eqs.~(\mref{eq:ncs341})-(\mref{eq:ncs3423})
 $$N_{4,6}=\left(\begin{matrix}
 a&-(a+b)& b\\
 a&-(a+b)&b\\
 c&d&-(c+d)\\
 \end{matrix}\right)\quad (a,b,c,d\in\bfk, b\neq 0, a+b\neq 0).$$

 \smallskip
 \noindent
 {\bf Subcase 2. $c_{13}=0$:} Then $c_{11}+c_{22}=0$. Denote $b=c_{12}$. Then $0\neq b\in\bfk$. Also denote $a=c_{11}$. Then $c_{22}=-a$ and Eq.~(\mref{eq:ncs341}) gives $c_{21}
 =-\frac{a^2}{b}$. By Eq.~(\mref{eq:ncs343}), we have $c_{23}=0$.
 \begin{enumerate}
 \item When $a+b\neq 0$, then by
 Eq.~(\mref{eq:ncs3419}), we get \begin{equation}
 b(c_{31}+c_{33})=a(a+b)-\frac{a^2}{b}
 c_{32}.\mlabel{eq:ncs34im1}
 \end{equation}
 By Eq.~(\mref{eq:ncs3420}), we have

 \begin{equation}
 b(c_{31}+c_{33})=(2a+b)c_{32}-b(a+b).
 \mlabel{eq:ncs34im2}
 \end{equation}
From Eqs.~(\mref{eq:ncs34im1}) and~(\mref{eq:ncs34im2}),
 we get $c_{32}=b$ and so $c_{31}+c_{33}=a$.  Denote $c=c_{31}$. Then $c_{33}=a-c$. So we get the solutions
 $$N_{4,7}=\left(
 \begin{matrix}
 a&b&0\\
 -\frac{a^2}{b}&-a&0\\
 c&b&a-c\\
 \end{matrix}\right)\quad (a,b,c\in\bfk, a+b\neq 0, b\neq 0).$$
 \item
When $a+b=0$, then we have $b=-a\neq 0$ and $c_{11}+c_{12}=0$.  By Eq.~(\mref{eq:ncs3419}), we have $c_{12}(c_{31}+c_{33})=c_{32}c_{21}$. Thus $c_{32}=-(c_{31}+c_{33})$. Denote $c=c_{31}$ and $d=c_{32}$. Then $c_{33}=-d-c$.  So we obtain the solutions
 $$N_{4,8}=\left(\begin{matrix}
 -b&b&0\\
 -b&b&0\\
 c&d&-c-d\\
 \end{matrix}\right)\quad (b,c,d\in\bfk, b\neq 0).$$
 \end{enumerate}

\subsubsection{The proof for $\bfk[NCS(5)]$}
\mlabel{ss:ncs5}
We next prove that the matrices of the Rota-Baxter operators on the semigroup algebra $\bfk[NCS(5)]$ are given by $N_{5,i}, 1\leq i\leq 22,$ in Table~\mref{ta:rbncs1}.

Applying the Cayley table of $NCS(5)$ in Eqs.~(\mref{eq:3eijm}) and then
taking $i=j=1$ with $1\leq m\leq 3$; $i=2,j=1$ with $1\leq m\leq 3$ and $i=3,j=1$ with $1\leq m\leq 3$, we obtain
\begin{eqnarray}
&c_{11}^2+c_{12}c_{21}+c_{13}c_{31}=0,
\mlabel{eq:ncs351}\\
&c_{11}c_{12}+c_{12}c_{22}+
c_{13}c_{32}=0,
\mlabel{eq:ncs352}\\
&c_{11}c_{13}+c_{12}c_{23}+c_{13}c_{33}=0, \mlabel{eq:ncs353}\\
&c_{21}c_{11}+c_{21}c_{22}+c_{23}c_{31}
=0,\mlabel{eq:ncs354}\\
&c_{21}c_{12}+c_{22}^2+c_{23}c_{32}=0,
\mlabel{eq:ncs355}\\
&c_{21}c_{13}+c_{22}c_{23}+c_{23}c_{33}
=0,\mlabel{eq:ncs356}\\
&c_{11}c_{31}+c_{21}c_{32}+c_{31}c_{33}
=0,\mlabel{eq:ncs357}\\
&c_{12}c_{31}+c_{22}c_{32}+c_{32}c_{33}
=0,\mlabel{eq:ncs358}\\
&c_{13}c_{31}+c_{23}c_{32}+c_{33}^2=0.
\mlabel{eq:ncs359}
\end{eqnarray}
From Eq.~(\mref{eq:ncs359}) we have $c_{33}^2=-(c_{13}c_{31}+c_{23}c_{32})$.
Denote $a=c_{13}$, $b=c_{32}$, $c=c_{31}$ and $d=c_{23}$. Then $c_{33}^2=-ac-bd$. So we have $c_{33}=\pm\sqrt{-ac-bd}$. Denote $F:=\sqrt{-ac-bd}$. So $F^2=-ac-bd$ and $c_{33}=\pm F$.

We divide the proof into five cases.

\smallskip
\noindent
{\bf Case 1. $a, b, c$ and $d$ are nonzero:} Then $c_{13}, c_{32}, c_{31}$ and $c_{23}$ are nonzero.  By Eqs.~(\mref{eq:ncs353}) and ~(\mref{eq:ncs357}), we have $$c_{31}c_{23}c_{12}=-c_{31}(c_{11}c_{13}
+c_{13}c_{33})=-c_{13}(c_{11}c_{31}
+c_{31}c_{33})=c_{13}c_{32}c_{21}.$$
Thus
\begin{equation}
c_{21}=\frac{cd}{ab}c_{12}.
\mlabel{eq:ncs35im1}
\end{equation}
Then by Eqs.~(\mref{eq:ncs357}) and ~(\mref{eq:ncs35im1}), we obtain \begin{equation}c_{11}=-c_{33}-\frac{c_{32}}{c_{31}}c_{21}
=-c_{33}-\frac{b}{c}c_{21}=-c_{33}-\frac{d}{a}c_{12}.
\mlabel{eq:ncs35im2}
\end{equation}
Applying Eqs.~(\mref{eq:ncs356}) and ~(\mref{eq:ncs35im1}), we have
\begin{equation}
c_{22}=-c_{33}-\frac{c_{13}}{c_{23}}c_{21}
=-c_{33}-\frac{a}{d}c_{21}=-c_{33}-\frac{c}{b}c_{12}.
\mlabel{eq:ncs35im3}
\end{equation}
In Eq.~(\mref{eq:ncs352}), by replacing $c_{11}$ and $c_{22}$ by $-c_{33}-\frac{d}{a}c_{12}$ and $-c_{33}-\frac{c}{b}c_{12}$, respectively, we get
\begin{equation}
c_{33}^2c_{12}^2-2abc_{33}c_{12}+(ab)^2=0.
\mlabel{eq:ncs35im4}
\end{equation}
If $ac+bd=0$, then $c_{33}=0$. Thus by Eq.~(\mref{eq:ncs35im4}), $ab=0$, a contradiction. Thus we have $ac+bd\neq 0$.  So $c_{33}\neq 0$ and $F\neq 0$. By Eq.~(\mref{eq:ncs35im4}) again, we get $(c_{33}c_{12}-ab)^2=0.$ Thus $c_{12}=\frac{ab}{c_{33}}$. Since $c_{33}=\pm F$,  $c_{12}=\pm\frac{ab}{F}$. Then by Eq.~(\mref{eq:ncs35im1}), $c_{21}=\pm\frac{cd}{F}$. By Eqs.~(\mref{eq:ncs35im2}) and ~(\mref{eq:ncs35im3})  and $F^2=-ac-bd$, we get $c_{11}=\pm\frac{ac}{F}$ and $c_{22}=\pm\frac{bd}{F}$. Thus we obtain the solutions
$$N_{5,1}=\left(\begin{matrix}
\frac{ac}{F}&\frac{ab}{F}&a\\
\frac{cd}{F}&\frac{bd}{F}&d\\
c&b&F\\
\end{matrix}\right)\quad\text{and}
\quad N_{5,2}=\left(\begin{matrix}
-\frac{ac}{F}&-\frac{ab}{F}&a\\
-\frac{cd}{F}&-\frac{bd}{F}&d\\
c&b&-F\\
\end{matrix}\right)\quad (a,b,c,d\in\bfk\setminus \{0\}).$$
\delete{
\smallskip
\noindent
{\bf Subcase 1.
 $c_{33}=F$:} Then $c_{12}=\frac{ab}{F}$. So by Eq.~(\mref{eq:ncs35im2}) and $F^2=-ac-bd$, we have $$c_{11}=-F-\frac{d}{a} \frac{ab}{F}=-\frac{F^2+bd}{F}=\frac{ac}{F}.$$
By Eq.~(\mref{eq:ncs35im3}), we get
$c_{22}=\frac{bd}{F}.$
Thus we obtain the solutions
$$N_{5,1}=\left(\begin{matrix}
\frac{ac}{F}&\frac{ab}{F}&a\\
\frac{cd}{F}&\frac{bd}{F}&d\\
c&b&F\\
\end{matrix}\right)(a,b,c,d\in\bfk\setminus \{0\}).$$

\smallskip
\noindent
{\bf Subcase 2. $c_{33}=-F$:}
Then we have $c_{12}=-\frac{ab}{F}.$ Similarly, by Eqs.~(\mref{eq:ncs35im2}) and ~(\mref{eq:ncs35im3}), we get $c_{11}=-\frac{ac}{F}$ and $c_{22}=-\frac{bd}{F}$.  Thus we obtain the solutions
$$N_{5,2}=\left(\begin{matrix}
-\frac{ac}{F}&-\frac{ab}{F}&a\\
-\frac{cd}{F}&-\frac{bd}{F}&d\\
c&b&-F\\
\end{matrix}\right)(a,b,c,d\in\bfk\setminus \{0\}).$$
}
\smallskip
\noindent
{\bf Case 2. Exactly one of $a,b,c,d$ is $0$:}  Then there are four subcases to consider.

\smallskip
\noindent
{\bf Subcase 1, $a=0$ and $b,c,d\neq 0$:} Then $c_{13}=0$, $c_{32}\neq 0, c_{31}\neq0$, and $c_{23}\neq0$.  By Eq.~(\mref{eq:ncs353}), we have $c_{12}c_{23}=0$. So $c_{12}=0$. Thus by Eq.~(\mref{eq:ncs351}), $c_{11}=0$. Then Eq.~(\mref{eq:ncs356}) gives $c_{22}+c_{33}=0$ and Eq.~(\mref{eq:ncs359}) gives $c_{33}=\pm\sqrt{bd}i$. So $c_{22}=\mp\sqrt{bd}i$. By Eq.~(\mref{eq:ncs354}), we have $c_{21}=-\frac{c_{23}c_{31}}{c_{22}}$.
So $c_{21}=\mp\frac{c\sqrt{d}i}{\sqrt{b}}$.
Thus we obtain the solutions
$$N_{5,3}=\left(
\begin{matrix}
0&0&0\\
-\frac{c\sqrt{d}i}{\sqrt{b}}&-\sqrt{bd}i
&d\\
c&b&\sqrt{bd}i\\
\end{matrix}\right)\quad \text{and}\quad
N_{5,4}=\left(
\begin{matrix}
0&0&0\\
\frac{c\sqrt{d}i}{\sqrt{b}}&\sqrt{bd}i
&d\\
c&b&-\sqrt{bd}i\\
\end{matrix}\right)\quad (a,b,c,d\in \bfk, b,c,d\neq 0).$$

\delete{
\begin{enumerate}
\item
When $c_{33}=\sqrt{bd}i$, then $c_{22}=-\sqrt{bd}i$. By Eq.~(\mref{eq:ncs354}), we have $c_{21}=-\frac{c_{23}c_{31}}{c_{22}}$.
So $c_{21}=-\frac{c\sqrt{d}i}{\sqrt{b}}$.
Thus we obtain the Rota-Baxter operators
$$N_{5,3}=\left(
\begin{matrix}
0&0&0\\
-\frac{c\sqrt{d}i}{\sqrt{b}}&-\sqrt{bd}i
&d\\
c&b&\sqrt{bd}i\\
\end{matrix}\right)\quad (b,c,d\in \bfk\setminus \{0\}).$$
\item
When $c_{33}=-\sqrt{bd}i$, we similarly obtain the solutions
$$N_{5,4}=\left(
\begin{matrix}
0&0&0\\
\frac{c\sqrt{d}i}{\sqrt{b}}&\sqrt{bd}i
&d\\
c&b&-\sqrt{bd}i\\
\end{matrix}\right)\quad (b,c,d\in \bfk\setminus \{0\}).$$
\end{enumerate}
}

\smallskip
\noindent
{\bf Subcase 2. $b=0$ and $a,c,d\neq 0$:} Then $c_{32}=0$, $c_{13}\neq 0$, $c_{31}\neq 0$ and $c_{23}\neq 0$. Then by Eq.~(\mref{eq:ncs358}), $c_{12}c_{31}=0$. So $c_{12}=0$ and then by Eq.~(\mref{eq:ncs355}), $c_{22}=0$. By Eq.~(\mref{eq:ncs353}), we have $c_{11}+c_{33}=0.$
Then Eq.~(\mref{eq:ncs359}) gives $c_{33}=\pm\sqrt{ac}i$. Thus $c_{11}=\mp\sqrt{ac}i$. By Eq.~(\mref{eq:ncs356}), we have $c_{21}=-\frac{c_{23}c_{33}}{c_{13}}.$
So $c_{21}=\mp\frac{d\sqrt{c}i}{\sqrt{a}}$.
Then we obtain the solutions
$$N_{5,5}=\left(\begin{matrix}
-\sqrt{ac}i&0&a\\
-\frac{d\sqrt{c}i}{\sqrt{a}}&0&d\\
c&0&\sqrt{ac}i\\
\end{matrix}\right)\quad \text{and}\quad N_{5,6}=\left(\begin{matrix}
\sqrt{ac}i&0&a\\
\frac{d\sqrt{c}i}{\sqrt{a}}&0&d\\
c&0&-\sqrt{ac}i\\
\end{matrix}\right)\quad (a,b,c,d\in \bfk, a, c, d\neq 0).$$

\smallskip
\noindent
{\bf Subcase 3. $c=0$ and $a,b,d\neq 0$:} Then By Eq.~(\mref{eq:ncs357}), we have $c_{32}c_{21}=0$. So $c_{21}=0$.  By Eq.~(\mref{eq:ncs351}), we have $c_{11}=0$. Thus Eq.~(\mref{eq:ncs358}) gives $c_{22}+c_{33}=0$ and Eq.~(\mref{eq:ncs359}) gives $c_{33}=\pm\sqrt{bd}i$. Thus $c_{22}=\mp\sqrt{bd}i$. By Eq.~(\mref{eq:ncs353}), we have $c_{12}=-\frac{c_{13}c_{33}}{c_{23}}
=\mp\frac{a\sqrt{b}i}{\sqrt{d}}$.
Then we obtain the solutions
$$N_{5,7}=\left(
\begin{matrix}
0&-\frac{a\sqrt{b}i}{\sqrt{d}}&a\\
0&-\sqrt{bd}i&d\\
0&b&\sqrt{bd}i\\
\end{matrix}\right)\quad\text{and}\quad
N_{5,8}=\left(
\begin{matrix}
0&\frac{a\sqrt{b}i}{\sqrt{d}}&a\\
0&\sqrt{bd}i&d\\
0&b&-\sqrt{bd}i\\
\end{matrix}\right)\quad (a,b,c,d\in \bfk, a,b,d\neq 0).$$

\smallskip
\noindent
{\bf Subcase 4. $d=0$ and $a,b,c\neq0$:} Then by Eq.~(\mref{eq:ncs356}), we have $c_{21}c_{13}=0$. So $c_{21}=0$. Thus by Eq.~(\mref{eq:ncs355}), we have $c_{22}=0$.  So Eq.~(\mref{eq:ncs353}) gives $c_{11}+c_{33}=0$. Furthermore, by Eq.~(\mref{eq:ncs359}), we have $c_{33}=\pm\sqrt{ac}i$. Thus $c_{11}=\mp\sqrt{ac}i$.
By Eq.~(\mref{eq:ncs358}), we get $c_{12}=-\frac{c_{32}c_{33}}{c_{31}}
=\mp\frac{b\sqrt{a}i}{\sqrt{c}}$.
Then  we obtain the solutions
$$N_{5,9}=\left(\begin{matrix}
-\sqrt{ac}i&-\frac{b\sqrt{a}i}{\sqrt{c}}&a\\
0&0&0\\
c&b&\sqrt{ac}i\\
\end{matrix}\right)\quad\text{and}\quad N_{5,10}=\left(\begin{matrix}
\sqrt{ac}i&\frac{b\sqrt{a}i}{\sqrt{c}}&a\\
0&0&0\\
c&b&-\sqrt{ac}i\\
\end{matrix}\right)\quad (a,b,c,d\in \bfk, a,b,c\neq 0).$$

\delete{
by Eq.~(\mref{eq:ncs358}), $c_{12}c_{31}=0$. So $c_{12}=0$ and then by Eq.~(\mref{eq:ncs355}), $c_{22}=0$. By Eq.~(\mref{eq:ncs353}), we have $c_{11}+c_{33}=0.$ Then Eq.~(\mref{eq:ncs359}) gives $c_{33}=\sqrt{ac}i$ or $-\sqrt{ac}i$.
\begin{enumerate}
\item
If $c_{33}=\sqrt{ac}i$, then by Eq.~(\mref{eq:ncs356}), we have $c_{21}=-\frac{c_{23}c_{33}}{c_{13}}.$
So $c_{21}=-\frac{d\sqrt{c}i}{\sqrt{a}}$.
Then we obtain the solutions
$$N_{5,5}=\left(\begin{matrix}
-\sqrt{ac}i&0&a\\
-\frac{d\sqrt{c}i}{\sqrt{a}}&0&d\\
c&0&\sqrt{ac}i\\
\end{matrix}\right)\quad (a, c, d\in \bfk\setminus \{0\}).$$
\item
If $c_{33}=-\sqrt{ac}i$, then we similarly obtain the solutions
$$N_{5,6}=\left(\begin{matrix}
\sqrt{ac}i&0&a\\
\frac{d\sqrt{c}i}{\sqrt{a}}&0&d\\
c&0&-\sqrt{ac}i\\
\end{matrix}\right)\quad (a,c,d\in \bfk\setminus \{0\}).$$
\end{enumerate}

\smallskip
\noindent
{\bf Subcase 3. $c=0$ and $a,b,d\neq 0$:}  Then By Eq.~(\mref{eq:ncs357}), we have $c_{32}c_{21}=0$. So $c_{21}=0$.  By Eq.~(\mref{eq:ncs351}), we have $c_{11}=0$. Thus Eq.~(\mref{eq:ncs358}) gives $c_{22}+c_{33}=0$ and Eq.~(\mref{eq:ncs359}) gives $c_{33}=\sqrt{bd}i$ or $-\sqrt{bd}i$.
\begin{enumerate}
\item
If $c_{33}=\sqrt{bd}i$, then $c_{22}=-\sqrt{bd}i$. By Eq.~(\mref{eq:ncs353}), we have $c_{12}=-\frac{c_{13}c_{33}}{c_{23}}
=-\frac{a\sqrt{b}i}{\sqrt{d}}$. Then we get the solutions
$$N_{5,7}=\left(
\begin{matrix}
0&-\frac{a\sqrt{b}i}{\sqrt{d}}&a\\
0&-\sqrt{bd}i&d\\
0&b&\sqrt{bd}i\\
\end{matrix}\right)\quad (a,b,d\in \bfk\setminus \{0\}).$$
\item
If $c_{33}=-\sqrt{bd}i$, then we similarly obtain
the solutions
$$N_{5,8}=\left(
\begin{matrix}
0&\frac{a\sqrt{b}i}{\sqrt{d}}&a\\
0&\sqrt{bd}i&d\\
0&b&-\sqrt{bd}i\\
\end{matrix}\right)\quad (a,b,d\in \bfk\setminus \{0\}).$$
\end{enumerate}

\smallskip
\noindent
{\bf Subcase 4. $d=0$ and $a,b,c\neq0$:} Then by Eq.~(\mref{eq:ncs356}), we have $c_{21}c_{13}=0$. So $c_{21}=0$. Thus by Eq.~(\mref{eq:ncs355}), we have $c_{22}=0$.  So Eq.~(\mref{eq:ncs353}) gives $c_{11}+c_{33}=0$. Furthermore, by Eq.~(\mref{eq:ncs359}), we have $c_{33}=\sqrt{ac}i$ or $-\sqrt{ac}i$.
\begin{enumerate}
\item
If
$c_{33}=\sqrt{ac}i$, then $c_{11}=-\sqrt{ac}i$. By Eq.~(\mref{eq:ncs358}), we get $c_{12}=-\frac{c_{32}c_{33}}{c_{31}}
=-\frac{b\sqrt{a}i}{\sqrt{c}}$. Then  we obtain the solutions
$$N_{5,9}=\left(\begin{matrix}
-\sqrt{ac}i&-\frac{b\sqrt{a}i}{\sqrt{c}}&a\\
0&0&0\\
c&b&\sqrt{ac}i\\
\end{matrix}\right)\quad (a,b,c\in \bfk\setminus \{0\}).$$
\item
If $c_{33}=-\sqrt{ac}i$, then we similarly obtain the solutions
$$N_{5,10}=\left(\begin{matrix}
\sqrt{ac}i&\frac{b\sqrt{a}i}{\sqrt{c}}&a\\
0&0&0\\
c&b&-\sqrt{ac}i\\
\end{matrix}\right)\quad (a,b,c\in \bfk\setminus \{0\}).$$
\end{enumerate}
}
\smallskip
\noindent
{\bf Case 3. Exactly two of $a,b,c,d$ are $0$:} There are six subcases to consider. But note that if $a=b=0$, $c\neq 0$ and $d\neq 0$, i.e. $c_{13}=c_{32}=0$, $c_{31}\neq 0$ and $c_{23}\neq0$, then by Eq.~(\mref{eq:ncs358}), we have $c_{12}c_{31}=0$. So $c_{12}=0$. Thus Eqs.~(\mref{eq:ncs351}) and ~(\mref{eq:ncs355}) give $c_{11}=0$ and $c_{22}=0$, respectively. Then by Eq.~(\mref{eq:ncs354}), we have $c_{23}c_{31}=cd=0$, a contradiction. Similarly, if $c=d=0$ and $a,b\neq 0$, then we can obtain $c_{13}c_{32}=ab=0$, a contradiction.  So there are four subcases left to consider.

\smallskip
\noindent
{\bf Subcase 1. $a=c=0$ and $b,d\neq 0$:} Then $c_{13}=c_{31}=0$, $c_{32}\neq0$ and $c_{23}\neq0$. Thus by Eq.~(\mref{eq:ncs353}), we have $c_{12}c_{23}=dc_{12}=0$. So $c_{12}=0$.  Thus by Eq.(\mref{eq:ncs351}), we have $c_{11}=0$.  Then Eq.~(\mref{eq:ncs358}) gives $c_{22}+c_{33}=0$. By Eq.~(\mref{eq:ncs359}), we have $c_{33}=\pm\sqrt{bd}i$.
Thus we get $c_{22}=\mp\sqrt{bd}i$. Then we  get the solutions
$$N_{5,11}=\left(\begin{matrix}
0&0&0\\
0&\sqrt{bd}i&d\\
0&b&-\sqrt{bd}i\\
\end{matrix}\right)
\text{ and }
N_{5,12}=\left(\begin{matrix}
0&0&0\\
0&-\sqrt{bd}i&d\\
0&b&\sqrt{bd}i\\
\end{matrix}\right)\quad (b, d\in \bfk\setminus \{0\}).$$

\smallskip
\noindent
{\bf Subcase 2. $a=d=0$ and $b,c\neq 0$:} Then $c_{13}=c_{23}=0$, $c_{32}\neq 0$ and $c_{31}\neq 0$.
Eq.~(\mref{eq:ncs359}) gives $c_{33}=0$. Then by Eq.~(\mref{eq:ncs357}), we obtain $c_{11}=-\frac{c_{32}}{c_{31}}c_{21}=
-\frac{b}{c}c_{21}$. From Eq.~(\mref{eq:ncs358}), we have $c_{12}=-\frac{c_{32}}{c_{31}}c_{22}=-\frac
{b}{c}c_{22}$.
\begin{enumerate}
\item
If $c_{21}=0$, then $c_{11}=0$. By Eq.~(\mref{eq:ncs355}) we have $c_{22}=0$. So $c_{12}=0$. Then we obtain the solutions
$$N_{5,13_1}=\left(\begin{matrix}
0&0&0\\
0&0&0\\
c&b&0\\
\end{matrix}\right)\quad (b, c\in\bfk\setminus\{0\}). $$
\item
If $c_{21}\neq 0$, then denote $e=c_{21}$.
Thus by $c_{11}=-\frac
{b}{c}c_{21}$, we have $c_{11}=-\frac{be}{c}$.
By Eq.~(\mref{eq:ncs354}), we have $c_{11}+c_{22}=0$. Thus $c_{22}=\frac{be}{c}$. So we have $c_{12}=-\frac{b}{c}c_{22}=-\frac{b^2e}{c^2}
$.Then we obtain the solutions $$N_{5,13_2}=\left(\begin{matrix}
-\frac{be}{c}&-\frac{b^2e}{c^2}&0\\
e&\frac{be}{c}&0\\
c&b&0\\
\end{matrix}\right)\quad (b,c, e\in\bfk\setminus\{0\}).
$$
\end{enumerate}
In summary, we obtain the solutions
$$N_{5,13}=\left(\begin{matrix}
-\frac{be}{c}&-\frac{b^2e}{c^2}&0\\
e&\frac{be}{c}&0\\
c&b&0\\
\end{matrix}\right)\quad (b,c\in\bfk\setminus \{0\}, e\in\bfk).
$$

\smallskip
\noindent
{\bf Subcase 3. $b=c=0$ and $a,d\neq 0$:} This subcase is similar to the above Subcase 2.
Denote $e=c_{21}$. Then we can obtain
the solutions
$$N_{5,14}=\left(\begin{matrix}
-\frac{ae}{d}&-\frac{a^2e}{d^2}&a\\
e&-\frac{ae}{d}&d\\
0&0&0\\
\end{matrix}\right)\quad (a, d, e\in \bfk, a,d\neq 0).
$$

\smallskip
\noindent
{\bf Subcase 4. $b=d=0$ and $a,c\neq 0$:} By Eq.~(\mref{eq:ncs359}),  $c_{33}=\pm\sqrt{ac}i$. Then by Eq.~(\mref{eq:ncs353}), we have $c_{11}=\mp\sqrt{ac}i$.
Then we can obtain the solutions

$$N_{5,15}=\left(\begin{matrix}
-\sqrt{ac}i&0&a\\
0&0&0\\
c&0&\sqrt{ac}i\\
\end{matrix}\right) \text{ and }
N_{5,16}=\left(\begin{matrix}
\sqrt{ac}i&0&a\\
0&0&0\\
c&0&-\sqrt{ac}i\\
\end{matrix}\right)\quad (a,c\in \bfk\setminus \{0\}).
$$

\delete{
Then $c_{13}=c_{23}=0$, $c_{32}\neq 0$ and $c_{31}\neq 0$.
Eq.~(\mref{eq:ncs359}) gives $c_{33}=0$. Then by Eq.~(\mref{eq:ncs357}), we obtain $c_{11}=-\frac{c_{32}}{c_{31}}c_{21}=
-\frac{b}{c}c_{21}$. From Eq.~(\mref{eq:ncs358}), we have $c_{12}=-\frac{c_{32}}{c_{31}}c_{22}=-\frac
{b}{c}c_{21}$. Denote $e=c_{21}$.
\begin{enumerate}
\item
If $c_{21}=0$, then $c_{11}=0$. By Eq.~(\mref{eq:ncs355}) we have $c_{22}=0$. So $c_{12}=0$. Then we obtain the solutions
$$N_{5,13_1}=\left(\begin{matrix}
0&0&0\\
0&0&0\\
c&b&0\\
\end{matrix}\right)\, (b, c \in\bfk\setminus \{0\}).
$$
\item
If $c_{21}\neq 0$, then  $e\neq 0$. So $c_{11}=-\frac{be}{c}$. By Eq.~(\mref{eq:ncs354}), we have $c_{11}+c_{22}=0$. Thus $c_{22}=\frac{be}{c}$. So we have $c_{12}=-\frac{b}{c}c_{22}=-\frac{b^2e}{c^2}
$.
Then we obtain the solutions $$N_{5,13_2}=\left(\begin{matrix}
-\frac{be}{c}&-\frac{b^2e}{c^2}&0\\
e&\frac{be}{c}&0\\
c&b&0\\
\end{matrix}\right)\, (b,c,e\in\bfk\setminus \{0\}).
$$

\end{enumerate}
In summary, we obtain the solutions

$$N_{5,13}=\left(\begin{matrix}
-\frac{be}{c}&-\frac{b^2e}{c^2}&0\\
e&\frac{be}{c}&0\\
c&b&0\\
\end{matrix}\right)\, (b,c\in\bfk\setminus \{0\}, e\in\bfk).
$$

\noindent
{\bf Subcase 3. $b=c=0$ and $a,d\neq 0$:} This subcase is similar to the above Subcase 2.
Denote $e=c_{21}$. Then we can obtain
the solutions

$$N_{5,14}=\left(\begin{matrix}
-\frac{ae}{d}&-\frac{a^2e}{d^2}&a\\
e&-\frac{ae}{d}&d\\
0&0&0\\
\end{matrix}\right)\, (a,d\in\bfk\setminus \{0\}, e\in\bfk).
$$

\smallskip
\noindent
{\bf Subcase 4. $b=d=0$ and $a,c\neq 0$:} This subcase is similar to the above Subcase 1.
Then we can obtain
the Rota-Baxter operators

$$N_{5,15}=\left(\begin{matrix}
-\sqrt{ac}i&0&a\\
0&0&0\\
c&0&\sqrt{ac}i\\
\end{matrix}\right) \text{ and }
N_{5,16}=\left(\begin{matrix}
\sqrt{ac}i&0&a\\
0&0&0\\
c&0&-\sqrt{ac}i\\
\end{matrix}\right)\quad (a,c\in \bfk\setminus \{0\}).
$$
}
\smallskip
\noindent
{\bf Case 4. Exactly three of $a,b,c,d$ are $0$:}
Then we divide into four subcases to consider.

\smallskip
\noindent
{\bf Subcase 1. $a=b=c=0$ and $d\neq 0$:}
Then $c_{13}=c_{32}=c_{31}=0$ and $c_{23}\neq0$.  Thus by Eq.~(\mref{eq:ncs353}), we have $c_{12}c_{23}=0$. So $c_{12}=0$ and then Eqs.~(\mref{eq:ncs351}) and ~(\mref{eq:ncs355}) give $c_{11}=c_{22}=0$. By Eq.~(\mref{eq:ncs359}), we have $c_{33}=0$. Denote $e=c_{21}$. Then we obtain the solutions
$$N_{5,17}=\left(\begin{matrix}
0&0&0\\
e&0&d\\
0&0&0\\
\end{matrix}\right)\quad (d, e\in\bfk, d\neq 0).$$

Similarly, we obtain the following solutions for the rest of the subcases.

\smallskip
\noindent
{\bf Subcase 2. $a=b=d=0$ and $c\neq 0$:} Denote $e=c_{21}$. Then we have
$$N_{5,18}=\left(\begin{matrix}
0&0&0\\
e&0&0\\
c&0&0\\
\end{matrix}\right)\quad (c, e\in\bfk, c\neq 0).$$

\smallskip
\noindent
{\bf Subcase 3. $a=c=d=0$ and $b\neq 0$:}
Denote $e=c_{12}$.  Then we have
$$N_{5,19}=\left(\begin{matrix}
0&e&0\\
0&0&0\\
0&b&0\\
\end{matrix}\right)\quad (b, e\in\bfk, b\neq 0).$$

\smallskip
\noindent
{\bf Subcase 4. $b=c=d=0$ and $a\neq 0$:}
Denote $e=c_{12}$, where $e\in\bfk$. Then we have
$$N_{5,20}=\left(\begin{matrix}
0&e&a\\
0&0&0\\
0&0&0\\
\end{matrix}\right)\quad (a, e\in\bfk, a\neq 0).$$

\smallskip
\noindent
{\bf Case 5. $a=b=c=d=0$:}
By Eq.~(\mref{eq:ncs359}), we have $c_{33}=0$.  Denote $e=c_{12}$ and $f=c_{21}$. Then by Eqs.~(\mref{eq:ncs351}) and ~(\mref{eq:ncs355}), we have $c_{11}=\pm\sqrt{ef}i$ and $c_{22}=\mp\sqrt{ef}i$.
\begin{enumerate}
\item
If $c_{12}\neq 0$ or $c_{21}\neq 0$, by Eqs.~(\mref{eq:ncs352}) and ~(\mref{eq:ncs354}), we have $c_{11}+c_{22}=0$. Then we can obtain the solutions
$$N_{5,21_1}=\left(\begin{matrix}
\sqrt{ef}i&e&0\\
f&-\sqrt{ef}i&0\\
0&0&0\\
\end{matrix}\right)\quad\text{and}\quad
N_{5,21_2}=\left(\begin{matrix}
-\sqrt{ef}i&e&0\\
f&\sqrt{ef}i&0\\
0&0&0\\
\end{matrix}\right)\quad (e\neq 0\, \text{or}\, f\neq 0)$$
\item
If $c_{12}=c_{21}=0$, then by Eqs.~(\mref{eq:ncs351}) and ~(\mref{eq:ncs355}), we have $c_{11}=c_{12}=0$.  Thus we obtain the zero solution
$N_{5,21_3}=\mathbf{0}_{3\times 3}.$

In summary, we obtain the solutions
$$N_{5,21}=\left(\begin{matrix}
\sqrt{ef}i&e&0\\
f&-\sqrt{ef}i&0\\
0&0&0\\
\end{matrix}\right)\quad\text{and}\quad
N_{5,22}=\left(\begin{matrix}
-\sqrt{ef}i&e&0\\
f&\sqrt{ef}i&0\\
0&0&0\\
\end{matrix}\right)\quad (e,f\in\bfk).$$
\end{enumerate}

\subsubsection{The proof for $\bfk[NCS(6)]$}
\mlabel{ss:ncs6}

We finally prove that the matrices for the Rota-Baxter operators on the semigroup algebra $\bfk[NCS(6)]$ are given by $N_{6,i}, 1\leq i\leq 9,$ in Table~\mref{ta:rbncs1}.
\mlabel{lem:ncs6}

Applying the Cayley table of $NCS(6)$ in Eqs.~(\mref{eq:3eijm}) and then
taking $i=j=1$ with $1\leq m\leq 3$; $i=1,j=2$ with $m=2$; $i=1,j=3$ with $m=2$; $i=2, 1\leq j\leq 3$ with $1\leq m\leq 3$;
 $i=3,j=1$ with $m=1,3$ and $i=3,j=3$ with $m=2$, we obtain
\begin{eqnarray}
&c_{11}^2+c_{13}c_{31}=0,
\mlabel{eq:ncs361}\\
&c_{12}^2+2c_{11}c_{12}+c_{12}c_{13}
+c_{13}c_{32}=0,\mlabel{eq:ncs362}\\
&c_{13}(c_{11}+c_{33})=0,
\mlabel{eq:ncs363}\\
&c_{12}(c_{11}+c_{21}+c_{22}+c_{23})
)+c_{13}c_{32}=0,\mlabel{eq:ncs365}\\
&c_{12}(c_{11}+c_{31}+c_{32}+c_{33})
+c_{13}c_{32}=0,\mlabel{eq:ncs368}\\
&c_{21}c_{13}+c_{22}c_{11}=c_{11}(c_{22}+c_{11})+c_{31}(c_{23}+c_{13}),\mlabel{eq:ncs3610}\\
&c_{12}(c_{21}+c_{22}+c_{11})+
c_{32}(c_{23}+c_{13})=0,
\mlabel{eq:ncs3611}\\
&c_{23}(c_{11}+c_{12}+c_{13})=
c_{13}(c_{21}+c_{11})+c_{12}c_{23}+
c_{33}(c_{23}+c_{13}),
\mlabel{eq:ncs3612}\\
&c_{21}(c_{21}+c_{23})=2c_{21}c_{11}
+2c_{23}c_{31},\mlabel{eq:ncs3613}\\
&c_{22}^2+2c_{21}c_{12}+2c_{23}c_{32}
=0,\mlabel{eq:ncs3614}\\
&c_{23}(c_{21}+c_{23})=2c_{21}c_{13}+
2c_{23}c_{33},\mlabel{eq:ncs3615}\\
&c_{21}(c_{31}+c_{33})+c_{22}c_{31}=
c_{11}(c_{21}+c_{31})+c_{31}(c_{22}+c_{23}
+c_{33}),\mlabel{eq:ncs3616}\\
&c_{12}(c_{21}+c_{31})+c_{32}(c_{22}+c_{23}
+c_{33})=0,\mlabel{eq:ncs3617}\\
&c_{23}(c_{31}+c_{32}+c_{33})=
c_{13}(c_{21}+c_{31})+c_{32}c_{23}+
c_{33}(c_{23}+c_{33}),\mlabel{eq:ncs3618}\\
&c_{31}(c_{11}+c_{33})=0,\mlabel{eq:ncs3619}\\
&c_{33}^2+c_{13}c_{31}=0,\mlabel{eq:ncs3621}\\
&c_{32}(c_{31}+c_{32}+2c_{33})+c_{31}
c_{12}=0.\mlabel{eq:ncs3626}
\end{eqnarray}
Note that  Eqs.~(\mref{eq:ncs361}) and ~(\mref{eq:ncs3621}) give $c_{11}^2=c_{33}^2$. Depending on whether $c_{12}=0$ or not, we divide the proof into two cases.

\smallskip
\noindent
{\bf Case 1. $c_{12}=0$:} Then by Eq.~(\mref{eq:ncs362}), we have $c_{13}c_{32}=0$. Denote $a=c_{13}$. Then there are two subcases to consider.

\smallskip
\noindent
{\bf Subcase 1. $a\neq 0$:}  Then we have $c_{32}=0$. Thus by Eq.~(\mref{eq:ncs3614}),  we have $c_{22}=0$.  Eq.~(\mref{eq:ncs363}) gives $c_{11}+c_{33}=0$. Denoting $b=c_{31}$, then by Eq.~(\mref{eq:ncs361}), we have $c_{11}=\pm\sqrt{ab}i$.
Thus $c_{33}=\mp\sqrt{ab}i$.
\begin{enumerate}
\item
Let $a+b=0$. Then $b=-a$. By Eq.~(\mref{eq:ncs361}), we have $c_{11}=a$ or $c_{11}=-a$.
\begin{enumerate}
\item
Let $c_{11}=a$.  Then $c_{33}=-a$. By Eqs.~(\mref{eq:ncs3610}) and ~(\mref{eq:ncs3612}), we get $c_{21}=c_{23}=0.$ Thus we obtain the solutions
$$N_{6,1}=\left(
\begin{matrix}
a&0&a\\
0&0&0\\
-a&0&-a\\
\end{matrix}\right)\quad (a\in \bfk, a\neq 0).$$
\item
Let $c_{11}=-a$. Then $c_{33}=a$. By Eqs.~(\mref{eq:ncs3610}) and ~(\mref{eq:ncs3612}), we get $c_{21}+c_{23}=0$. Denote $c=c_{23}$. Then $c_{21}=-c$. Thus we obtain the solutions
$$N_{6,2}=\left(\begin{matrix}
-a&0&a\\
-c&0&c\\
-a&0&a\\
\end{matrix}\right)\quad (a, b\in \bfk, a\neq0).
$$
\end{enumerate}

\item
Let $a+b\neq0$.
Then by Eq.~(\mref{eq:ncs3618}), we get $c_{23}c_{31}=c_{13}c_{21}+(c_{13}c_{31}+
c_{33}^2)$. Using Eq.~(\mref{eq:ncs3621}), we have $c_{23}c_{31}=c_{13}c_{21}$. Since $c_{13}=a\neq 0$, we have $c_{21}=\frac{b}{a}c_{23}$. From Eq.~(\mref{eq:ncs3616}), we get $c_{21}(c_{31}+c_{33}-c_{11})=c_{31}c_{23}
+c_{31}(c_{11}+c_{33})$. So by Eq.~(\mref{eq:ncs3619}), we have $c_{21}(c_{31}+c_{33}-c_{11})=c_{31}c_{23}.$
Together with $c_{21}=\frac{b}{a}c_{23},$
we have
\begin{equation}
\frac{b}{a}\,c_{23}(c_{13}+c_{11}-(c_{31}+
c_{33}))=0.
\mlabel{eq:ncs36im1}
\end{equation}
Assume $c_{13}+c_{11}=c_{31}+c_{33}$, using $c_{11}=-c_{33}$,  then we have $2c_{33}=c_{13}-c_{31}=a-b.$ So $4c_{33}^2=(a-b)^2$. Since $c_{33}^2+c_{13}c_{31}=0,$ then $-4ab=(a-b)^2$. Thus $(a+b)^2=0$, a contradiction. Thus $c_{13}+c_{11}\neq c_{31}+c_{33}$. Then we have the following two cases.
\begin{enumerate}
\item
Let $b\neq 0$.  Then by Eq.~(\mref{eq:ncs36im1}), we have $c_{23}=0$, and so $c_{21}=0$. Since $c_{11}=\pm\sqrt{ab}i$ and $c_{11}+c_{33}=0$, we have $c_{33}=\mp \sqrt{ab}$. Thus we obtain the solutions
$$N_{6,3}=\left(
\begin{matrix}
\sqrt{ab}i& 0&a\\
0&0&0\\
b&0&-\sqrt{ab}i\\
\end{matrix}\right)\text{ and }
N_{6,4}=\left(
\begin{matrix}
-\sqrt{ab}i& 0&a\\
0&0&0\\
b&0&\sqrt{ab}i\\
\end{matrix}\right)\quad (a, b\in \bfk, a ,b\neq 0).$$
\item
Let $b=0$. Then $c_{31}=0$. By Eq.~(\mref{eq:ncs361}), we have $c_{11}=0$. Since $c_{11}^2=c_{33}^2$, we have $c_{33}=0$. Then Eq.~(\mref{eq:ncs3618}) gives $c_{21}=0$. So Eq.~(\mref{eq:ncs3615}) gives $c_{23}=0$. Thus we obtain the solutions
$$N_{6,5}=\left(
\begin{matrix}
0&0&a\\
0&0&0\\
0&0&0\\
\end{matrix}\right)\quad (a\in \bfk, a\neq 0).$$
\end{enumerate}
\end{enumerate}

\smallskip
\noindent
{\bf Subcase 2 $a=0$.} Then by Eq.~(\mref{eq:ncs361}), we have $c_{11}=0$, and so $c_{33}=0$. Eqs.~(\mref{eq:ncs3610}) and~(\mref{eq:ncs3611}) give $c_{31}c_{23}=0$ and $c_{32}c_{23}=0$, respectively. So by Eq.~(\mref{eq:ncs3614}), we have $c_{22}=0$.
\begin{enumerate}
\item
Let $c_{31}=0$.
Then by Eq.~(\mref{eq:ncs3626}), we have $c_{32}=0$. By Eqs.~(\mref{eq:ncs3613}) and ~(\mref{eq:ncs3615}), we get $c_{21}(c_{21}+c_{23}=0$ and $c_{23}(c_{21}+c_{23})=0$. Thus $c_{21}+c_{23}=0$. Take $d=c_{21}$.  Then $c_{23}=-d$. Thus we obtain the solutions
$$N_{6,6}=\left(\begin{matrix}
0&0&0\\
d&0&-d\\
0&0&0\\
\end{matrix}\right)\quad (d\in \bfk).$$
\item
Let $c_{31}\neq 0$.  Then $c_{23}=0$. By Eq.~(\mref{eq:ncs3613}), we have $c_{21}=0$. Then Eq.~(\mref{eq:ncs3626}) gives $c_{32}(c_{31}+c_{32})=0$.  Denote $d=c_{31}$. Then $0\neq d\in\bfk$. Then we get $c_{32}=0$ or $c_{32}=-d$. Thus we get
the solutions
$$N_{6,7}=\left(\begin{matrix}
0&0&0\\
0&0&0\\
d&0&0\\
\end{matrix}\right)\quad \text{and} \quad N_{6,8}=\left(\begin{matrix}
0&0&0\\
0&0&0\\
d&-d&0\\
\end{matrix}\right)\quad (d\in \bfk, d\neq 0).$$
\end{enumerate}

\smallskip
\noindent
{\bf Case 2 $c_{12}\neq 0$.} Denote $a=c_{12}$ and $b=c_{32}$. Assume $c_{13}=0$. Then by Eq.~(\mref{eq:ncs361}) we have $c_{11}=0$. So by Eq.~(\mref{eq:ncs362}), we have $c_{12}=0$, a contradiction. Thus $c_{13}\neq 0$. By Eqs.~(\mref{eq:ncs362}), (\mref{eq:ncs365}) and ~(\mref{eq:ncs368}), we have $c_{11}+c_{12}+c_{13}=c_{21}+c_{22}+c_{23}
=c_{31}+c_{32}+c_{33}$. Then by adding Eqs.~(\mref{eq:ncs361})-(\mref{eq:ncs363}), we have $(c_{11}+c_{12}+c_{13})^2=0$. So $c_{11}+c_{12}+c_{13}=0$. Thus Eq.~(\mref{eq:ncs362}) gives $c_{11}c_{12}+c_{13}c_{32}=0$. Since $a=c_{12}\neq 0$, we have $c_{11}=-\frac{b}{a} c_{13}$. Since $c_{11}+c_{12}+c_{13}=0$,  we have $c_{11}+c_{13}=-a$. Then $c_{13}-\frac{b}{a}c_{13}=-a$, and so
$\frac{a-b}{a}c_{13}=-a$.
Assume $a=b$. Then $c_{11}+c_{13}=0$. This means that $c_{12}=0$, a contradiction. Thus $a-b\neq0$. So we have $c_{13}=\frac{a^2}{b-a}$, and then $c_{11}=\frac{ab}{a-b}$. Since $c_{13}\neq0$, Eq.(\mref{eq:ncs363}) gives $c_{11}+c_{33}=0$. Then we get $c_{33}=\frac{ab}{b-a}$. Since $c_{31}+c_{32}+c_{33}=0$, we have $c_{31}=-c_{32}-c_{33}=\frac{b^2}{a-b}$.
 By Eq.~(\mref{eq:ncs3617}) and $c_{22}+c_{23}=-c_{21}$, we get $(c_{12}-c_{32})c_{21}=-(c_{12}c_{31}+c_{
32}c_{33})=0$. Since¡¡$c_{12}-c_{32}=a-b\neq 0$, we have $c_{21}=0$.
There are two subcases.

\smallskip
\noindent
{\bf Subcase 1. $b\neq 0$:} Then $c_{31}=\frac{b^2}{a-b}\neq 0$.
By Eq.~(\mref{eq:ncs3613}), we have $c_{23}c_{31}=0$. So $c_{23}=0$. Thus by Eq.~(\mref{eq:ncs3614}), we have $c_{22}=0$.
Then we obtain the solutions
$$N_{6,9_1}=
\left(
\begin{matrix}
\frac{ab}{a-b}&a& \frac{a^2}{b-a}\\
0&0&0\\
\frac{b^2}{a-b}&b&\frac{ab}{b-a}\\
\end{matrix}\right)\quad (a,b\in\bfk\setminus\{0\}).$$

\smallskip
\noindent
{\bf Subcases 2. $b=0$:}
Then $c_{32}=0$.  From above discussion, we have $c_{11}=c_{31}=c_{33}=0$ and $c_{13}=-a$.
By $c_{21}=0$, Eq.~(\mref{eq:ncs3614}) gives $c_{22}=0$. Since $c_{21}+c_{22}+c_{23}=0$, we have $c_{23}=0$.
\delete{
By Eq.~(\mref{eq:ncs3623}), we have $c_{31}c_{12}=0$. Since $c_{12}\neq 0$, we have $c_{31}=0$. Thus Eq.~(\mref{eq:ncs361}) and ~(\mref{eq:ncs3621}) give $c_{11}=0$ and $c_{33}=0$, respectively. Since $c_{11}+c_{12}+c_{13}=0$, we have $c_{13}=-a$.}
Then we get the solutions
$$N_{6,9_2}=
\left(
\begin{matrix}
0&a& -a\\
0&0&0\\
0&0&0\\
\end{matrix}\right)\quad (a\in \bfk, a\neq 0).$$

In summary, we get the solutions
$$
N_{6,9}=\left(
\begin{matrix}
\frac{ab}{a-b}&a& \frac{a^2}{b-a}\\
0&0&0\\
\frac{b^2}{a-b}&b&\frac{ab}{b-a}\\
\end{matrix}\right)\quad (a, b\in\bfk, a\neq 0).$$
This completes the proof for $\bfk[NCS(6)]$.
\smallskip

Now the proof of Theorem~\mref{thm:rbonc3} is completed.

\section{The method of computer algebra}
\mlabel{sec:comp}

In this section, we describe the computer algebra procedure (in Mathematica) to
compute the Rota-Baxter operators on semigroup algebras of semigroup of order
$3$ in Tables~\mref{ta:s3} and~\mref{ta:ncs3}. This procedure serves both for
guiding and verifying the manual proofs of the classification theorems carried
out in the main body of the paper.

The Mathematica code and accompanying syntax definitions are given in
Figure~\mref{ta:pro} in page \pageref{ta:pro}. The function~\texttt{RBA} with four arguments creates the
equations in Eq.~(\mref{eq:neijm}) for a fixed pair of elements, which are then
instantiated by all generator pairs. For added clarity, we have also displayed
the general form of these equations for the generic $2\times2$ Cayley table
defined at the beginning. The main function for finding Rota-Baxter operators
is~\texttt{FindRBO}, which simply solves the equations created
by~\texttt{RBA}. In the generic case, this yields only the trivial operator. For
converting a given Cayley table to the structure constants~$r^m_{kl}$ used
in Eq.~(\mref{eq:neijm}), the simple function~\texttt{SGM} is used.

  \begin{figure}[h!]
    \fbox{%
      \includegraphics[width=\textwidth,height=0.9\textheight]{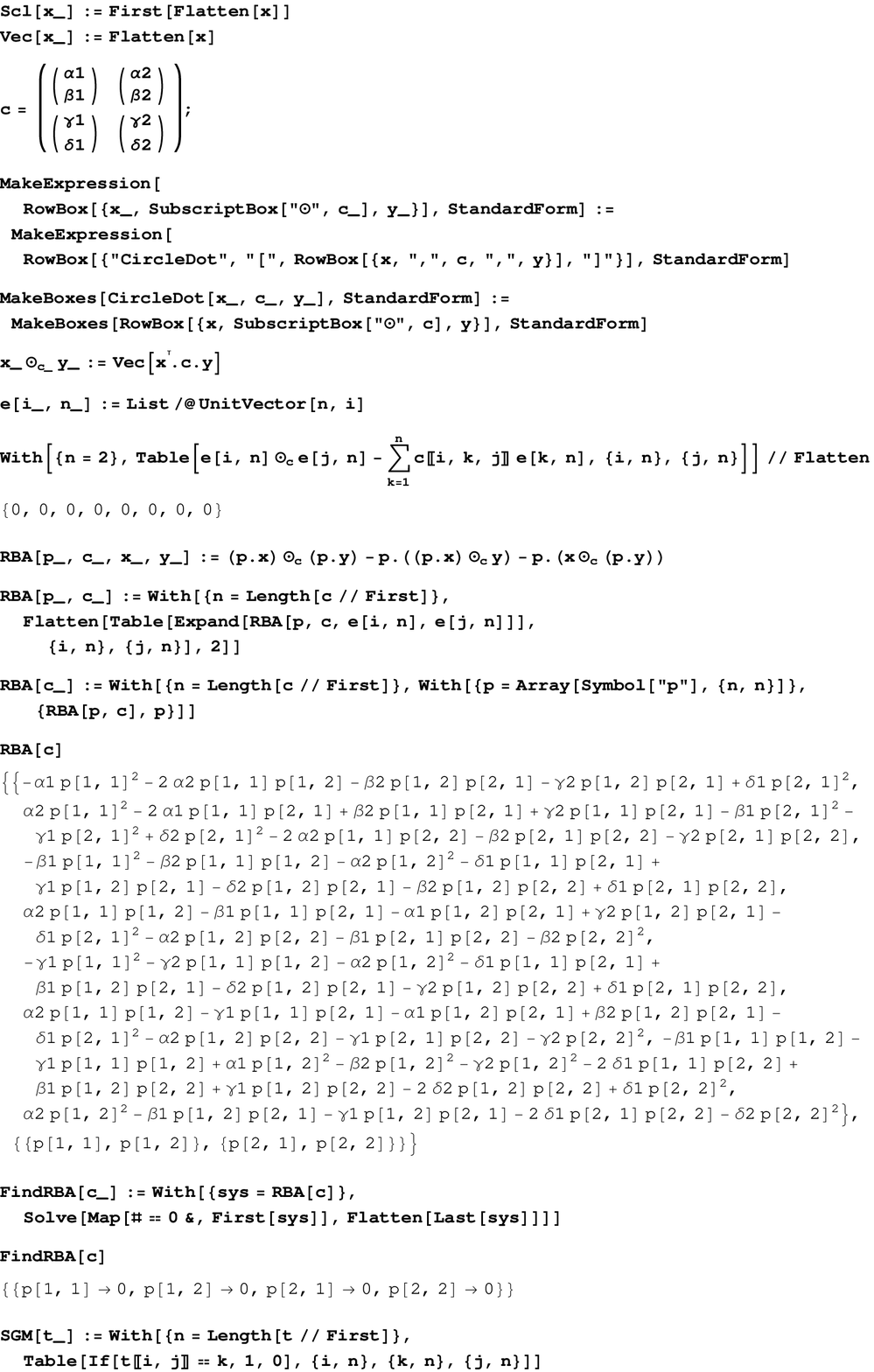}
    }
    \caption{The procedure to compute Rota-Baxter operators on semigroup
      algebras\mlabel{ta:pro}}
  \end{figure}

We illustrate these functions by considering the first semigroup of
  Table~\mref{ta:s3}, so let $t=CS(1)$. Here the underlying set $\{e_1,e_2,
  e_3\}$ of $CS(1)$ is simplified to $\{1, 2, 3\}$. The above Mathematica code
  above yields the results given in Figure~\mref{ta:pro2} in page~\pageref{ta:pro2}. In fact, the output
  gives two Rota-Baxter operators for~$CS(1)$ but the second is a special case
  of the first.  Let $p(1,1)=a, p(2,1)=b, p(1,2)=c, p(2,2)=d, p(1,3)=e,
  p(2,3)=f$, where $a,b,c,d,e,f\in\bfk$. Then $p(3,1)=-a-b, p(3,2)=-c-d$ and
  $p(3,3)=-e-f$.  So we get the matrix
$$\left(\begin{matrix}
a&c&e\\
b&d&f\\
-a-b&-c-d&-e-f\\
\end{matrix}\right)\quad (a,b,c,d,e,f\in\bfk).$$
The transpose of the above matrix is $C_{1,1}$ in Table~\mref{ta:rbcs}.

  \begin{figure}[h!]
    \fbox{%
      \includegraphics[width=\textwidth,height=0.2\textheight]{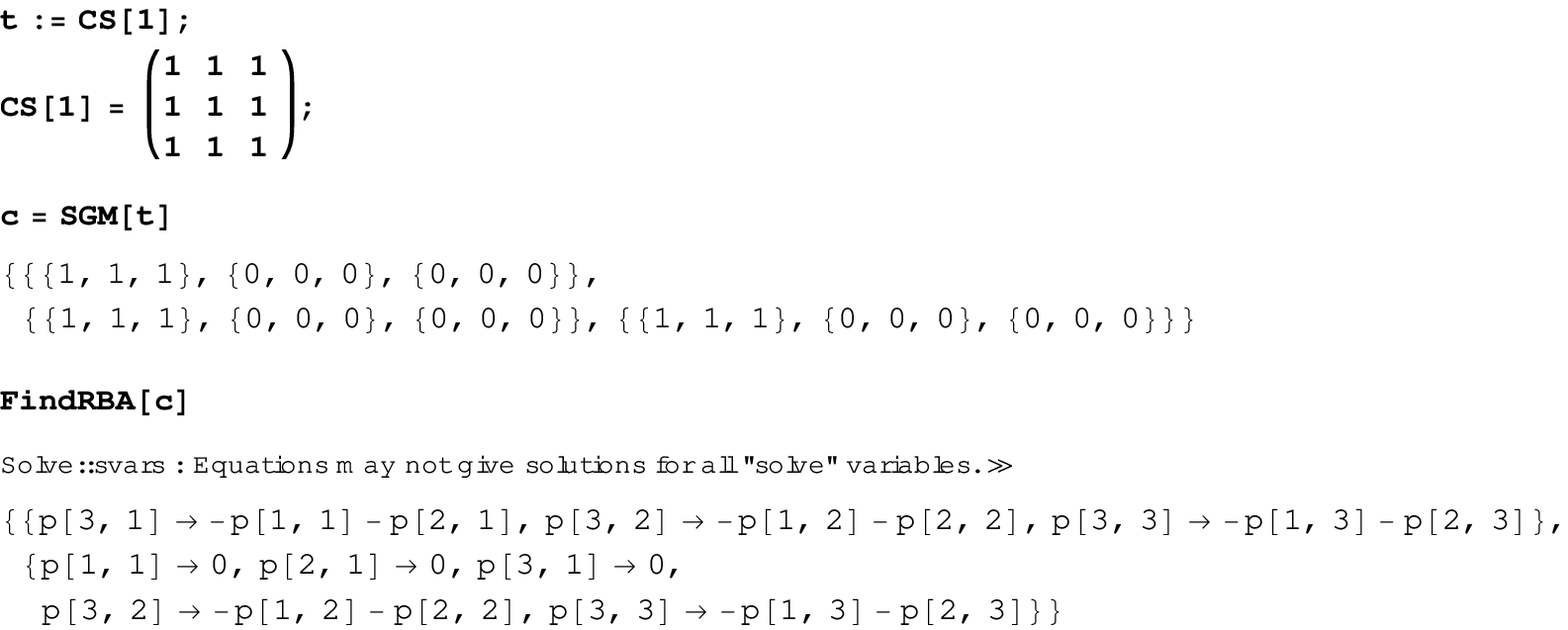}
    }
    \caption{Results of the Mathematica code on $CS(1)$\mlabel{ta:pro2}}
  \end{figure}

\noindent {\bf Acknowledgements}: This work was supported by the National
Science Foundation of US (Grant No. DMS~1001855), the
Engineering and Physical Sciences Research Council of UK (Grant No. EP/I037474/1) and the National Natural Science Foundation of China (Grant No. 11371178).

\end{document}